\documentclass[12pt, twoside, a4paper]{article} 

\usepackage[utf8]{inputenc}
\usepackage[T1]{fontenc}
\usepackage{geometry}
\usepackage{lmodern}
\usepackage{amsmath, amssymb, amsthm}
\usepackage{mathtools}
\usepackage{microtype}
\usepackage{xcolor}
\usepackage{hyperref}
\usepackage{titlesec}
\usepackage{enumitem}
\usepackage{fancyhdr}
\usepackage{tikz}
\usetikzlibrary{arrows.meta}

\usepackage{multirow}
\usepackage{booktabs}  
\usepackage{multirow}
\usepackage{arydshln}

\usepackage{tocloft}

\mathtoolsset{showonlyrefs} 
\numberwithin{equation}{section}
\allowdisplaybreaks

\setlist[itemize]{itemsep=3pt, topsep=0pt, parsep=0pt}
\setlist[enumerate]{itemsep=3pt, topsep=0pt, parsep=0pt}

\hypersetup{
	colorlinks=true,
	allcolors=blue, 
}

\theoremstyle{plain}
\newtheorem{theorem}{Theorem}[section]
\newtheorem{lemma}[theorem]{Lemma}
\newtheorem{proposition}[theorem]{Proposition}
\newtheorem{corollary}[theorem]{Corollary}
\newtheorem{condition}[theorem]{Condition}  

\theoremstyle{definition}
\newtheorem{definition}[theorem]{Definition}

\newtheorem{remark}[theorem]{Remark}
\newtheorem{notation}[theorem]{Notation}

\DeclareMathOperator{\Esp}{E}
\DeclareMathOperator{\Prob}{P}
\DeclareMathOperator{\Qrob}{Q}
\DeclareMathOperator{\IR}{\mathbb{R}}
\DeclareMathOperator{\IN}{\mathbb{N}}
\DeclareMathOperator{\bF}{\mathcal{F}}

\DeclareMathOperator{\erfc}{erfc}
\DeclareMathOperator{\sgn}{sgn}
\DeclareMathOperator{\lc}{l.c.}
\DeclareMathOperator{\rc}{r.c.}

\DeclareMathOperator{\Lip}{Lip.}
\DeclareMathOperator{\dist}{dist}
\DeclareMathOperator{\Mills}{Mills}
\DeclareMathOperator{\Int}{Int}

\newcommand{\mfk}[1]{\mathfrak}
\newcommand{\mc}[0]{\mathcal}
\newcommand{\wh}[0]{\widehat}
\newcommand{\wt}[0]{\widetilde}
\newcommand{\un}[0]{u_n}
\newcommand{\params}[2]{({#1},{#2})}
\newcommand{\xlra}[2]{ \xrightarrow[#2]{#1} }
\newcommand{\ra}[0]{ \rightarrow }

\newcommand{\convergence}[1]{ \xlongrightarrow[n\ra \infty]{#1} }
\newcommand{\ucp}[0]{\text{\textup{ucp}}}
\newcommand{\law}[0]{\text{\textup{law}}}
\newcommand{\rd}{\mathrm{d}}
\newcommand{\vd}{\,\mathrm{d}}
\newcommand{\process}[1]{(#1)_{t\ge 0}}
\newcommand{\indic}[1]{\mathbf{1}_{#1}}
\newcommand{\indicB}[1]{\mathbf{1}_{\{{#1}\}}}
\newcommand{\loct}[3]{L^{#2}_{#3}(#1)}

\newcommand{\rloct}[3]{{\ell}^{#2}_{#3}(#1)}

\newcommand{\occt}[3]{{\mc O}^{#2}_{#3}(#1)}
\newcommand{\hfprocess}[3]{{#1}^{#2}_{\frac{{#3}}{n}}}

\newcommand{\braces}[1]{ \left({#1}\right) } 
\newcommand{\sqbraces}[1]{ \left[{#1}\right] } 
\newcommand{\cubraces}[1]{ \left\{{#1}\right\}} 
\newcommand{\xnorm}[2]{ \left\|{#1}\right\|_{#2}} 
\newcommand{\abs}[1]{ \left|{#1}\right|} 
 
\newcommand{\qv}[1]{ \left\langle {#1} \right\rangle }

\newenvironment{enuroman}{\begin{enumerate}[label=(\roman*)]}{\end{enumerate}}
\newenvironment{xenumerate}[1]{\begin{enumerate}[label= {\bfseries {#1}\arabic*.}, ref= { \bfseries  {#1}\arabic*}]}{\end{enumerate}}
\newcommand{\hemail}[1]{\href{mailto:{#1}}{\textit{#1}}}  

\newenvironment{cdate}{
	\noindent 
	\textbf{\small This version:}
	\small
}{\par}

\newenvironment{keys}{
	\noindent
	\textbf{\small Keywords:}
	\small
}{\par}

\newenvironment{msc}{
	\noindent
	\textbf{\small MSC 2020 codes:}
	\small
}{\par}

\usepackage{authblk}

\title{Sticky-Threshold Diffusions, Local Time \\ Approximation and Parameter Estimation}

\author[1]{Alexis Anagnostakis\thanks{\hemail{alexis.anagnostakis@yandex.com}, \hemail{aanagnostakis@cmm.uchile.cl}}}
\author[2]{Sara Mazzonetto\thanks{\hemail{sara.mazzonetto@univ-lorraine.fr}}}
\affil[1]{Centro de Modelamiento Matemático (CNRS IRL2807), Universidad de Chile, Santiago, Chile}
\affil[2]{Universit\'e de Lorraine, CNRS, IECL, F-54000 Nancy, France}

\date{ }

\newcommand{\shorttitle}{Local time approximation for sticky diffusions}
\newcommand{\shortauthors}{A. Anagnostakis and S. Mazzonetto}

\begin{document}
	
	\maketitle

	\begin{abstract}
We study a class of high-frequency path functionals for one-dimensional diffusions with singular thresholds or boundaries, allowing for skewness, discontinuities in the diffusion coefficient, and stickiness or sticky reflection. These functionals, originally developed for local-time approximation in non-singular diffusions, are constructed from a test function and a diverging normalizing sequence. We establish convergence to local time, regardless of the nature of the singular threshold features, thereby extending several recent results on specific cases. Notably, our framework allows for any normalizing diverging sequence that is $o(n) $, where $n $ is the observation frequency, and thresholds at which several singular behaviors occur simultaneously. Combining the local-time approximation with occupation-time approximations, we construct consistent estimators of the stickiness and skewness parameters that remain valid in the presence of jump-discontinuities in the diffusion coefficient; this solves the estimation problem under the simultaneous presence of all these singular features.
	\end{abstract}

    \begin{quote}
    \begin{cdate}
		\today
	\end{cdate}
	
	\begin{msc}
		62F12; 60J55; 60J60.
	\end{msc}
	
	\begin{keys}
		Sticky Brownian motion, Skew Brownian motion, reflection, Oscillating Brownian motion, parameter estimation, high frequency statistics, local time approximation, occupation time approximation.
	\end{keys}
    \end{quote}
	
	\bigskip

	
	\tableofcontents
	
	\section{Introduction}
	\label{sec_introduction}

	Consider a diffusion $X$ that solves the one-dimensional time-homogeneous It\^o stochastic differential equation (SDE) $\rd X_t = b(X_t) \vd t + \sigma(X_t) \vd B_t $ with smooth coefficients $ b$, $\sigma $, and driven by the standard Brownian motion $B $. 
	It is known that the local time of $X$ at some threshold $\zeta$ can be consistently approximated as $n\to \infty $ by the empirical process
	\begin{equation}
		\label{eq_intro_loct_approximation}
		t \mapsto \frac{u_n}{n} \sum_{i=1}^{[nt]} g(u_n (\hfprocess{X}{}{i-1} - \zeta)),
	\end{equation}
	where $g$ is a bounded integrable function normalized to $\int g(x) dx = 1$, $(u_n)_n$ is a diverging sequence such that $u_n = o(n) $, $[nt]$ is the integer part of $nt $, and $X_0, X_{1/n}, \dots, X_{i/n},\dots$ are discrete observations of the process, see~\cite[Theorem 1.1]{Jac98}.  
	Such approximations are pivotal for estimating local characteristics of a diffusion;
	by diffusion we mean a continuous strong Markov process that is not necessarily solution to a classical SDE. 
	Examples are the estimation of the diffusion coefficient, when regular, \cite{FlorensZmirou1993}, of the stickiness parameter \cite{Anagnostakis2022}, of the skewness parameter~\cite{Lejay23,Lejay2019,Maz19}, and of jump sizes in the diffusion coefficient~\cite{Maz19, Robert2023}.
	
	Diffusions can exhibit a wide variety of threshold-behaviors, some of which cannot be replicated by classical SDE solutions. 
	In this paper we focus on processes that combine one or more of the following features: A threshold is \emph{sticky} if the diffusion spends a positive amount of time there, with the duration governed by a \emph{stickiness parameter}; it is \emph{skew} if it induces asymmetric reflection, characterized by a \emph{skewness parameter}; and it is \emph{oscillating} if it admits a finite jump-discontinuity in the diffusion coefficient, leading to regime-switching dynamics. 
	Such stochastic motions arise in quantum mechanics, describing particle behavior near emission sources (e.g., \cite{davies1994browniansticky}),
	particle dynamics in layered media (e.g., \cite{ramirez2013, Karatzas2016}),
	colloids (e.g., \cite{BouRabee2020,Stell91}), 
	financial modeling of stock price movements (e.g., \cite{Bass, criens2022separating,Pigato2019,Robert2023}),
	and queuing systems (e.g., \cite{Harrison1981storage}).
	
	Elementary examples include sticky Brownian motion~\cite{EngPes}, skew Brownian motion~\cite{Harrison1981skew}, oscillating Brownian motion \cite{Keilson1978}, and reflected Brownian motion.  
	Processes can also combine several such features like the oscillating-sticky Brownian motion in \cite{touhami2023oscillating}, where the threshold is both oscillating and sticky, the skew-sticky Brownian motion ({\bf SkS-BM}), studied in \cite{Touhami2021,zhang2021onproperties}, where the threshold is both skew and sticky, and the skew-oscillating-sticky Brownian motion ({\bf SOS-BM}) defined in Section~\ref{ssec_loctime}, where the threshold may be simultaneously skew, oscillating, and sticky ({\bf SOS}) with sticky reflection included as a limiting skew-sticky case. 
	The last two processes play a fundamental role in this work.
	
	\bigskip 
	
	In this paper, we study the asymptotic behavior of the statistic~\eqref{eq_intro_loct_approximation} 
	involving processes that combine the SOS features at a unique threshold, say at $0$, and that solve SDE systems of the form 
	\begin{equation}
		\label{eq_SkSID_pathwise_char}
		\begin{dcases}
			\vd X_t = \nu_t \indicB{X_t \not = 0}  \vd t + \sigma(X_t) \indicB{X_t \not = 0}  \vd B_t
			+ \beta \vd \loct{X}{0}{t},\\
			\indicB{X_t = 0}\vd t = \rho \vd \loct{X}{0}{t},
		\end{dcases}
	\end{equation}
	where $B$ is a Brownian motion, the process $(t,y) \mapsto \loct{X}{y}{t} $ is the symmetric local time field of $X$ (definition recalled in Section~\ref{sec_notation}), 
	$\nu $ is a measurable process,
	$\sigma$ may have different one-sided limits at $0$, and $\rho \ge 0 $ and $\beta \in [-1,1] $ are the parameters that determine the sticky and skew features. 
	Intuitively, the first equation governs the dynamics away from the origin, while the second relates occupation time and local time at $0$. The latter induces a time-delay at the threshold $0$ which encodes the sticky behavior.
	
	Our result extends several recent works on singular-threshold processes in two directions: it allows multiple types of singular behavior to occur simultaneously, and it covers a broader class of normalizing sequences. Existing results for skew or oscillating thresholds are restricted to $u_n=\sqrt n$, while the sticky case has previously been treated only for $u_n=o(\sqrt n)$. 
	It is important to cover all regimes for $(u_n)_n$, both to obtain a robust consistency result and to prepare the study of optimal convergence rates.
	For example, in the smooth SDE diffusion case (without stickiness, skew, or jump discontinuities of $\sigma $) it is known that \( u_n = \sqrt{n} \) is the \emph{critical regime} that achieves the optimal convergence rate for local time approximations as~\eqref{eq_intro_loct_approximation} (c.f.~\cite[Theorem~2.1]{Jac98}). 
	Accordingly, sequences satisfying $u_n=o(\sqrt n)$ form the \emph{subcritical regime}, $u_n=\sqrt n$ is the \emph{critical regime}, and the remaining sequences satisfying $u_n=o(n)$ form the \emph{supercritical regime}.
	
	Let us briefly review the results that directly concern the one-point statistic~\eqref{eq_intro_loct_approximation}. For regular diffusions, the seminal work~\cite{Jac98} establishes convergence, the corresponding rates of convergence, and the limiting law for every diverging sequence $(u_n)_n$ satisfying $u_n=o(n)$. 
	For skew and/or oscillating Brownian motion, in the regime $u_n=\sqrt n$, \cite{Lejay2019} proves consistency, while \cite{Maz19,Robert2023} establish the convergence rates and limiting laws. 
	For sticky Brownian motion, \cite{Anagnostakis2022} proves consistency for diverging sequences
	$u_n=o(\sqrt n)$, under the additional assumption that the test
	function vanishes on a neighborhood of the threshold.
	
	Thus, before the present work, the sticky case was not covered for the critical regime $u_n=\sqrt n$ or for larger sequences satisfying $u_n=o(n)$, while the skew and oscillating cases were not covered for general diverging sequences $u_n=o(n)$. 

    \medskip
	
	Our contributions are therefore the following. 
	We establish local-time approximations at a threshold that may be simultaneously skew, oscillating, and sticky, with sticky reflection included as a limiting skew-sticky case. 
	We prove consistency for every diverging sequence $(u_n)_n$ satisfying $u_n=o(n)$. For sticky thresholds, this covers the critical and supercritical regimes that were not treated in~\cite{Anagnostakis2022}. For skew and oscillating non-sticky thresholds, it extends the one-point approximation from the regime $u_n=\sqrt n$ considered in~\cite{Maz19} (and, for a skew threshold, in~\cite{Lejay2019}) to all diverging sequences $u_n=o(n)$. 
	We also weaken the condition imposed on the test function in the sticky case: instead of requiring $g$ to vanish in a neighborhood of the threshold, we only require $g(0)=0$. 
	
	Our contributions with respect to the existing literature on the local time approximation problem at a singular threshold are summarized in Table~\ref{tab: contributions}.
	
	\begin{table}
		\raggedright
		\begin{tabular}{r c | c c c}
			& \textbf{Paper}
			& \textbf{Type of singularity}
			& \textbf{Sequences}
			& \textbf{Statistic(s)} \\
			\hline
			
			& Jacod~\cite{Jac98}
			& regular diffusion
			& $o(n)$
			& \eqref{eq_intro_loct_approximation} as special case \\
			
			\cmidrule{2-5}
			
			& Mazzonetto~\cite{Maz19}
			& skew-oscillating
			& $\sqrt n$
			&  \eqref{eq_intro_loct_approximation} as special case \\
			
			\cmidrule{1-5}
			
			\multirow{4}{*}{
				\rotatebox[origin=c]{90}{
					\parbox{2cm}{
						\centering
						\textbf{\footnotesize On sticky thresholds}
					}
				}
			}
			
			& \multirow{2}{*}{
				Anagnostakis~\cite{Anagnostakis2022}
			}
			& \multirow{2}{*}{sticky}
			& \multirow{2}{*}{$o(\sqrt n)$}
			& \eqref{eq_intro_loct_approximation} with $g$ \\
			& & & & vanishing around $0$ \\
			
			\cmidrule{2-5}
			
			& \multirow{2}{*}{Present paper}
			& skew-oscillating-sticky
			& \multirow{2}{*}{$o(n)$}
			& \eqref{eq_intro_loct_approximation} with \\
			& & or sticky-reflected
			& & $g(0)=0$ \\
		\end{tabular}		
		\caption{Results directly covering the one-point local-time approximation statistic \eqref{eq_intro_loct_approximation}. Regular diffusions are special cases of sticky and skew-oscillating diffusions, and that sticky diffusions are themselves special cases of skew-oscillating-sticky diffusions.}
		\label{tab: contributions}
	\end{table}

The proof is organized around SkS-BM, which serves as the canonical process reproducing the skew and sticky behavior of the threshold. We first prove the critical-regime result. A rescaling argument based on this result then yields convergence whenever $\log n = o(u_n)$. 

In the sticky-BM case, this also provides an alternative proof for a substantial part of the subcritical regime considered in~\cite{Anagnostakis2022} (whose argument itself builds on~\cite{Jac98}), while at the same time covering the critical and supercritical regimes. A complementary argument, following the approach of~\cite{Jac98,Anagnostakis2022}, covers the remaining subcritical cases $u_n = o(\log n)$. These arguments rely on new estimates for the transition kernel of SkS-BM.
	
	We then transfer the approximation to sticky-reflected Brownian motion, to SOS-BM through a piecewise-linear transformation, and finally to general SOS-diffusions through a Lamperti transformation and local equivalence of laws. 
	The dependencies between these reductions are displayed in Figure~\ref{fig:dependence}.
	
	We then use the local-time approximation to construct consistent estimators of the stickiness and skewness parameters. The stickiness estimator combines the occupation-local time relation in~\eqref{eq_SkSID_pathwise_char} with a discrete occupation-time approximation. For SOS-BM, we additionally construct consistent estimators of the two volatility levels; combined with the preceding estimators, these yield joint consistent estimation of the skewness, stickiness, and volatility parameters.
	
	We do not determine convergence rates or optimal asymptotic distributions in the sticky setting. Related crossing statistics, open questions concerning convergence rates, possible one-dimensional extensions, and connections with multidimensional sticky diffusions are discussed in Section~\ref{sec_comments_extensions}.

	\bigskip
	The paper is organized as follows. 
	Section~\ref{sec_notation} introduces the notions of convergence, local-time conventions, and notation used throughout the paper.
	Section~\ref{sec_results} introduces the processes of interest and states the main results on local-time approximation and parameter estimation.     Section~\ref{sec_comments_extensions} discusses related crossing statistics, open questions concerning convergence rates, and possible extensions of the main results. 
	Section~\ref{sec_preliminary} establishes preliminary results: bounds on the semigroup of SkS-BM and a preliminary local time approximation based on the Tanaka formula. 
	Lengthy proofs for this section are deferred to Appendix~\ref{app_proofsSemigroup}.
	In Section~\ref{sec_SkSBM}, we prove the main result for SkS-BM by differentiating between three cases that together cover the subcritical, critical, and supercritical regimes. 
	Section~\ref{sec_reduction} proves the remaining statements on local time approximation: for sticky-reflected BM, SOS-BM, and SOS-diffusion results by reduction to the SkS-BM case.
	Section~\ref{sec_estimation} proves all results related to parameter estimation.
    
	The appendices other than Appendix~\ref{app_proofsSemigroup} are organized as follows: 
	Appendix~\ref{app_singIto} derives explicit symmetric stochastic calculus versions of the It\^o--Tanaka formula and Girsanov theorem for SOS-diffusions.
	Appendix~\ref{app_occupation} proves an approximation of the occupation time,  used in Section~\ref{sec_estimation} to establish consistent parameter estimations.
	
	\section{Notions of convergence, notations, conventions}
	\label{sec_notation}

	Before stating the main results of this paper we introduce some preliminary notions, conventions, and notations that we adopt throughout this paper.
	
	First, we introduce the two notions of convergence in terms of which our main results are expressed: the locally uniform in time convergence, in probability, and the conditional convergence, in probability. 
	Our local time approximation results hold in locally uniform convergence in time, in probability, and our estimators converge in probability, conditionally on the event that the threshold of interest is reached.
	Next, we introduce the two notions of local time that we use: the symmetric and right local times.
    We primarily state our results in terms of symmetric local time; the right local time is used only in proofs.
	Last, we state our convention on filtrations and introduce some notations.
	
	\subsection{Convergence}
	
	Let $(A^{n})_{n\ge0} $ be a sequence of processes defined on the probability space $(\Omega, \mathcal F, \Prob) $.
	
	\begin{definition}
		[ucp convergence]
		We say that $(A^{n})_{n\ge1} $
		\emph{converges locally uniformly in time, in probability} or \ucp{} to $A^{0}$ if 
		for every $t>0 $, we have that
		\begin{equation}
			\sup_{s\le t} \abs{A^{n}_s - A^{0}_s} \xlra{\Prob}{n\ra \infty} 0.
		\end{equation}
		We denote this convergence with $A^{n} \xlra{\Prob\text{-ucp}}{n\ra \infty} A^{0}  $. 
	\end{definition}
	
	\begin{lemma}[see \cite{JacPro}, \S 2.2.3] 
		\label{lem_ucp_convergence_condition}
		If $A^n$ and $A^{0}$ have increasing paths and $A^{0}$ is continuous, then 
		\begin{equation}
			A^{n}_{t} \xlra{\Prob}{n\ra \infty} A^{0}_t,\;  \forall\, t \in D, \text{ with D dense in } \IR^+
			\; \Rightarrow \;
			A^{n} \xlra{\Prob\text{-ucp}}{n\ra \infty} A^{0}.
		\end{equation}
	\end{lemma}
	
	For the next definition assume $(Z_n)_{n\ge 0} $ is a sequence of random variables defined on a probability space $(\Omega, \mc F, \Prob) $ and $\mc L $ an event ($\mc L \in \mc F $) of positive probability ($\Prob(\mc L)>0 $).
	
	\begin{definition}
		[conditional convergence]
		We say that a sequence of random variables $(Z_n)_{n\ge 1} $ converges to $Z_0$ in probability, conditionally on $\mc L $,
		if $Z_n \to Z $ in $\Prob^{\mc L}$-probability.
		We denote this convergence with
		\begin{equation}
			Z_n \xlra{\Prob^{\mc L}}{n\ra \infty} Z_0.
		\end{equation}
	\end{definition}
	
	The next result is direct consequence of Bayes' rule.
	
	\begin{lemma}
		The sequence $(Z_n)_{n\ge 1}$ converges to $Z_0 $ in $\Prob^{\mc L}$-probability if and only if $(Z_n)_{n\ge 1} $ converges to $Z_0 $ in probability on $\mc L$, i.e.,
		for all $\varepsilon>0 $,
		\begin{equation}
			\lim_{n\to \infty} \Prob \braces{|Z_n - Z_0|>\varepsilon ;\,\mc L} = 0.
		\end{equation}
	\end{lemma}
	
	\subsection{
    Local time notions}
	
	Let $X$ be a continuous real semimartingale. 
	We denote by $(t,y) \mapsto \rloct{X}{y}{t} $ the \emph{right local time field} of $X$ defined in \cite[Theorem~VI.1.2]{RevYor} via the Tanaka formula: for every $y\in \IR $, $ \rloct{X}{y}{}  $ is the unique (up to indistinguishability) continuous non-decreasing process starting at $0$, that satisfies
	\begin{equation}
		\label{eq:tanaka}
		|X_t - y| = |X_0 - y| + \int_{0}^{t} \sgn_{\lc}(X_s-y) \vd X_s + \rloct{X}{y}{t} , \quad \forall\, t\ge 0,
	\end{equation}
	where $\sgn_{\lc} $ is the left-continuous sign function ($\sgn_{\lc}(x) = -1 $ for $x\le 0 $, and $\sgn_{\lc}(x) = 1 $ for $x>0 $).
	The \emph{left local time field} is the continuous non-decreasing process starting at $0$ that satisfies~\eqref{eq:tanaka}, but where $\sgn_{\lc} $ is replaced by the right-continuous sign function $\sgn_{\rc} $  (it is the function defined as $\sgn_{\rc}(x) = -1 $ for $x< 0 $, and $\sgn_{\rc}(x) = 1 $ for $x\geq 0 $; the symbol $\sgn $ is used for the symmetric sign function, which satisfies $\sgn(0)=0 $). 
	It is identified as the left-limit:
	\begin{equation}
		\rloct{X}{y-}{t} = \lim_{h \downarrow 0} \rloct{X}{y-h}{t}, \quad \forall\, (t,y).
	\end{equation}
	We denote by $(t,y) \mapsto \loct{X}{y}{t} $ the \emph{symmetric local time field} of the process $X$, which is the average of left and right local time fields: defined for all $t\ge 0 $ and $y\in \IR $ by 
	\begin{equation}
		\label{eq_txt_right_left_local_time_relation}
		\begin{aligned}
			\loct{X}{y}{t} &:= \frac{\rloct{X}{y}{t}+\rloct{X}{y-}{t}}{2}.
		\end{aligned}
	\end{equation}
	
	Right and symmetric local times are characterized by the following a.s. limits
	\begin{equation}
		\loct{X}{y}{t} = 
		\lim_{h\to 0} \frac{1}{2h}\int_{0}^{t} \indic{|X_s - y| < h} \vd \qv{X}_s, 
		\quad \rloct{X}{y}{t} = 
		\lim_{h\to 0} \frac{1}{h}\int_{0}^{t} \indic{y \le X_s < y+h} \vd \qv{X}_s,
		\quad \forall \, t\ge 0.
	\end{equation}
	
	It is not always the case that the right and symmetric local time fields are equal.
	An example is the skew BM, see e.g.~\cite[Theorem~2.1]{Salminen2019}.
	Typically, pathwise characterizations of diffusions are expressed in terms of the symmetric local time, see e.g.~\cite{Sal2017}.
	Conversely, results such as the It\^o--Tanaka formula \cite[Theorem~VI.1.5]{RevYor}, the occupation times formula \cite[Corollary~VI.1.6]{RevYor},
	and the representation of martingale diffusions as time-changed Brownian motion \cite[Theorem~V.47.1]{RogWilV2}, are typically expressed in terms of the right local time.

	\subsection{Conventions, notations}
	
	\subsubsection{On probability spaces and solutions to systems of SDEs} \label{sec:notions_existence}
	In this paper we always assume that (filtered) probability spaces satisfy the usual conditions of right-continuity and completeness.
	This means the filtered probability space $\mc P := (\Omega, \bF, \process{\bF_t},\Prob) $ satisfies $\bigcap_{s \ge t} \bF_s = \bF_t $, for all $t\ge 0 $, and that $\mc F_0 $ contains all $\Prob$-negligible sets.
	
	When a system of SDEs admits a unique (in law) weak solution $X$ it means that there exist a filtered probability space $\mc P$ and a BM $B$ on $\mc P$ 
	such that $(X,B) $ jointly solve the system. 
	When we wish to indicate that $X_0 = x $ a.s., we use the notation $(\mc P_x,\Prob_x) $ for the probability space and measure.
	To abbreviate the above, we say that: 
	\begin{equation}
		X\; \text{is the solution to the system of SDEs on}\; \mc P_x. 
	\end{equation}
	
	\subsubsection{Integration of functions with respect to specific measures.} 
	Let $m$ be a measure on $I$, a Borel subset of $\IR$. 
	For every $m $-integrable function $f $, every measurable function $h $, and
	every $\gamma \ge 0 $, let
	\begin{equation}
		\label{eq: measure_integral_notation}
		m(f):=\int_{I}f(x)\,m(\rd x), \quad
		m^{(\gamma)}(h):=\int_{I}|x|^{\gamma}|h(x)|\, m(\rd x).
	\end{equation}
	Throughout this paper, the reference measure in these expressions will always be defined on $I = \IR $ or $I = \IR_+ $. 
	
	By $\lambda $ we denote the Lebesgue measure on $\IR $. For functions $f $, $h $, and real $\gamma \ge 0 $, we denote by $\lambda(f) $, $\lambda^{(\gamma)}(h) $ the respective quantities defined in~\eqref{eq: measure_integral_notation}, but for $\lambda$ instead of $m$.

	\section{Main results}
	\label{sec_results}

	Before considering the general equation~\eqref{eq_SkSID_pathwise_char}, and specifying the assumptions on its coefficients ($\nu $, $\sigma $, $\beta $, $ \rho$), let us consider some particular cases of interest.  
	
	\subsection{SkS and SOS Brownian motions: introducing the singular behaviors}
	\begin{itemize}
		\item
		The skew BM is the (strong) solution to the equation 
		\begin{equation}
			\vd X_t =  \vd B_t + \beta \vd \loct{X}{0}{t}, \quad  t\ge 0,
		\end{equation}
		for $\beta\in [-1,1]$. 
		The term $\beta \vd \loct{X}{0}{t} $ generates a partial reflection effect at $0$, called skew, whose preferred direction depends on the sign of $\beta $, with limit cases $\beta \in \{-1,1\} $ corresponding to upper and lower reflection. When $\beta \in (-1,1)$ the process takes values on all $\mathbb R$, and not on a half-line as it would be if $|\beta|=1$. 
		The parameter $\beta$ is called {\it skewness}. 

		\item
		The oscillating BM is the (strong) solution to the equation 
		\begin{equation}
			\vd X_t =  \sigma_0(X_t) \vd B_t, \quad t\ge 0,
		\end{equation}
		where \[\sigma_0(x) := \indicB{x\le 0}\sigma_- + \indicB{x>0} \sigma_+,\] with $\sigma_-, \sigma_+>0 $.
		This leads the easiest example of oscillating threshold, that is $0$, i.e.~the discontinuity point of the diffusion coefficient. 
		Basically, when staying on $(0,+\infty)$, the process behaves a BM with volatility $\sigma_+$, when on $(-\infty,0)$ as a BM with volatility $\sigma_-$.

		\item
		The sticky BM is the (weak) solution to the system of equations
		\begin{equation}
			\begin{dcases}
				\vd X_t =  \indicB{X_t \not = 0}  \vd B_t & t\ge 0,\\
				\indicB{X_t = 0}\vd t = \rho \vd \loct{X}{0}{t} , & t\ge 0.
			\end{dcases}
		\end{equation}
		The equation $\indicB{X_t = 0}\vd t = \rho \vd \loct{X}{0}{t} $,  combined with vanishing coefficients, introduces a delay at $0$ called stickiness. This causes the process to spend a positive amount of time at $0$, and with intensity proportional to $\rho $.
		The reason we only consider weak solutions is that in presence of stickiness ($\rho>0$), 
		no strong solution exists, see~\cite{EngPes}.
		The parameter $\rho$ is called {\it stickiness}. 
	\end{itemize}
	
	These three processes illustrate the different threshold behavior we consider: skew (including reflecting), oscillating, and sticky. 
	These behaviors can coexist in an SOS threshold. For BM, this leads to SOS-BM, which satisfies the system of equations 
	\begin{equation}
		\label{eq_SOSBM_pathwise_char}
		\begin{dcases}
			\vd X_t = \sigma_0(X_t) \indicB{X_t \not = 0}  \vd B_t
			+ \beta \vd \loct{X}{0}{t},\\
			\indicB{X_t = 0}\vd t = \rho \vd \loct{X}{0}{t}.
		\end{dcases}
	\end{equation}
	This system is nothing else than~\eqref{eq_SkSID_pathwise_char} with $\nu \equiv 0 $, $\beta\in [-1,1]$ and diffusion coefficient  $\sigma_0$. 
    Weak existence and uniqueness in law hold for this system. For $\rho=0$ this follows, for instance, from the strong existence and pathwise uniqueness results in~\cite[Theorem~4.5]{BassChen2005}, while for $\rho>0$ it follows from~\cite[Theorem~3.4]{Sal2017}.
	In case $\sigma_0 \equiv 1 $, the SOS-BM is an SkS-BM. 
	In case $\beta=1 $ and $\sigma_0 \equiv 1 $, the SOS-BM is called \textit{sticky-reflected Brownian motion} and its state-space is $[0,\infty) $. From now on, when we consider SkS-BM we refer to the case $\beta\in (-1,1)$, unless explicitly specified.
	
	The fact that we can overlap the singular behaviors of the threshold is related to the fact that these processes are one-dimensional diffusions and can be expressed as transformations of a Brownian time-change (see \cite[Exercise~VII.3.18]{RevYor} and \cite[Theorem~V.47.1]{RogWilV2}). 
	In particular, if $X$ is the SOS-BM above with $\beta \in (-1,1)$, there exist a function $s_0$ and a measure $m_0$, called \emph{scale function} and \emph{speed measure} such that the process $s_0(X)$ is a time-changed Brownian motion $\process{B_{\gamma_0(t)}}$, where the time-change $\gamma_0$ is the right inverse of $t\mapsto \int_{\IR} \rloct{X}{x}{t} m_0(\rd x)$.
	The pair $(s_0,m_0) $ is given explicitly by
	\begin{equation}
		\begin{aligned}
			&s_0(x) := \frac{x}{a(x)},\quad 
			m_0(\rd x) := \frac{a(x)}{\sigma^{2}_0(x)} \vd x + \rho \delta_{0}(\rd x),\\   \text{where} \quad 
			& a(x):= 1 + \sgn(x)\beta. \label{eq_sm_SOSBM}
		\end{aligned}
	\end{equation}
	Note that $a$ takes values in $\{1-\beta,1+\beta\} $. 
	If $\beta=1$ (resp.~$\beta=-1$), the process is a sticky-reflected BM whose state space is $[0,+\infty)$ (resp.~$(-\infty,0]$), and $(s_0, m_0)$ above are supported on the state space: for instance if $\beta=1$, then
	\begin{equation}
		s_0(x) := \frac{x}{2} \indic{[0,+\infty)}, \quad 
		m_0(\rd x) := \frac{2}{\sigma^{2}_0(x)} \indic{[0,+\infty)} \vd x + \rho \delta_{0}(\rd x).				
	\end{equation}
	Clearly, in the latter case only $\sigma_+$ matters in $\sigma_0$. 
	Let us remark that the scale function transforms the diffusion into a local martingale, and the speed measure provides the time-change through the Brownian path of such local martingale.
	
	\bigskip
	
	\noindent
	\emph{Comments on the parametrization.} 
	The SDE path-space description of diffusions with singular thresholds involves several parametrization choices that vary across the literature (c.f. \cite{Anagnostakis2022, anagnostakis2025weak, EngPes, Sal2017, Touhami2021}).
	To avoid ambiguity, we adopt the following conventions.
	For the scale and speed, we use the same convention as in \cite{anagnostakis2025weak}: the Brownian motion has scale $s(x)=x$ and speed measure $m(\rd x)=\rd x$. 
	For the singular parameters, we set the stickiness parameter $\rho$ as the coefficient of the symmetric local time in the occupation--local time relation,
	and the skewness parameter $\beta$ as the coefficient of the symmetric local time in the SDE describing the motion away from $0$.
	
	\subsection{The SOS-diffusions considered}
	We are now ready to consider the system of equations~\eqref{eq_SkSID_pathwise_char} which generalizes the one of SOS-BM. 
	%
	We assume that the couple drift and diffusion coefficient $(\nu,\sigma) $ is such that~\eqref{eq_SkSID_pathwise_char} admits a unique (in law) weak solution that is, as specified in Section~\ref{sec:notions_existence}, 
	\begin{equation}
		X\; \text{is the solution to~\eqref{eq_SkSID_pathwise_char} on}\; \mc P_x. 
	\end{equation}
	We also refer to this process as an \emph{SOS-diffusion}. Sufficient conditions for existence and uniqueness of such processes are provided, for instance, in~\cite{Sal2017}. 
	
	Throughout the paper we impose the following additional assumptions on the pair $(\nu,\sigma)$.  
	\begin{condition}
		\label{cond_drift_volatility_functions}
		\begin{xenumerate}{c}
			\item \label{item_cond_sigma} 
			The diffusion coefficient $\sigma$ ensures that the SDE 
			\[
			\vd Y_t = \sigma(Y_t)\, \vd B_t + \beta \vd\loct{Y}{0}{t},
			\]
			admits a unique-in-law solution (a diffusion) that almost surely does not explode in finite time. 
			
			Its state space $J$ is either of the form $[0,b) $ (if $\beta=1 $),
			$(a,0] $ (if $\beta=-1 $), or $(a,b) $ (otherwise), with $-\infty \le a < 0 < b \le +\infty $. 
			
			The coefficient $\sigma$ is strictly positive on $J$ and is $C^1 $ on either side of $0$ in $J $ (on $J \cap (-\infty,0) $ and on $J \cap (0,\infty) $) with an at most finite jump at $0$ for both $\sigma$ and its derivative. It follows that the left and right limits exist: 
			\[
			0<\sigma(0-)<\infty,
			\qquad
			0<\sigma(0+)<\infty.
			\]
			If the state-space is an interval of the form $J := (\ell, 0] $ or $[0,r) $, then we assume $\sigma > 0 $, $\sigma\in C^1(J) $, and $\sigma(0) $ is finite.
			
			\item \label{item_cond_nu} The (law of the) solution to~\eqref{eq_SkSID_pathwise_char} is locally equivalent to the (law of the) solution to 
			\begin{equation}
				\label{eq_SkSID_equiv}
				\begin{dcases}
					\vd Y_t = \frac12{\sigma(Y_t) \sigma'(Y_t)}\indicB{Y_t \not = 0}  \vd t + \sigma(Y_t) \indicB{Y_t \not = 0}  \vd B_t
					+ \beta \vd \loct{Y}{0}{t}, & t\geq 0,\\
					\indicB{Y_t = 0} \vd t = \rho \vd \loct{Y}{0}{t}, & t\geq 0,
				\end{dcases}
			\end{equation}
			with the same $\sigma$, $\beta$, and $\rho$ as in~\eqref{eq_SkSID_pathwise_char}.
			Let us say that $X$ solves \eqref{eq_SkSID_pathwise_char} on $\mc P_x$ and $Y$ solves \eqref{eq_SkSID_equiv} on $\wt{\mc P}_x=(\wt \Omega, \wt \bF, \process{\wt \bF_t}, \wt \Prob_x)$. 
			Local equivalence here means that for every $T>0 $, 
			the law of $X|_{[0,T]} $ under $\Prob_x $ 
			is equivalent to the law of $Y|_{[0,T]} $ under $\wt \Prob_x $. 
		\end{xenumerate}
	\end{condition}
	
	The first assumption in Condition~\ref{cond_drift_volatility_functions} ensures that the driftless non-sticky analogue of~\eqref{eq_SkSID_pathwise_char} (where $\nu \equiv 0$) admits a unique (in law) weak solution and determines the state-space. 
	The second assumption allows us, after some additional work, to reduce the problem for the solution to~\eqref{eq_SkSID_pathwise_char} to one for the simpler base process SOS-BM. 
	
	Let us further comment on the second assumption, and discuss sufficient conditions for it. 
    The system~\eqref{eq_SkSID_equiv}, which appears in such assumption, satisfies uniqueness in law since, as we will see later,  its solutions can be bijectively transformed into SOS-BM.
    Fix $T>0$ and define for all $t\in[0,T]$
    \begin{equation}
	   \label{eq_girsanov_density_condition}
	   \theta_t :=
       \begin{cases}
		   0, & \text{if } X_t=0,
           \\
           \frac{\sigma'(X_t)}{2}-\frac{\nu_t}{\sigma(X_t)},
		& \text{if } X_t\neq0,
    	\end{cases}
	   \qquad
	   \mc E_t
	   :=
	   \exp\left\{
		  \int_0^t\theta_s\,\vd B_s
		-\frac12\int_0^t\theta_s^2\,\vd s
	   \right\}.
    \end{equation}
	By Lemma~\ref{lem_sticky_girsanov}, a sufficient criterion for the local equivalence in Condition~\ref{cond_drift_volatility_functions}-\ref{item_cond_nu} is that, for every $T>0$, $(\mc E_t)_{0\leq t\leq T}$ is a $\Prob_x$-martingale. 
	Indeed, let $\Prob_x^T$ denote the restriction of $\Prob_x$ to $\bF_T$, and define a probability measure $\Qrob_x^T$ on $(\Omega,\bF_T)$ by 
    \[
	   \frac{\vd\Qrob_x^T}{\vd\Prob_x^T}
	   =
	   \mc E_T.
    \]
    Lemma~\ref{lem_sticky_girsanov} shows that
	$
	\widetilde B_t
	:=
	B_t-\int_0^t\theta_s\,\vd s$, $0\leq t\leq T$, 
	is a Brownian motion under $\Qrob_x^T$. Moreover, under
	$\Qrob_x^T$, the process $X$ satisfies~\eqref{eq_SkSID_equiv} with respect to the Brownian motion $\widetilde B$. 
	Consequently, the law of $X|_{[0,T]}$ under $\Qrob_x^T$ coincides with the law of $Y|_{[0,T]}$ under $\widetilde\Prob_x$. Since $\Qrob_x^T\sim\Prob_x^T$, the two original laws are locally equivalent: the law of $X|_{[0,T]}$ under $\Prob_x$ is equivalent to the law of $Y|_{[0,T]}$ under $\widetilde\Prob_x$.
	
	A sufficient condition for the martingale property of
	$(\mc E_t)_{0\leq t\leq T}$ is Novikov's condition
	\begin{equation}
		\label{eq_novikov_condition}
		\Esp_x\left[
		\exp\left(
		\frac12\int_0^T
		\theta_s^2 \,\vd s
		\right)
		\right]
		<\infty. 
	\end{equation}
which holds, for example, whenever $\theta$ is a.s. bounded on $[0,T]$. This setting allows for path-dependent and time-inhomogeneous drifts: for instance, if $\sigma$ is uniformly elliptic and $\sigma'$ is bounded (i.e., there exist $c,C>0$ such that $c\leq\sigma\leq C$ and $|\sigma'|\leq C$), then Novikov's condition is satisfied for both $\nu_t = \sin(2\pi t)$ and $\nu_t:= \tanh\left( X_t-\frac1t\int_0^tX_s\,\vd s \right)$, $t>0$.
	
We note that the choice $\theta_t=0$ on $\{X_t=0\}$ is made for simplicity; other extensions are possible provided the stochastic exponential remains a martingale.

	\subsection{Local time approximation}
	\label{ssec_loctime}
	
	All our results hold for any function $g$ and sequence $\un$ satisfying the following assumption.
	
	\begin{condition} \label{cond:fun}
		Let $g$ be a bounded integrable function, with $g(0)=0$ whenever $\rho>0$, and let $(u_n)_n$ be a sequence that diverges slower than $n$, i.e.,  
		\begin{equation}
			\label{eq_intro_un_condition}
			\lim_{n\to \infty} \un = \infty,
			\qquad 
			\lim_{n\to \infty} \frac{\un}{n} = 0.  
		\end{equation}
	\end{condition}
	
		
		We first state the main result for the non-reflected SkS-BM.
		
		\begin{theorem}
			\label{thm_main_SkSBM} 
			Let $T$ be a function such that 
			\begin{equation}\label{eq_condition_T}
            \begin{split}
				& T\in C^1(\IR) \cap C^2(\IR\setminus\{0\}), \quad  T(0)=0, \quad T'(0)=1, \quad T''(0^\pm) \in \IR,
                \\
                & \text{and there exists } \varepsilon>0 \text{ such that } \forall\, y\in \IR\setminus \{0\}: \quad \varepsilon \le T'(y) \le 1/\varepsilon, \quad \abs{T''(y)} \le 1/\varepsilon.
			\end{split}
            \end{equation}
			Under Conditions~\ref{cond:fun}, 
			\begin{align}
				\label{eq_thm_limit_gnT_loctime_SkSBM}
                \frac{\un}{n} \sum_{i=1}^{[n\cdot]} g(\un T(\hfprocess{X}{}{i-1}))  &\xlra{\Prob_x\text{-ucp}}{n\ra \infty}
				\left( \int_{\IR} g \vd m_{0}\right) \loct{X}{0}{}
			\end{align}
				when $X$ is the SkS-BM solution to~\eqref{eq_SOSBM_pathwise_char} on $\mc P_x$ with $\sigma_0\equiv 1$ and $\beta\in(-1,1)$. 		
				\\
				Note that, in this case, $\int_{\IR} g \vd m_{0}=  \int_{-\infty}^{0} g(x) (1-\beta) \vd x+ \int_{0}^{\infty} g(x) (1+\beta) \vd x$. 
		\end{theorem}

        The proof of Theorem~\ref{thm_main_SkSBM} is given in Section~\ref{sec_SkSBM}. 
        The functions $T$ satisfying~\eqref{eq_condition_T} are Lamperti-like transforms. 
        In this paper, we also use the following notation to rewrite~\eqref{eq_thm_limit_gnT_loctime_SkSBM}:  Consider the functions $(g_n[T])_n$ defined as 
        \begin{equation}
            g_n[T](y) := g(\un T(y/\un))
        \end{equation}
        for all $y\in \IR $ and $n\in \IN $. Then,
        \begin{equation}
            \frac{\un}{n} \sum_{i=1}^{[n\cdot]} g_n[T](\un \hfprocess{X}{}{i-1})= \frac{\un}{n} \sum_{i=1}^{[n\cdot]} g(\un T(\hfprocess{X}{}{i-1})).
        \end{equation}
        This notation helps in the understanding that considering observation of $T(X)$ through a, say, test function $g$ which is independent on the sample size $n$, corresponds to considering observations of $X$ through a test function $g_n[T]$ which is dependent on $T$ and the sample size $n$. Clearly, when $T$ is the identity function $g=g_n[T]$. 
        
        \medskip
		
        Next, we state the result for sticky-reflected BM and SOS-BM. Note that both sticky-reflected BM and SkS-BM are special cases of SOS-BM. 
		
			%
		
		
		\begin{theorem}
			\label{thm_main_SOSBM}
			Theorem~\ref{thm_main_SkSBM} remains valid both when
			\begin{enumerate}
				[label=(\roman*), ref=(\roman*), leftmargin=*]
				\item \label{thm_item_Refl}
				$X$ is the sticky-reflected BM solution to~\eqref{eq_SOSBM_pathwise_char} on $\mc P_x$ with $\sigma_0\equiv 1$ and $\beta\in\{-1,1\}$.
				\\ 
				Note that, in the case $\beta=1$, $\int_{\IR} g \vd m_{0}=  2 \int_{0}^{\infty} g(x) \vd x$. 
				\item \label{thm_item_SOS}
				$X$ is the SOS-BM solution to~\eqref{eq_SOSBM_pathwise_char} on $\mc P_x$. \\
				Note that the speed measure $m_0$ is given in~\eqref{eq_sm_SOSBM}, 
			\end{enumerate}
		\end{theorem}
		
		The latter result for SOS-BM is then used to prove the following local time approximation for SOS-diffusions. 
		
		\begin{theorem}
			\label{thm_main}
			Let $X$ be the solution to~\eqref{eq_SkSID_pathwise_char} on $\mc P_x$.  
			Under Conditions~\ref{cond_drift_volatility_functions}-\ref{cond:fun}, it holds that 
			\begin{equation}
				\label{eq_main_result_SkSBM}
				\frac{\un}{n} \sum_{i=1}^{[n\cdot ]} g(\un \hfprocess{X}{}{i-1}) \xlra{\Prob_x\text{-ucp}}{n\ra \infty}
				\braces{\int_{\IR} g \vd m_0 } \loct{X}{0}{},
			\end{equation}
			where $m_0$ is the measure defined in~\eqref{eq_sm_SOSBM} with 
			\begin{equation}
				\label{eq_sigma0}
				\begin{aligned}
					\sigma_0(y) :=
					\begin{cases}
						\sigma(0+), &y>0, \\
						\sigma(0-), &y\le 0.
					\end{cases}
				\end{aligned}
			\end{equation}
		\end{theorem}

		\begin{remark}
			The condition $g(0)=0$ is specific to the sticky setting $\rho>0$.
			Indeed, the speed measure contains an atom at the threshold,
			of the form $\rho\,\delta_0$. Hence
			\[
			\int_{\mathbb R} g\,\rd m_0
			=
			\int_{\mathbb R\setminus\{0\}} g\,\rd m_0
			+ \rho g(0).
			\]
			Requiring $g(0)=0$ removes the direct contribution of the atomic
			part of the speed measure. This explains why the statistic
			retains a local-time limit of the same form as in the
			non-sticky case, despite the positive occupation time at the
			threshold.
			
			In the reflected cases, \(\sigma_0\) is only relevant on the state space. Its value on the inaccessible half-line may be chosen arbitrarily, since the corresponding part of the speed measure vanishes. 
		\end{remark}
		
		The proofs of 
		Theorem~\ref{thm_main_SOSBM} and Theorem~\ref{thm_main} are given in Section~\ref{sec_reduction}.
		The proofs of Theorems~\ref{thm_main_SOSBM} and~\ref{thm_main} use the following transformations:
		\begin{equation}
			\label{eq_T1_T2_def}
			\begin{aligned}
				T_1(x) &:= \frac{x}{\sigma_0(x)},
				& T_2(x) &:= \int_{0}^{x} \frac{\sigma_0(y)}{\sigma(y)} \vd y,
				& x& \in \IR.
			\end{aligned}
		\end{equation}
		
		The function $T_1$ transforms an SOS-BM into an SkS-BM, while the transformation $T_2$ (along with local equivalence of laws) intervenes in reducing the problem from an SOS-diffusion to an SOS-BM. 
		Indeed $T_2$ is a Lamperti transform that allows to transform into an SOS-BM some SOS-diffusions that solve equations of the kind~\eqref{eq_SkSID_equiv}. 
		
\begin{figure}[htbp]
			\centering
			\includegraphics[scale = 0.7]{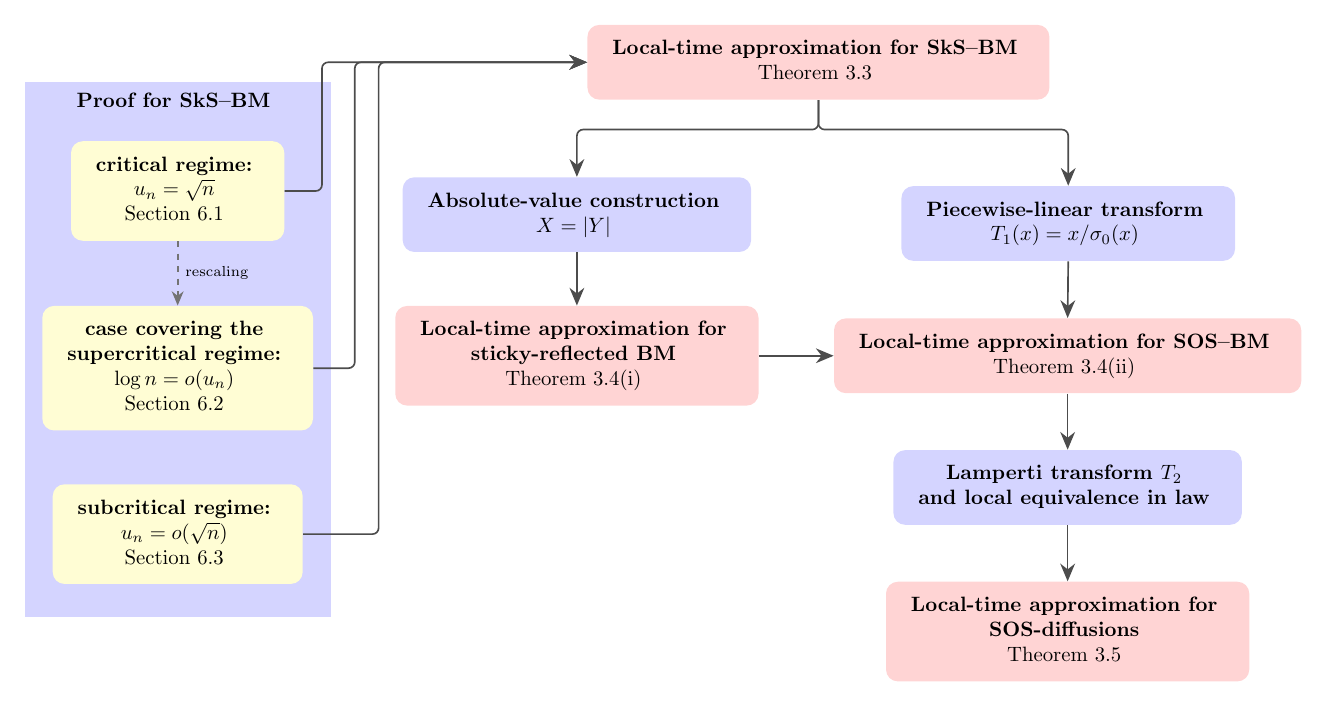}
			\caption{Proof strategy and dependencies between the main
				local-time approximation results.}
			\label{fig:dependence}
		\end{figure}

		The proof strategy and the dependencies between the local-time approximation results are summarized in Figure~\ref{fig:dependence}.
		
		\subsection{Parameter estimation}
		
		Let $t>0$ be fixed. Assume we observe an SOS-diffusion $X $ defined on $\mc P_x $ at increasing frequencies $n$: $X_0=x,X_{1/n},\ldots, X_{i/n}, \ldots, X_{[nt]/n}$.  
		We address two distinct estimation problems. 
		The first is the estimation of the parameters $(\beta, \rho)$ of an SOS-diffusion when the volatility levels $(\sigma(0-), \sigma(0+))$ are known. The second is the joint estimation of the full parameter set $(\beta, \rho, \sigma(0-), \sigma(0+))$ of an SOS-BM.
		
		Note that with respect to our local time approximation result, the second estimation problem faces an identifiability issue: The volatility gap $\sigma(0+)-\sigma(0-)$ and the skewness parameter $\beta$ are confounded in their asymptotic behavior of the statistic~\eqref{eq_intro_loct_approximation}, as noted in~\cite{Maz19}. We solve this issue in Proposition~\ref{prop_SOSBM_estimation} for the case of SOS-BM by means of additional (occupation times) statistics inspired by quadratic variation estimates for Oscillating BM in~\cite{Lejay2018}, paralleling the treatment in \cite{Maz19} (both references consider $\rho=0$). 
		For this, we need an occupation times approximation result,  Lemma~\ref{lem_sosdiff_occtime_approximation}, which extends and corrects a related result for sticky BM of~\cite{Anagnostakis2022}.
		Its proof is deferred to Appendix~\ref{app_occupation}. 
        For related general results on the discrete approximation of occupation time functionals, see also~\cite{altmeyer2021approximation}.
		
		\begin{lemma}
			\label{lem_sosdiff_occtime_approximation}
			Let $X$ be the SOS-diffusion solving~\eqref{eq_SkSID_pathwise_char} on $\mc P_x $ and let $U$ be any interval of $\IR $, singletons included. Then, the following convergence holds
			\begin{equation}
				\frac{1}{n} \sum_{i=1}^{[n\cdot]} \indic{U}(\hfprocess{X}{}{i-1})
				\xrightarrow[n\to \infty]{\Prob_x\text{-}\ucp} 
				\int_{0}^{\cdot} \indic{U}(X_s) \vd s =: \occt{X}{U}{\cdot}. 
			\end{equation}
		\end{lemma}
        We will often denote by $\occt{X}{U}{}$ the \emph{occupation time of $U$ by $X$}.  
			
        \begin{corollary} \label{cor:occ:time}
			The conclusion of Lemma~\ref{lem_sosdiff_occtime_approximation} holds for all Borel subsets $U$ of $\IR $. 
    	\end{corollary}

		\medskip
        
        For any non-trivial function $g$, let 
		$g_{>0}:= \indic{(0,\infty)}g $ and $g_{<0}:=  \indic{(-\infty,0)}g$. 
		Also, for a process $X$ and a function $g$ that satisfies $\lambda (g_{>0}), \lambda(g_{<0}) \not =0 $, let  $S_{n}^{g+}(X) $ and $S_{n}^{g-}(X)$ be statistics defined, for all $n$, as
		\begin{align}
			S_{n}^{g+}(X)&:= \frac{\un}n\frac{\sigma^{2}(0+)}{ \int_{\IR} g_{>0}(x) \vd x } \sum_{i=1}^{[nt]} g_{>0}(\un \hfprocess{X}{}{i-1}),
			\\
			S_{n}^{g-}(X)&:= \frac{\un}n \frac{\sigma^{2}(0-)}{\int_{\IR} g_{<0}(x) \vd x} \sum_{i=1}^{[nt]} g_{<0}(\un \hfprocess{X}{}{i-1} ).
		\end{align}
		Note that in order to compute these statistics, one needs to know 
		$\sigma(0-) $ and $\sigma(0+) $. 
		
		We now present the first estimation result which yields 
		consistent estimators of stickiness and skewness parameters.
		
		\begin{proposition}
			\label{thm_estimators_SkSBM}
			Let $X$ be the SOS-diffusion that solves~\eqref{eq_SkSID_pathwise_char} on $\mc P_x$, 
			and let a non-trivial function $g$ and a sequence $(u_n)_n$ both satisfying the assumptions of 
			Theorem~\ref{thm_main}. 
			We consider the following estimators
			\begin{equation}
				\label{eq_thm_estimators_SkSBM}
				\wh\rho_n(X):=\frac{2}{n}\frac{\sum_{i=1}^{[nt]}\indicB{\hfprocess{X}{}{i-1}=0}}{S_{n}^{g+}(X) + S_{n}^{g-}(X)},
				\quad
				\wh\beta_n(X) := \frac{S_{n}^{g+}(X) - S_{n}^{g-}(X)}
				{S_{n}^{g+}(X) + S_{n}^{g-}(X)}.   
			\end{equation}
			If $\mc L_t := \{\tau^{X}_{0}<t\}$, with $\tau^{X}_{0} $ the hitting time of $0$ by $X$, then $\wh\rho_n(X) $ and $\wh\beta_n(X) $ are respectively $\Prob_{x}^{\mc L_t} $-conditional consistent estimators of
			$\rho $ and skew $\beta $, i.e.,
			\begin{equation}
				\label{eq_consistency}
				\wh\rho_n(X) \xlra{\Prob_{x}^{\mc L_t}}{n\ra \infty} \rho,
				\quad \text{and}\quad 
				\wh\beta_n(X) \xlra{\Prob_{x}^{\mc L_t}}{n\ra \infty} \beta. 
			\end{equation}
		\end{proposition}
		
		Let us note that the estimators can only be consistent on $\mc L_{t} $. 
		Indeed, on $\mc L_{t}^{c} $ we have that $ \occt{X}{0}{t} = \loct{X}{0}{t}=0$, therefore the estimators are asymptotically $0/0 $
		divisions. 
		
		For SOS-BM, we additionally construct consistent estimators of the two volatility levels. Together with the stickiness and skewness estimators, these yield joint consistent estimation of all four parameters $(\rho,\beta,\sigma_-,\sigma_+) $ of SOS-BM. 
		
		\begin{proposition}
			\label{prop_SOSBM_estimation}
			Let $X$ be the SOS-BM that solves~\eqref{eq_SOSBM_pathwise_char} on $\mc P_x $
			and $\tau_+,\tau_- $ be the following hitting times:
			\begin{equation}
				\begin{aligned}
					\tau_{+} &:= \inf \{s\ge 0: X_s>0\},
					&
					\tau_{-} &:= \inf \{s\ge 0: X_s<0\},
				\end{aligned}
			\end{equation}
			$\mc L_{t}^{+}:=\{\tau_+ <t\} $ and $\mc L_{t}^{-}:=\{\tau_- <t\} $.
			Then, for all $t>0 $, in $\Prob_x$-probability, 
			\begin{align}
				\sqrt{\frac{\sum_{i=1}^{[nt]} (\hfprocess{X}{+}{i}-\hfprocess{X}{+}{i-1})^{2}}{\frac{1}{n} \sum_{i=1}^{[nt]} \indic{(0,\infty)}(\hfprocess{X}{}{i-1})}}
				\xlra{\Prob_{x}^{\mc L^{+}_t}}{n\ra \infty}
				\sigma_{+},
				\quad
				\sqrt{\frac{\sum_{i=1}^{[nt]} (\hfprocess{X}{-}{i}-\hfprocess{X}{-}{i-1})^{2}}{\frac{1}{n} \sum_{i=1}^{[nt]} \indic{(-\infty,0)}(\hfprocess{X}{}{i-1})}}
				\xlra{\Prob_{x}^{\mc L^{-}_t}}{n\ra \infty} 
				\sigma_{-},
			\end{align}
			where $X^+=\max\{0,X\}$ and $X^-=\max\{0,-X\}$.
			Whenever one of the denominators in the preceding display vanishes, the corresponding estimator is defined to be zero.
		\end{proposition}
		
		\section{Discussion, related statistics, and extensions}
		\label{sec_comments_extensions}
		
		Our main results concern local-time approximation and parameter
		estimation at a single SOS threshold for a one-dimensional
		diffusion. 
		In this section, we first discuss related crossing statistic and the unresolved question of convergence
		rates in the sticky setting. 
		We then describe two possible one-dimensional
		extensions: to multiple isolated
		thresholds and to semi-time-inhomogeneous coefficients. Finally, we
		relate our results to the existing literature on multidimensional
		sticky diffusions and discuss the additional difficulties arising
		in that setting.
		
		\subsection{Related statistics and the critical regime}
		
		The present paper concerns the one-point statistic \eqref{eq_intro_loct_approximation}. More general high-frequency statistics based on two consecutive observations include crossing counts. For an account of this literature, we refer to the authors later work~\cite{Anagnostakis_Mazzonetto_2026} and the references therein.
		
		In~\cite{Anagnostakis_Mazzonetto_2026}, the authors distinguish three notions of crossing at a sticky threshold.
		Strict crossings, for which the two consecutive observations lie
		in opposite open half-lines, are negligible after division by any
		diverging sequence. Crossings that allow one observation to equal
		the threshold, but exclude pairs for which both observations equal
		the threshold, converge after normalization by $n^{-1/2}$ to
		the local time up to a constant (which does not depend on the stickiness). Finally, crossings that also count pairs of consecutive observations both equal to the threshold converge after normalization by $n^{-1}$ to the occupation time at the threshold. 
		The critical-regime local-time approximation established in
		Proposition~\ref{thm_the_case_a_05} and the occupation times approximations in Lemma~\ref{lem_sosdiff_occtime_approximation} are the main ingredients in
		these arguments.
		
		The present paper proves consistency of
		\eqref{eq_intro_loct_approximation} but does not determine its
		optimal convergence rate. For regular, skew, and oscillating
		non-sticky diffusions, the normalizing sequence $u_n=\sqrt n$ is the critical regime and
		yields the optimal rate; see \cite{Jac98,Maz19}. Whether the same
		optimality statement holds in the sticky setting remains open. 
		
		Numerical experiments in~\cite{Anagnostakis2022} investigate the
		local-time approximation and the associated stickiness estimator
		for normalizing sequences of the form $u_n=n^\alpha$. 
		Therein, numerical evidence beyond the
		subcritical range covered theoretically in that paper and motivate
		the conjecture that consistency holds for every
		$\alpha\in(0,1)$. The results of the present paper prove this
		conjecture, and more generally cover every diverging sequence
		$u_n=o(n)$.
		
		The paper~\cite{Anagnostakis_Mazzonetto_2026} also contains a numerical study comparing its crossing-based stickiness estimator with the local-time-based estimator of~\cite{Anagnostakis2022}. These
		experiments examine finite-sample behavior, including empirical
		variance and apparent convergence speed, but do not establish a
		theoretical optimal rate.

		\subsection{Possible extensions within one dimension}
		
		We briefly indicate two directions in which the localization and
		transformation arguments of the present paper may be extended.
		These extensions are not used elsewhere in the paper, and complete
		proofs are beyond its scope.
		
		\paragraph{Multiple isolated thresholds.}
		The approximation in Theorem~\ref{thm_main} is local around the threshold, in the sense that the spatial rescaling of the statistic and the $L^1$-integrability of the test function concentrate its asymptotic contribution in a shrinking neighborhood of the threshold. 
		This suggests an extension to diffusions with a finite
		or locally finite collection of isolated SOS thresholds. For a
		fixed threshold, one may stop the process upon leaving a
		neighborhood containing no other singular point and then apply the
		corresponding one-threshold approximation. A rigorous formulation
		would require verifying that this localization is compatible with
		the statistic and with the chosen normalization of local time. We
		leave this extension for future work.
		
		\paragraph{Semi-time-inhomogeneous coefficients.}
		A second possible extension concerns sticky SDEs with
		time-dependent coefficients for which the diffusion coefficient at
		the threshold is constant in time:
		\[
		\sigma(t,0)=\eta,
		\qquad t\geq0.
		\]
		Under suitable smoothness and non-degeneracy assumptions, the
		time-dependent Lamperti transform
		\[
		T(t,x)
		=
		\int_0^x\frac{\rd y}{\sigma(t,y)}
		\]
		suggests a reduction to sticky Brownian motion after an appropriate
		change of probability measure. A complete proof would have to control the additional drift
		generated by $\partial_tT$, the discretization error caused by the
		time dependence of the transformation, and the corresponding
		transformation of local time. Besides being of independent
		interest, such an extension would provide a natural intermediate
		step toward the multidimensional setting discussed below, where
		the dynamics normal to a sticky hypersurface generally involve
		time- and state-dependent coefficients.

		\subsection{Multidimensional sticky diffusions}
		
		Sticky diffusions have also been studied in multidimensional
		settings, where the stickiness occurs on a lower-rank surface instead of a threshold. 
		Grothaus and Vo{\ss}hall~
		\cite{grothaus2017_stochastic_differential} construct diffusions in
		bounded domains that spend positive time on the boundary and may,
		while remaining there, evolve according to a diffusion on the
		boundary itself. Related models on smooth bounded domains,
		including fractional dynamic boundary conditions, are considered
		by D'Ovidio~\cite{dovidio2026_fractional_boundary}.
		Multidimensional sticky-reflected Brownian motions also arise in
		wedges and orthants, for example as scaling limits of interacting
		particle systems. Further classes of nonlinear diffusions with
		sticky reflection are studied in
		\cite{graham1995_homogenization_propagation}.
		
		Extending the local-time approximation of the
		present paper to such processes is nevertheless not immediate. 
		In dimension greater than one, the singular set is typically a boundary or a hypersurface rather than an isolated point. The local time at a point must therefore be replaced by a boundary or surface
		local time, and the statistic must be formulated in terms of the distance to, or local coordinates around, the singular set.
		
		For a smooth hypersurface $\Gamma$, a possible approach is to
		flatten $\Gamma$ locally and separate the normal coordinate from
		the tangential coordinates. In special cases where the normal
		component is autonomous and the statistic depends only on the
		signed distance to $\Gamma$, the problem may reduce directly to a
		one-dimensional sticky diffusion.
		
		In the general case, however, the tangential components may
		continue to evolve while the process remains on $\Gamma$, and they
		may interact with the normal dynamics. Flattening the hypersurface
		would then produce time- and state-dependent coefficients,
		curvature terms, and coupling between the normal and tangential
		coordinates. A multidimensional extension would consequently
		require, as an intermediate step, a weighted and
		time-inhomogeneous version of the local-time approximation proved
		here.
		
		Determining the appropriate normalization, identifying the
		limiting surface-local-time constant, and constructing estimators
		of the corresponding stickiness parameters constitute interesting
		open problems.

		\section{Preliminary results on SkS-BM}
		\label{sec_preliminary}

		In this section we present several results on SkS-BM that are key for our derivations.
		First, we present the time-space scaling property of the SkS-BM kernel.
		Second, we present some bounds on SkS semigroup.
		Last, we prove a preliminary local time approximation based on a sequence of test functions $ \wh g_n$ and present an asymptotic of $\wh g_n $.
		All proofs (except that of the local time approximation) are deferred
		to Appendix~\ref{app_proofsSemigroup}.

		Given $\rho \ge 0 $ and $\beta\in (-1,1) $, let $X^{\params{\rho}{\beta}} $ be the SkS-BM,
		that is, the unique (in law) weak solution to the following system, which is a special case to~\eqref{eq_SkSID_pathwise_char} (it is \eqref{eq_SOSBM_pathwise_char} with $\sigma_0 \equiv 1 $):
		\begin{equation}
			\label{eq_SkSBM_pathwise_char}
			\begin{dcases}
				\vd X^{\params{\rho}{\beta}}_t = \indicB{X^{\params{\rho}{\beta}}_t \not = 0}  \vd B_t
				+ \beta \vd \loct{X^{\params{\rho}{\beta}}}{0}{t},\\
				\indicB{X^{\params{\rho}{\beta}}_t = 0}\vd t = \rho \vd \loct{X^{\params{\rho}{\beta}}}{0}{t},
			\end{dcases}
		\end{equation}
		where $B$ is a standard BM and \( L^0(X^{\rho,\beta}) \) is the local time of \( X^{(\rho,\beta)} \) at zero.
		
		Let us introduce the following notations.
		
		\begin{notation} \label{notation:SkS}
			Let
			\begin{itemize}
				\item \( m_{(\rho,\beta)} \) be the speed measure of \( X^{(\rho,\beta)} \), which is special case of~\eqref{eq_sm_SOSBM} and  admits the decomposition
				\begin{equation}\label{eq_speed_measure_decomp}
					m_{(\rho,\beta)}(\mathrm{d} x) = m_{(0,\beta)}(\mathrm{d} x) + \rho \, \delta_0(\mathrm{d} x).
				\end{equation}
				The measure \( m_{(0,\beta)} \) is absolutely continuous with respect to Lebesgue measure \( \lambda \): $m_{(0,\beta)}(\rd x) = (1+\sgn(x)\beta) \vd x$;
				
				\item \( p_{(\rho,\beta)}(t,x,y) \) be the transition density of \( X^{(\rho,\beta)} \) with respect to \( m_{(\rho,\beta)} \), whose existence follows from \cite[p.~151]{ItoMcKean96}. We provide it explicitly in the next section.
				
				\item If $X^{(\rho,\beta)} $ is the solution to~\eqref{eq_SkSBM_pathwise_char} on $\mc P_x $, then the transition semigroup $\process{P^{(\rho,\beta)}_t} $ reads  
				\begin{equation}
					P^{(\rho,\beta)}_t f(x)
					:= \Esp_x \sqbraces{f(X^{(\rho,\beta)}_t)}
					= \int_{\IR} p_{(\rho,\beta)}(t,x,y) f(y)\, m_{(\rho,\beta)} (\rd y),
				\end{equation}
				for all $t>0 $, $x\in \IR $, and every suitable $f$.
			\end{itemize}
		\end{notation}

		\subsection{The transition density of SkS-BM}
		
		The probability transition kernel $ p_{\params{\rho}{\beta}} $ of the $\params{\rho}{\beta}{}$-SkS-BM with respect to $ m_{\params{\rho}{\beta}} $ was  computed in \cite[Theorem 2.4]{Touhami2021}:
		\begin{equation}\label{eq_def_SkSBM_kernel_factorization}
			p_{\params{\rho}{\beta}}(t,x,y) 
			= 
			\frac{1}{1+ \sgn(y)\beta} \braces{u_1(t,x,y) - u_2(t,x,y)} + v_{\rho}(t,x,y),
		\end{equation}
		for all $t> 0 $ and $x,y\in \IR $, where $a$  is the function defined in~\eqref{eq_sm_SOSBM}, and $u_1,u_2,v_{\rho} $ are defined for all 
		$t>0 $, $x,y\in \IR $ by 
		\begin{equation}\label{eq_def_SkSBM_kernel_factors}
			\begin{dcases}
				u_1(t,x,y) =  \braces{2 \pi t}^{-1/2}  e^{-(x-y)^2/2t} ,\\
				u_2(t,x,y) = \braces{2 \pi t}^{-1/2} e^{-(|x|+|y|)^2/2t},\\
				v_{\rho}(t,x,y) = \frac{1}{\rho} e^{2(|x|+|y|)/\rho + 2 t/\rho^2} \erfc \braces{\frac{|x|+|y|}{\sqrt{2t}} + \frac{\sqrt{2t}}{\rho}}, &\forall \rho>0,\\
				v_{0}(t,x,y) = u_2(t,x,y).
			\end{dcases}
		\end{equation}
		
		\begin{remark}
			\begin{enumerate}
				\item For $\rho>0$ and $xy\le 0 $, the expression of the probability transition kernel simplifies to   
				$p_{\params{\rho}{\beta}}(t,x,y) = v_{\rho}(t,x,y)$.
				\item The term $v_{\rho}(t,x,y) $ does not appear in the probability transition kernel of the standard BM and skew BM.  
				The  probability transition kernel of the skew BM for $\beta\not = 0$ and  of the standard BM with respect to their speed measure
				($m_{(0,\beta)}(\rd x) = (1+ \sgn(x)\beta) \vd x $ and $m_{(0,0)}(\rd x) = \vd x $), are 
				respectively
				\begin{align}
					& p_{\params{0}{\beta}}(t,x,y)=
					\frac{1}{1+ \sgn(y)\beta} \braces{u_1(t,x,y) - u_2(t,x,y)} + u_2(t,x,y)
					\\
					\text{and} \quad & p_{\params{0}{0}}=
					u_1(t,x,y).
				\end{align}
			\end{enumerate}
		\end{remark}
		
		\subsection{Scaling of the SkS-BM kernel}
		\label{ssec: scaling}
		
		Here, we present the time-scaling property of the SkS-BM kernel and  semigroup.
		The proof of the kernel scaling is deferred to Appendix~\ref{aapp: proof_TS_scaling}.
		
		\begin{proposition}
			[Kernel scaling]
			\label{cor_semi_group_SOS-BM_timescaling_ptk}
			For every sufficiently integrable $h:\IR^2\to \IR $,
			\begin{equation}
				\int_{\IR} h(x,y) p_{\params{\rho}{\beta}}(ct,x,y) m_{\params{\rho}{\beta}}(\vd y)
				= 
				\int_{\IR} h(x,\sqrt{c}y) p_{\params{\rho/\sqrt{c}}{\beta}}
				\braces{t,\frac{x}{\sqrt{c}},y}m_{\params{\rho/\sqrt{c}}{\beta}}(\vd y).
			\end{equation}
		\end{proposition}
		
		\begin{corollary}[Semigroup scaling]\label{cor_semi_group_SOS-BM_timescaling}
			For every sufficiently integrable function $h:\IR \mapsto \IR$, we have that
			\begin{equation}
				P_{t}^{\params{\rho\sqrt n}{\beta}} h(x \sqrt n) = 
				\Esp_x \braces{h(\sqrt n \hfprocess{X}{\params{\rho}{\beta}}{t}) } .
			\end{equation} 
		\end{corollary}
		
		\begin{proof}
			It is a special case in  Proposition~\ref{cor_semi_group_SOS-BM_timescaling_ptk}.
		\end{proof}
		
		\subsection{Bounds on the SkS-BM semigroup}
		
		Here, we present bounds on the SkS-BM kernel and semigroup.
		Proofs are deferred to Appendix~\ref{aapp: proof_kernel_bounds}.

		\begin{lemma}\label{lem_kernel_bound}
			There exists a constant $K>0$ such that for all $t>0 $ and  $x,y\in \IR $, 
			\begin{equation}\label{eq_lem_kernel_bound}
				p_{\params{\rho}{\beta}}(t,x,y) \le K u_1(t,x,y),
			\end{equation}
			where $K$ does not depend on $\rho$, but depends on $\beta$. 
			In particular there exists a constant $K>0$ (not depending on $\rho$, nor on $
        \beta$) such that  
			\begin{equation}
				v_{\rho}(t,x,y)  \le K u_1(t,x,y),
			\end{equation}
			for all $t>0, \rho\geq 0$ and  $x,y\in \IR $.
		\end{lemma}

		\begin{lemma}
			\label{prop_semig_bound}
			There exists a constant $K>0$ (depending on $\beta$ but not on $\rho$) such that 
			\begin{equation}
				|P^{\params{\rho}{\beta}}_t h(x)| 
				\le K  \frac{\rho \sqrt{2 t}}{ \rho |x|/2 + 2 t} |h(0)| + \frac{K}{\sqrt t}m_{\params{0}{\beta}}(|h|)
				\le \frac{K}{\sqrt t} m_{\params{\rho}{\beta}}(|h|),
			\end{equation}
			for all $h\in L^{1}(\IR)$, $t>0 $, and $x\in \IR $. 
		\end{lemma}
		
		\begin{lemma} \label{lem:semig:bound:strong}
			For every $\gamma \ge 0 $, there exists a constant $K >0 $, that depends on $\beta $ and $\gamma$ (and not on $\rho $), such that
			\begin{equation}\label{eq_sticky_semigroup_bound2}
				\abs{P^{\params{\rho}{\beta}}_th(x) - m_{\params{\rho}{\beta}}(h)p_{\params{\rho}{\beta}}(t,x,0)} 
				\le K \frac{1}{t}
				\braces{m_{\params{0}{\beta}}^{(1)}(h) + \frac{m_{\params{0}{\beta}}^{(1)}(h)}{1 + |x/\sqrt t|^{\gamma}} + \frac{m_{\params{0}{\beta}}^{(\gamma)}(h)}{1+ |x|^{\gamma}}},
			\end{equation}
			for all $h\in L^{1}(\IR)$, $t>0 $, and $x\in \IR $. 
		\end{lemma}
		
		We now present the last semigroup estimates we need for our main proof (of Theorem~\ref{thm_main}).
		The result is expressed in terms of the following functional, which aggregates the action of the semigroup over the sampling data. This functional arises naturally when we consider moments of our statistics.  
		It is defined for every suitable (e.g.~bounded or positive) measurable function $h$, $t>0 $, $n\in \IN $, and $x\in \IR $, as
		\begin{equation}
			\label{eq:gamma_semigroup}
			\gamma^{\params{\rho}{\beta}}_n[h](t,x) := \sum_{i=2}^{[nt]} \Esp_x \braces{ h(\sqrt n \hfprocess{X}{\params{\rho}{\beta}}{i-1}) }
			= \sum_{i=2}^{[nt]} P^{\params{\sqrt n \rho}{\beta}}_{i-1} h(\sqrt n x),
		\end{equation}	
		where the last equality is consequence of the time-space scaling of $X^{\params{\rho}{\beta}} $, see  Corollary~\ref{cor_semi_group_SOS-BM_timescaling}.

		\begin{lemma} 
			\label{lem_gamma_bounds_SOS-BM}
			\begin{enuroman}
				\item \label{eq_def_gamma_bound_A_os_SOS-BM}
				There exists a positive constant $K>0$ (depending on $\beta$ and not on $\rho$) such that,
				for every $t>0 $, $n\in \IN $, and suitable function $h$:
				\begin{equation}\label{eq_lem_aggregate_action_bound1}
					|\gamma^{\params{\rho}{\beta}}_n[h](x,t)|  \le 
					K m_{\params{\rho \sqrt n}{\beta}}(|h|) \sqrt{n t},
				\end{equation}
				
				\item \label{eq_def_gamma_bound_A_os_SOS-BM2} 
				Assume that $m_{\params{\sqrt n \rho}{\beta}}(h)=0$, for all $n\in \IN $. Then, there exists a positive constant $K >0$ (depending on $\beta $ and not on $\rho$) such that,
				for every $t>0 $, $n\in \IN $:
				\begin{equation}\label{eq_def_gamma_bound_B_os_SOS-BM}
					|\gamma^{\params{\rho}{\beta}}_n[h](x,t)|  \le K  m_{\params{0}{\beta}}^{(1)}(h)  (1+ \max(0,\log(n t))).
				\end{equation}
			\end{enuroman}
		\end{lemma}
		
		\subsection{A preliminary local time approximation}
		\label{subsec: preliminary_ltime}
		
		In this section we establish a preliminary local time approximation result based on the Tanaka formula.
		In Proposition~\ref{prop_mngn_convergence}, we also present a property of this approximation, that will prove useful for our further derivations, and whose proof is deferred to Appendix~\ref{aap: proof_lta_limit}.
		
		The statistic used for our first approximation of the local time 
		is the rescaled mean absolute displacement of the process.
		Given the parameters $\rho\ge 0 $ and $\beta\in (-1,1) $, the statistics is constructed from the sequence of functions 
		\begin{equation}
			\label{eq_ghn_first_def}
			\wh g_{n}(x) := \Esp_{x} \sqbraces{|X^{\params{\rho \sqrt n}{\beta}}_{1}| - |x|},
			\quad n\in \IN,
		\end{equation}
		where $X^{\params{\rho \sqrt n}{\beta}}$ is the SkS-BM on $\mc P_x$ for all $x\in \IR$. Recall that the state space is $\IR$, since $\beta\in (-1,1)$.  
		\begin{proposition}
			\label{prop_first_loct_approximation}
			Let $X^{\params{\rho }{\beta}}$ be an SkS-BM on $\mc P_x$. The following convergence holds:
			\begin{equation}
				\frac{1}{\sqrt n} \sum_{i=1}^{[n\cdot]} \widehat{g}_{n}(\sqrt n \hfprocess{X}{\params{\rho }{\beta}}{i-1}) 
				\xlra{\Prob_x\text{-ucp}}{n\ra \infty} \loct{X^{\params{\rho }{\beta}}}{0}{}.
			\end{equation}
		\end{proposition}
		\begin{proof}
			Since the process $ \process{X^{\params{\rho }{\beta}}_t}$ is a semimartingale,
			the (symmetric) Tanaka formula
			(see e.g. It\^o--Tanaka formula in Lemma \ref{lem:ito-tanaka}) ensures that
			\begin{equation}
				|X^{\params{\rho }{\beta}}_t|-|x|
				= \int_{0}^{t} \sgn(X^{\params{\rho }{\beta}}_s) \indicB{X_s \not = 0} \vd B_s
				+  \loct{X^{\params{\rho }{\beta}}}{0}{t}.
			\end{equation}
			Since $\int_{0}^{t} \sgn(X^{\params{\rho }{\beta}}_s) \indicB{X_s \not = 0} \vd B_s $ is a martingale,
			\begin{equation}
				\Esp_x \sqbraces{|X^{\params{\rho }{\beta}}_{t}| - |X^{\params{\rho }{\beta}}_{s}|} = \Esp_x \sqbraces{\loct{X^{\params{\rho }{\beta}}}{0}{t}-\loct{X^{\params{\rho }{\beta}}}{0}{s}}.
			\end{equation}
			Also, from the scaling property (Corollary~\ref{cor_semi_group_SOS-BM_timescaling}),
			\begin{equation}
				\frac{1}{\sqrt n}\widehat{g}_{n}(\sqrt n \hfprocess{X}{\params{\rho }{\beta}}{i-1}) = \Esp_x \sqbraces{|\hfprocess{X}{\params{\rho }{\beta}}{i}| - |\hfprocess{X}{\params{\rho }{\beta}}{i-1}| \big| \hfprocess{\bF}{}{i-1}}.
			\end{equation}
			Thus, from \cite[Lemma 2.14]{Jac84},
			\begin{equation}
				\frac{1}{\sqrt n} \sum_{i=1}^{[nt]} \widehat{g}_{n}(\sqrt n \hfprocess{X}{\params{\rho }{\beta}}{i-1}) 
				= \sum_{i=1}^{[nt]}	\Esp_x \sqbraces{\loct{X^{\params{\rho }{\beta}}}{0}{\frac{i}{n}}-\loct{X^{\params{\rho }{\beta}}}{0}{\frac{i-1}{n}}\big| \hfprocess{\bF}{}{i-1}} \convergence{\Prob_x } \loct{X^{\params{\rho }{\beta}}}{0}{t}
			\end{equation}
			and, from Lemma~\ref{lem_ucp_convergence_condition}, the convergence is locally uniform in time, in probability (\ucp).
		\end{proof}
		
		\begin{proposition}
			\label{prop_mngn_convergence}
			It holds that $m_{\params{\rho \sqrt n}{\beta}} (\wh g_n) \to 1 $
			as $n\to \infty $.
		\end{proposition}
		The latter proposition is proven in Appendix~\ref{aap: proof_lta_limit}. 
		
		\section{Proof of the main result for SkS-BM, Theorem~\ref{thm_main_SkSBM}}
		\label{sec_SkSBM}
		
		Throughout this section, let ${X}^{\params{\rho}{\beta}} $ be a SkS-BM (i.e.,~solution to~\eqref{eq_SkSBM_pathwise_char}) on $\mc P_x$,  
		and $m_{\params{\rho}{\beta}}$ be its speed measure (see~\eqref{eq_sm_SOSBM}).
		Recall that $\beta\in (-1,1)$.  
		The assumptions for $T$ and $\un$ and the notation $g_n[T]$ of Theorem~\ref{thm_main_SkSBM} remain unchanged.  
		We now prove the convergence 
		\begin{align}
			\label{eq_thm_limit_gnT_inline}
			\frac{\un}{n} \sum_{i=1}^{[n\cdot]} g_n[T](\un \hfprocess{X}{\params{\rho }{\beta}}{i-1}) &\xlra{\Prob_x\text{-ucp}}{n\ra \infty}
			\left( \int_{\IR} g \vd m_{\params{0}{\beta}}\right) \loct{{X}^{\params{\rho}{\beta}}}{0}{}
		\end{align}
		separately different regimes of $(u_n)_n $, clearly among those satisfying~\eqref{eq_intro_un_condition}, that is $\un \to \infty$ and $\un = o(n)$. 
		\begin{itemize}
			\item In Section~\ref{ssec_un_case1}, we show the above convergence, if $T$ is the identity, for the \emph{critical regime} $u_n = \sqrt{n} $. 
            We remind the reader that this result is new already in the simpler case of sticky BM. 
			\item In Section~\ref{ssec_un_case2}, we prove convergence under the condition
			$ \log n = o(u_n)$,  
			which includes all the power regimes (\emph{critical}, \emph{subcritical}, \emph{supercritical}) $\un=n^\alpha$ for all $\alpha \in (0,1)$ and in particular the \emph{supercritical power regimes} \(\un=n^\alpha\) with \(\alpha\in(1/2,1)\), which have never been treated before.  
			\item In Section~\ref{ssec_un_case3}, we prove convergence under the condition 
			${\un^2}/{n}\to 0$, which includes the \emph{subcritical regime} $\un=n^\alpha$ with $\alpha\in (0,1/2)$ and, in particular, all sequences satisfying \(\un=O(\log n)\). The latter do not fall into the setting of Section~\ref{ssec_un_case2}. 
		\end{itemize}
        
		Even though the \emph{subcritical power regimes} fall into both the second and the third case, the corresponding proofs are quite different.
        For sticky BM, the third case has been considered in~\cite{Anagnostakis2022}.  

        The conditions in the second and third cases exhaust all admissible diverging sequences $(u_n)_n$ satisfying $u_n=o(n)$, up to subsequence extraction. 
        Indeed, let $(n_k)_k$ be an arbitrary subsequence.
        If $(u_{n_k}/\log n_k)_k$ is unbounded, one may extract a further subsequence, still denoted by $(n_k)_k$, such that
        ${\log n_k}/{u_{n_k}}\to0$. 
        Hence the condition of Section~\ref{ssec_un_case2} holds along this subsequence.
        Otherwise, $(u_{n_k}/\log n_k)_k$ is bounded. Thus, after possibly passing to a further subsequence, there exists $C>0$ such that 
        $u_{n_k}\leq C\log n_k$. Consequently,
        ${u_{n_k}^2}/{n_k} \leq C^2 {(\log n_k)^2}/{n_k} \to 0$, so the condition of Section~\ref{ssec_un_case3} holds.
        The general conclusion then follows from the subsequence criterion for convergence in probability, applied to the supremum norm on every compact time interval.
		
		The ideas of the proofs and their main steps are detailed in the respective subsections. 
		
		\subsection{The critical regime: $u_n = \sqrt n $}\label{ssec_un_case1}
		
		In this section we show a result that is stronger than the convergence~\eqref{eq_thm_limit_gnT_inline} in the critical regime $u_n = \sqrt{n} $, when $T$ is the identity function.
		
		For convenience, we fix $x$, and we let ${X}^{\params{\rho}{\beta}} $ be solution to~\eqref{eq_SkSBM_pathwise_char} on $\mc P_x$.
		Recall the notation introduced in~\eqref{eq: measure_integral_notation}   
		and introduce the sequence of functionals
		$(B_n)_n$ defined for all bounded measurable $f$, and $n\in \IN $ as
		\begin{align}
			B_n[f](x) :={}& \frac{f^2(\sqrt n x)}{n} + 
			\frac{m_{\params{\rho \sqrt n}{\beta}}(f^{2})}{\sqrt n}  
			\\ & +  \frac{  \braces{1 + \log(n)} m_{\params{0}{\beta}}^{(1)}(f) \abs{ f(\sqrt n x) }}{n}
			+   \frac{{  \braces{1 + \log(n)} m_{\params{0}{\beta}}^{(1)}(f) m_{\params{\rho \sqrt n}{\beta}}(|f|)}}{\sqrt n}. 
			\label{eq_connd_gn_aggreg}
		\end{align}
		
		The main result of this section, Proposition~\ref{thm_the_case_a_05}, works under the following condition on the sequence of test functions $(g_n)_n$. 
		
		\begin{condition}
			\label{cond_B_vanish}
			The real sequence $n\mapsto B_n[g_n](x)$ vanishes as $n\to \infty $.
		\end{condition}
		
		\begin{remark}\label{rmk_condition_sum_of_fun_sequences}
	        If a sequence of functions $(f_n)_n$ satisfies Condition~\ref{cond_B_vanish}, then so do $(f_n^+)_n$ and $(f_n^-)_n$, where $f_n^+(y)=\max\{f_n(y),0\}$ and $f_n^-(y)=\max\{-f_n(y),0\}$. 
            The condition is also preserved under multiplication by a bounded sequence of real numbers.
            If a measurable function $f$ is bounded, integrable, $\int_{-\infty}^{+\infty}|f(x)x|\,\vd x<\infty$, and $f(0)=0$, then the constant sequence $(f)_n$ satisfies Condition~\ref{cond_B_vanish}.
		\end{remark}
		
		We can now present the result, which reads as follows.
		
		\begin{proposition}
			\label{thm_the_case_a_05}
			Let $(g_{n})_{n}$ be a sequence of functions that satisfies Condition~\ref{cond_B_vanish} and $\lim_{n\to \infty}m_{\params{\sqrt n\rho}{\beta}}(g_n) = M$.
			Then, for all $t\ge 0 $, 
			\begin{equation}\label{eq_the_case_a_05_plain}
				\frac{1}{\sqrt n} \sum_{i=1}^{[nt]}g_n(\sqrt n \hfprocess{X}{\params{\rho}{\beta}}{i-1}) \xlra{\Prob_x}{n\ra \infty} 
				M \loct{X^{\params{\rho}{\beta}}}{0}{t}.
			\end{equation}
			Also, if $\sup_n \cubraces{m_{\params{\sqrt n\rho}{\beta}}(|g_n|)} <\infty$, then the above convergence is $\Prob_x$-ucp.
		\end{proposition}
		
		Let us provide a step-by-step guide to the proof and the intermediary results.
		
		\begin{enumerate}
			\item For every $n$, consider the linearization of $g_n $ with respect to the preliminary local time approximation statistic given in  Proposition~\ref{prop_first_loct_approximation}:
			\begin{equation}
				\label{eq_proof_linearization_wrt_loct0}
				\begin{aligned}
					h_{n}(y) &= g_{n}(y) - \frac{m_{\params{\sqrt n\rho}{\beta}}(g_n)}{m_{\params{\sqrt n\rho}{\beta}}(\widehat{g}_n)} \widehat{g}_n(y),\qquad 
					&
					\wh g_{n}(y) &= \Esp_{y} \braces{|X^{\params{\rho \sqrt n}{\beta}}_{1}| - |y|}.
				\end{aligned}
			\end{equation}
			By construction, $m_{\params{\sqrt n\rho}{\beta}}(h_n)=0$. 
        	\item We then establish the two estimates needed to control the centered remainder $h_n$.
            Lemma~\ref{lem_whgn_sqtisfies_condition_38} shows that $(\widehat g_n)_n$ satisfies Condition~\ref{cond_B_vanish}, while 	Lemma~\ref{lem_L2_funct_convergence} shows that any centered sequence satisfying this condition gives a negligible statistic in $L^2(\Prob_x)$.
            \\
            The proof of Lemma~\ref{lem_L2_funct_convergence} is the only place where the new kernel bounds 	Lemmas~\ref{lem:semig:bound:strong} and 	\ref{lem_gamma_bounds_SOS-BM}(ii) are used.

	\item A direct verification using Condition~\ref{cond_B_vanish} and Lemma~\ref{lem_whgn_sqtisfies_condition_38} shows that $(h_n)_n$ satisfies Condition~\ref{cond_B_vanish}. 
    Since, moreover
	$m_{\params{\sqrt n\rho}{\beta}}(h_n)=0$,
	Lemma~\ref{lem_L2_funct_convergence} yields, for every $t>0$,
	\[
		\frac{1}{\sqrt n}
		\sum_{i=1}^{[nt]}
		h_n\left(
			\sqrt n\,
			\hfprocess{X}{\params{\rho}{\beta}}{i-1}
		\right)
		\xrightarrow{\Prob_x}0.
	\]
	Hence the statistic associated with $g_n$ is asymptotically
	equivalent to
	\[
		\frac{
			m_{\params{\sqrt n\rho}{\beta}}(g_n)
		}{
			m_{\params{\sqrt n\rho}{\beta}}(\widehat g_n)
		}
		\frac{1}{\sqrt n}
		\sum_{i=1}^{[nt]}
		\widehat g_n\left(
			\sqrt n\,
			\hfprocess{X}{\params{\rho}{\beta}}{i-1}
		\right).
	\]

	\item Finally, Proposition~\ref{prop_first_loct_approximation}
	gives
	\[
		\frac{1}{\sqrt n}
		\sum_{i=1}^{[n\cdot]}
		\widehat g_n\left(
			\sqrt n\,
			\hfprocess{X}{\params{\rho}{\beta}}{i-1}
		\right)
		\xrightarrow{\Prob_x\text{-ucp}}
		\loct{X}{0}{},
	\]
	while Proposition~\ref{prop_mngn_convergence}, together with the
	assumption on $m_{\params{\sqrt n\rho}{\beta}}(g_n)$, identifies
	the limit of the multiplicative coefficient.
	This yields the limit $ M\loct{X}{0}{t}$.
\end{enumerate}

		Before proving Proposition~\ref{thm_the_case_a_05}, we establish the two preliminary results used in the above argument. 
		The first result is an adaptation to this context of~\cite[Lemma~4.2]{Jac98}.		 
		
		\begin{lemma}\label{lem_L2_funct_convergence}
			Let $(g_{n})_{n}$ be a sequence of functions that satisfies Condition~\ref{cond_B_vanish} and for all $n\in \IN$: 
			$ m_{\params{\sqrt n \rho}{\beta}}(g_n) =0$. Then,  	
			\begin{equation}\label{eq_L2_funct_convergence}
				\lim_{n\to \infty}	\Esp_x \braces{
					\abs{
						\frac{1}{\sqrt n} \sum_{i=1}^{[nt]}g_n(\sqrt n \hfprocess{X}{\params{\rho}{\beta}}{i-1})
					}^2
				} = 0.
			\end{equation}
		\end{lemma}

		\begin{proof}
			In this proof we use the notation for the functional $\gamma^{\params{\rho}{\beta}}_n$ given by~\eqref{eq:gamma_semigroup}. 
			It holds that 
			\begin{align}
				\Esp_x & \braces{
					\abs{
						\sum_{i=1}^{[nt]}g_n(\sqrt n \hfprocess{X}{\params{\rho}{\beta}}{i-1})
					}^2
				} \le{}  g_{n}^2(\sqrt n x) + \gamma^{\params{\rho}{\beta}}_n [g^{2}_n](t,x) 
				\\ & + 2 \braces{\sup_{y\in \IR; s\le t}  \cubraces{ |\gamma^{\params{\rho}{\beta}}_n [g_n](s,y)| } } 
				\braces{ \abs{ g_n(\sqrt n x)}  + \abs{ \gamma^{\params{\rho}{\beta}}_n [g_n](t-1/n,x)}}.
			\end{align}
			From Lemma~\ref{lem_gamma_bounds_SOS-BM}, there exists a constant $K >0 $ such that, for every $t>0 $, 
			\begin{align}
				\sup_{y\in \IR; s\le t}  \cubraces{ |\gamma^{\params{\rho}{\beta}}_n [g_n](s,y)| } & \le 
				K  \sqrt{t}  m^{(1)}_{\params{0}{\beta}}(g_n) \braces{1 + \log(n)},\\
				\sup_{y\in \IR; s\le t} \cubraces{ |\gamma^{\params{\rho}{\beta}}_n[g_n](s,x)| } &\le  K  \sqrt{t} m_{\params{\sqrt n \rho}{\beta}}(|g_n|) \sqrt{n },\\
				\sup_{y\in \IR; s\le t} \cubraces{ |\gamma^{\params{\rho}{\beta}}_n[g^{2}_n](s,x)| } &\le K  \sqrt{t} m_{\params{\sqrt n \rho}{\beta}}(g^{2}_n) \sqrt{n }.
			\end{align}
			Therefore,
			\begin{equation}
				\Esp_x \braces{
					\abs{
						\frac{1}{\sqrt n} \sum_{i=1}^{[nt]}g_n(\sqrt n \hfprocess{X}{\params{\rho}{\beta}}{i-1})
					}^2
				} \le K  \sqrt{t} B_n[g_n](x),
			\end{equation}
			which converges to $0$ as $n\to \infty $.
		\end{proof}

		\begin{lemma}
			\label{lem_whgn_sqtisfies_condition_38}
			The sequence of functions $(\wh g_n)_n $ defined in~\eqref{eq_ghn_first_def} satisfies Condition~\ref{cond_B_vanish}. Moreover, 
			\begin{equation}
				\lim_{n\to\infty}\wh g_n(0)=0.
			\end{equation}
		\end{lemma}
		
		\begin{proof}
			In this proof, with abuse of notation, we make vary the probability space $\mc P_y$, $y\in \IR$, to underline when the process we are considering starts at $y$. By~\cite[p.~277]{RogWilV2} or \cite[above Proposition~2.2]{Touhami2021}, for all $n$, there exists a skew BM $X^{\params{0}{\beta}}$ defined on an extension of the probability space such that 
			$X^{\params{\rho}{\beta}} = \process{X^{\params{0}{\beta}}_{\gamma_n(t)}}$, where the time change $\gamma_n$ is the right-inverse of $t \mapsto t + c \rho\sqrt n \, \loct{X^{\params{0}{\beta}}}{0}{t}$, with $c$ a non-negative constant and it is strictly increasing, a.s.~continuous, and 
			$0 < \gamma_{n}(t) \le t $, for all $t\ge 0 $. 
			By the above we have that 
			\begin{equation}
				0 \le \wh g_{n}(y) = \Esp_{y} \braces{|X^{\params{\rho \sqrt n}{\beta}}_{1}| - |y|} 
				= \Esp_{y} \braces{|X^{\params{0}{\beta}}_{\gamma_n(1)}| - |y|}.
			\end{equation}
			Recall the scale function $s$ of $X^{\params{0}{\beta}} $, defined in~\eqref{eq_sm_SOSBM}. 
			We observe that $s(X^{\params{0}{\beta}}) $ is a martingale and $z \mapsto |s^{-1}(z)| $ is convex.  
			This qualifies 
			the process $|X^{\params{0}{\beta}}|= |s^{-1}\circ s(X^{\params{0}{\beta}})|$ as a submartingale. 
			The fact that $s(X^{\params{0}{\beta}})$ is a martingale can be proven directly by applying the Itô--Tanaka formula (Lemma~\ref{lem:ito-tanaka}) or by reminding the properties of the scale function.  
			The fact that $|X^{\params{0}{\beta}}|$ is a submartingale, and since $\gamma_n(1)\leq 1$ is a bounded stopping time, the optional sampling theorem yields  
			\begin{equation}
				0 \le \wh g_{n}(y)  \le \Esp_{y} \braces{|X^{\params{0}{\beta}}_{1}| - |y|} 
				\le 
				K \Esp_{y} \braces{|B_{1}| - |y|} =: K \wh g (y).
				\label{lem_cond_gn_Gwhgn_bound}
			\end{equation}
			The latter bound can be justified with Lemma~\ref{lem_kernel_bound} or by the fact that the absolute value of a skew-BM is a reflected BM and so $\Esp_{y} \braces{|X^{\params{0}{\beta}}_{1}|}= \Esp_{y} \braces{|B_{1}|}$. 
			
			The function $\wh g$ satisfies Condition~\ref{cond_B_vanish} with $\rho=0$ in $B_n[\wh g](x)$ (see the proof of \cite[Theorem~4.1]{Jac98}). 
			To show that $(\wh g_{n})_n$ satisfies Condition~\ref{cond_B_vanish} for $\rho \neq 0$, 
			it remains to prove that 
			\begin{equation}
				\label{eq_lem_proof_whgn_convergences_pre}
				\lim_{n\to \infty}\wh g^{2}_n(0)= 0 
				\quad \text{and}\quad  
				\lim_{n\to \infty}(1+\log n)  \wh g_n(0) = 0.
			\end{equation}
			Indeed, since $|y| m_{\params{\rho \sqrt n}{\beta}} (\rd y) = |y| m_{\params{0}{\beta}} (\rd y) = |y|(1+\sgn(y)\beta) \vd y$ and $p_{\params{\rho \sqrt n}{\beta}}(1,0,y)= v_{\rho \sqrt n}(1,0,y)$ by the expressions of these quantities~\eqref{eq_def_SkSBM_kernel_factorization}, we have 
			\begin{align}
				\wh g_n(0) 
				&=  \Esp_{0} \braces{|X^{\params{\rho \sqrt n}{\beta}}_{1}|}
				=  \int_{\IR} |y| p_{\params{\rho \sqrt n}{\beta}}(1,0,y) \, m_{\params{\rho \sqrt n}{\beta}} (\rd y)
				\\ &\le 2 \int_{\IR} |y| v_{\rho \sqrt n}(1,0,y) \vd y
				\\ &= \frac{2}{\rho \sqrt n} \int_{\IR} |y| 
				e^{2|y|/\rho \sqrt n + 2/\rho^2 n} \erfc \braces{\frac{|y|}{\sqrt{2}} + \frac{\sqrt{2}}{\rho\sqrt n}}
				\vd y.
			\end{align}
			Therefore, $\braces{(1+\log n)  \wh g_n(0)} \to 0 $ as $n\to \infty $.
			Similarly,
			\begin{align}
				\left(\wh g_n(0)\right)^2 
				\le{}&  \Esp_{0} \braces{|X^{\params{\rho \sqrt n}{\beta}}_{1}|^{2}}
				\le 2 \int_{\IR} |y|^{2} v_{\rho \sqrt n}(1,0,y) \vd y
				\\ ={}& \frac{2}{\rho \sqrt n}   
				\int_{\IR} |y|^{2} 
				e^{2|y|/\rho \sqrt n + 2/\rho^2 n} \erfc \braces{\frac{|y|}{\sqrt{2}} + \frac{\sqrt{2}}{\rho\sqrt n}}
				\vd y.
			\end{align}
			Therefore, $(\wh g_n(0))^2 \to 0$ as $n\to \infty $ and  Condition~\ref{cond_B_vanish} is satisfied. 
			This finishes the proof.
		\end{proof}
		
		We are now ready to address the main proof of this section.
		
		\begin{proof}[Proof of Proposition \ref{thm_the_case_a_05}]
			Let $(\widehat{g}_n)_{n \in \IN}, (h_{n})_{n\in \IN} $ be the sequences of functions defined for all
			$y\in \IR $, $n\in \IN $ by 
			\begin{equation}
				\label{eq_proof_linearization_wrt_loct}
				\begin{aligned}
					h_{n}(y) &= g_{n}(y) - \frac{m_{\params{\sqrt n\rho}{\beta}}(g_n)}{m_{\params{\sqrt n\rho}{\beta}}(\widehat{g}_n)} \widehat{g}_n(y),\qquad 
					&
					\wh g_{n}(y) &= \Esp_{y} \braces{|X^{\params{\rho \sqrt n}{\beta}}_{1}| - |y|}.
				\end{aligned}
			\end{equation}
			Note that $\wh g_n$ was already defined in~\eqref{eq_ghn_first_def}.
			We now show that the sequence $(h_n)_n $ satisfies Condition~\ref{cond_B_vanish}.
            Set
\begin{equation}
	c_n:=
	\frac{m_{\params{\sqrt n\rho}{\beta}}(g_n)}
	     {m_{\params{\sqrt n\rho}{\beta}}(\widehat g_n)},
	\qquad
	h_n=g_n-c_n\widehat g_n.
\end{equation}
By Proposition~\ref{prop_mngn_convergence}, $(c_n)_n$ is bounded and,
for all sufficiently large $n$,
$
	|c_n|
	\leq
	K\,m_{\params{\sqrt n\rho}{\beta}}(|g_n|).
$ 
Moreover, the fact that $g_n$ satisfies Condition~\ref{cond_B_vanish} and the boundedness of
$m_{\params{\sqrt n\rho}{\beta}}(g_n)$ imply 
\begin{equation}
	\frac{1+\log n}{\sqrt n}
	m_{\params{\sqrt n\rho}{\beta}}(|g_n|)
	\convergence{}0.
\end{equation}
Using these bounds, Lemma~\ref{lem_whgn_sqtisfies_condition_38},
and the definition of $B_n$, a direct expansion of the four terms of
$B_n[h_n](x)$ shows that $B_n[h_n](x)\to 0$. 
Thus $(h_n)_n$ satisfies Condition~\ref{cond_B_vanish}.
            
			Also, for all $n\in \IN$, $m_{\params{\sqrt n\rho}{\beta}}(h_n)=0 $.
			Therefore, by Lemma \ref{lem_L2_funct_convergence}, 
			\begin{equation}
				\Esp_x \braces{
					\abs{
						\frac{1}{\sqrt n} \sum_{i=1}^{[nt]}h_n(\sqrt n \hfprocess{X}{\params{\rho}{\beta}}{i-1})
					}^2
				}
				\convergence{}  0, \quad \text{for all } t>0.
			\end{equation}
			Convergence in probability follows, i.e.,
			\begin{equation}
				\label{eq_mod1_inprobability}
				\frac{1}{\sqrt n} \sum_{i=1}^{[nt]}h_n(\sqrt n \hfprocess{X}{\params{\rho}{\beta}}{i-1})
				\xlra{\Prob_x}{n\to \infty} 0, 
				\quad \text{for all }t>0.
			\end{equation}
			This, combined with Propositions~\ref{prop_first_loct_approximation} and~\ref{prop_mngn_convergence}, proves \eqref{eq_the_case_a_05_plain}.
			
			It remains to prove that the convergence is $\Prob_x$-ucp under the additional assumption that $\sup_n \cubraces{m_{\params{\sqrt n\rho}{\beta}}(|g_n|)} <\infty$.
			If $(g_{n})_n $ are all positive, the processes $\frac{1}{\sqrt n} \sum_{i=1}^{[nt]} g_n(\sqrt n  \hfprocess{X}{\params{\rho}{\beta}}{i-1}) $ are non-decreasing and have a continuous limit, $\Prob_x$-almost surely. 
			Thus, by Lemma~\ref{lem_ucp_convergence_condition}, the convergence is locally uniform in time, in probability (\ucp{}).
            More precisely, for an arbitrary $(g_n)_n$ satisfying the conditions of Proposition~\ref{thm_the_case_a_05}, let $g_n=g_n^+-g_n^-$, where $g_n^+(y)=\max\{g_n(y),0\}$ and $g_n^-(y)=\max\{-g_n(y),0\}$. 
            Under the additional assumption $\sup_n m_{\params{\sqrt n\rho}{\beta}}(|g_n|)<\infty$,
            the two sequences $m_{\params{\sqrt n\rho}{\beta}}(g_n^+)$ and $m_{\params{\sqrt n\rho}{\beta}}(g_n^-)$ are bounded. Hence, from any subsequence, one may extract a further subsequence along which
            $
	           m_{\params{\sqrt n\rho}{\beta}}(g_n^+)\to M_+$ and $m_{\params{\sqrt n\rho}{\beta}}(g_n^-)\to M_-$.
            Then
			\begin{align}
				&\frac{1}{\sqrt n} \sum_{i=1}^{[n\cdot]} g^{+}_n(\sqrt n \hfprocess{X}{\params{\rho}{\beta}}{i-1}) \convergence{\Prob_x\text{-}\ucp} M_+ \loct{X^{(\rho,\beta)}}{0}{}
				\\
				\text{and}\qquad &\frac{1}{\sqrt n} \sum_{i=1}^{[n\cdot]} g^{-}_n(\sqrt n \hfprocess{X}{\params{\rho}{\beta}}{i-1}) \convergence{\Prob_x\text{-}\ucp} M_- \loct{X^{(\rho,\beta)}}{0}{}.
			\end{align}

            Since 
            $ m_{\params{\sqrt n\rho}{\beta}}(g_n)
	=
	m_{\params{\sqrt n\rho}{\beta}}(g_n^+)
	-
	m_{\params{\sqrt n\rho}{\beta}}(g_n^-)
	\to M,
    $
    we have $M_+-M_-=M$.
    Along this further subsequence, applying the positive case and Lemma~\ref{lem_ucp_convergence_condition} to $g_n^+$ and $g_n^-$ yields the ucp convergence of their difference to $M\loct{X^{(\rho,\beta)}}{0}{}$.
    The conclusion for the full sequence follows from the subsequence criterion for convergence in probability, applied to the supremum norm on every compact time interval.
			This completes the proof. 
		\end{proof}

        Now that we have proven Proposition \ref{thm_the_case_a_05}, we can provide an additional local time approximation result, which is not directly applied in this paper, but it is interesting on its own and the simple proof already contains steps used in subsequent proofs, like in Proposition~\ref{lem_result_superrates} and in the one of Lemma~\ref{lem_aggregate_bound}. 
		
				\begin{lemma}\label{thm_the_case_a_01}
						Let $(g_{n})_{n}$ be a sequence of functions such that 
						\begin{equation}
								\lim_{n\to \infty} m_{\params{\sqrt n \rho}{\beta}}(|g_n|) + g_n(\sqrt{n} x)/\sqrt{n}= 0.
							\end{equation}
						Then, as $n\to \infty $
						\begin{equation}\label{eq_the_case_a_01}
								\frac{1}{\sqrt n} \sum_{i=1}^{[n\cdot]}g_n(\sqrt n \hfprocess{X}{\params{ \rho}{\beta}}{i-1}) \to 0,
							\end{equation}
						locally uniformly in time, in $L^1(\Prob_x)$.
					\end{lemma}

				\begin{proof}
						We observe that 
						\begin{equation}
								\Esp_x \braces{
										\sup_{s\in [0,t]} \abs{
												\sum_{i=1}^{[n s]}g_n(\sqrt n \hfprocess{X}{\params{ \rho}{\beta}}{i-1})
											}
									}
								\le
								|g_n(x \sqrt{n})| + \gamma^{\params{\rho}{\beta}}_n[|g_n|](x,t),
							\end{equation}
		  			where $\gamma^{\params{\rho}{\beta}}_n$ is the functional in~\eqref{eq:gamma_semigroup}. 
						Thus, from Item~\ref{eq_def_gamma_bound_A_os_SOS-BM} in Lemma~\ref{lem_gamma_bounds_SOS-BM}, 
						\begin{equation}
								\Esp_x \braces{
										\sup_{s\in [0,t]} \abs{
												\frac{1}{\sqrt n}\sum_{i=1}^{[n s]}g_n(\sqrt n \hfprocess{X}{\params{ \rho}{\beta}}{i-1})
											}
									}
								\le \frac{|g_n(x \sqrt{n})|}{\sqrt n} +  K m_{\params{ \sqrt n \rho}{\beta}}(|g_n|) \sqrt{ t},
							\end{equation}
						which vanishes as $n\to \infty $. This completes the proof.
					\end{proof}

		\subsection{A case that covers the supercritical regime}\label{ssec_un_case2}
		
		The statement in this case is the following.
		
		\begin{proposition} \label{lem_result_special}
			The convergence~\eqref{eq_thm_limit_gnT_inline} holds under the additional assumption that
			\begin{equation}
				\log n/u_n \convergence{} 0.
			\end{equation}
		\end{proposition}
		
		The proof of this mode reduces to applying Proposition~\ref{thm_the_case_a_05} after a proper re-scaling of the test function. 
		More precisely, we first prove Proposition~\ref{lem_result_superrates} which allows to pass from considering $\un=\sqrt{n}$ (for which Proposition~\ref{thm_the_case_a_05} holds) to general $\un$. 
		The proof of Proposition~\ref{lem_result_superrates} is based on showing that the asymptotic behavior of the same statistics associated with a sequence $(g_n)_n$ is similar to that of the statistics with $g_n \indic{(-r,r)}$, $r>0$, instead of $g_n$. 
		Next, in the proof of Proposition~\ref{lem_result_special}, we show that the sequence $(g_n[T])_n$ satisfies the assumptions of Proposition~\ref{lem_result_superrates}. 
		
		Note that in the following results and their proofs we keep using the notation in~\eqref{eq: measure_integral_notation} both with the speed measure and the Lebesgue measure, the latter denoted by $\lambda$.

		\begin{proposition}\label{lem_result_superrates}
			Let $(g_n)_n$ be a sequence of functions and $(u_n)_{n} $ a sequence of real numbers such that  
			\begin{enumerate}
                \item if $\rho>0$, $g_n(0)=0$ for all $n\in \IN$,
				\item
				as $n\to \infty $:  $\un \to \infty $, $\log n/u_n \to 0$, and $u_n/n \to 0$,
				\item $(g_n)_n $ a sequence of functions such that
				\begin{align}
					\sup \lambda(|g_n|)&<\infty ,
					&
					\sup_{n,y} |g_n(y)| &< \infty,
					\\
					\lim_{n\to \infty} m_{\params{0}{\beta}}(g_n)&=: M,
					&
					\lim_q \limsup_n \int_{|y|>q}|g_n(y)|\vd y &= 0.
				\end{align}
			\end{enumerate}
			Then,
			\begin{equation}
				\label{eq: supercritical_convergence}
				\frac{u_n}{n} \sum_{i=1}^{[n\cdot]}g_n(u_n \hfprocess{X}{\params{\rho}{\beta}}{i-1})
				\xlra{\Prob_x\text{-ucp}}{n\ra \infty}  M \loct{X^{\params{\rho}{\beta}}}{0}{}.
			\end{equation}
		\end{proposition}

		\begin{proof}
			Let $(k_n)_n$ be the sequence of functions defined for all
			$y\in \IR $ by $k_n(y) = \frac{u_n}{\sqrt n} g_n \braces{ \frac{u_n}{\sqrt n} y } $.
			Assume first that $\sup_n \lambda^{(1)}(g_n) < \infty$.
			Then, for some $K>0 $,
			\begin{align}\label{eq_def_transp_measures}
				& |k_n(x)| \le K \frac{\un}{\sqrt n},\qquad 
				m_{\params{0}{\beta}}(k_n) = m_{\params{0}{\beta}}(g_n), \qquad
				m_{\params{0}{\beta}}(|k_n|) = m_{\params{0}{\beta}}(|g_n|) \le K, 
				\\
				& m^{(1)}_{\params{0}{\beta}}(k_n) = \frac{\sqrt n}{\un} m^{(1)}_{\params{0}{\beta}}(g_n)
				\le \frac{\sqrt n}{\un}K, \qquad
				\lambda(k^{2}_n)= \frac{\un}{\sqrt n} \lambda(g^{2}_n)
				\le \frac{\un}{\sqrt n}K.
			\end{align}	
			From the hypotheses made on $ (\un)_n$ and $(g_n)_n $,
			we have also that $(k_n)_n $ satisfies Condition~\ref{cond_B_vanish}.
			Thus, from Proposition~\ref{thm_the_case_a_05},
			\begin{equation}\label{eq_proof_avar_the_case_a_05}
				\frac{1}{\sqrt n} \sum_{i=1}^{[n\cdot]}k_n(\sqrt n \hfprocess{X}{\params{\rho}{\beta}}{i-1}) 
				\xlra{\Prob_x\text{-ucp}}{n\ra \infty}
				M \loct{X^{\params{\rho}{\beta}}}{0}{},
			\end{equation}
			which is exactly~\eqref{eq: supercritical_convergence}. 
			This completes the proof under the condition $\sup_n \lambda^{(1)}(g_n) < \infty$.
			
			Consider now $\sup_n \lambda^{(1)}(g_n) = \infty$. 
			For all $r \geq 1$ consider $h_{n,r}:=g_n \indic{(-r,r)}$. 
			Then, $m^{(1)}_{\params{0}{\beta}}(h_{n,r}) \leq r m_{\params{0}{\beta}}(|g_n|) <\infty$.
			Moreover, for every $r$, up to considering a subsequence, 
			$ m_{\params{0}{\beta}}(h_{n,r}) \to \alpha_r$, 
			where $\alpha_r \to M$ as $r\to\infty$.
			Hence, from \eqref{eq_proof_avar_the_case_a_05}, 
			\begin{equation}
				\frac{\un}{n} \sum_{i=1}^{[n\cdot]}h_{n,r}(u_n \hfprocess{X}{\params{\rho}{\beta}}{i-1}) 
				\xlra{\Prob_x\text{-ucp}}{n\ra \infty}
				\alpha_r \loct{X^{\params{\rho}{\beta}}}{0}{}.
			\end{equation}
			Since $\alpha_r \to M$ as $r \to \infty$, it remains to show that
			\begin{equation}\label{eq_proof_diff_gn_hnr_A}
				\lim_{r\to \infty}\limsup_{n\to\infty} \Esp_x \braces{
					\sup_{s\in [0,t]}
					\abs{
						\frac{u_n}{n} \sum_{i=1}^{[ns]} g_n( u_n \hfprocess{X}{\params{\rho}{\beta}}{i-1}) 
						- \frac{u_n}{n} \sum_{i=1}^{[ns]} h_{n,r}( u_n \hfprocess{X}{\params{\rho}{\beta}}{i-1})
					}
				}
				=  0.
			\end{equation}
			If $\overline h_{n,r}(y) = g_n (y)- h_{n,r}(y) = g_n(y) \indicB{|y|>r} $,
			then the expectation in \eqref{eq_proof_diff_gn_hnr_A} is bounded by
			\begin{equation}
				\frac{\un}{n}
				\sum_{i=1}^{[nt]} \Esp_x \braces{ \abs{\overline h_{n,r}(\un \hfprocess{X}{\params{\rho}{\beta}}{i-1})}} 
				\le
				\frac{u_n}{n} \overline h_{n,r} (u_n x) +
				\frac{1}{\sqrt n} \gamma^{\params{\rho}{\beta}}_n [\wt h_{n,r}](x,t),
			\end{equation}
			where $\wt h_{n,r}(y) := \frac{\un}{\sqrt n} \overline h_{n,r} (\frac{\un}{\sqrt n} y) = |k_n|(y) \indicB{|y|>r \sqrt{n}/u_n}$
			and $\gamma^{\params{\rho}{\beta}}_n[\tilde h_{n,r}](y,t) $ defined in \eqref{eq:gamma_semigroup}.
			Since $|k_n| $ is bounded by $K \un/\sqrt n$, 
			\begin{equation}
				\lim_{n\to \infty} \frac{u_n}{n} \overline h_{n,r} (u_n x) = 0.
			\end{equation}
			It remains to check that  
			\[\limsup_{r\to\infty} \limsup_{n\to \infty} \frac1{\sqrt n} \gamma^{\params{\rho}{\beta}}_n[\tilde h_{n,r}](x,t)  =0.\]
			Indeed, since $\wt h_{n,r}(0)=0 $, Lemma~\ref{lem_gamma_bounds_SOS-BM}(i) yields, for some constant $K>0$,
			\begin{equation}
				\frac1{\sqrt n} \gamma^{\params{\rho}{\beta}}_n[\tilde h_{n,r}](x,t)  \le K  m_{\params{0}{\beta}}(|\tilde h_{n,r}|) \sqrt{t},
			\end{equation} 
			which converges to 0 as $r \to \infty$.
		\end{proof}

		\begin{proof}[Proof of Proposition~\ref{lem_result_special}]
			We observe that
			\begin{align}
				\sup_{n,y} \abs{ g_n[T](y)} &\le \xnorm{g}{\infty},
				&
				\lambda(\abs{g_n[T]}) \le \lambda(|g|)\xnorm{\tfrac{1}{T'}}{\infty} \le \tfrac{\lambda(|g|)}{\varepsilon},\\
				\lim_{n\to \infty} m_{\params{0}{\beta}}(g_n[T]) &=
				m_{\params{0}{\beta}}(g).
			\end{align}
			Also, using the change of variable $z=u_nT(y/u_n)$ and the fact that $T'\geq\varepsilon$, 
            \begin{equation}
	           \int_{|y|>q}|g_n[T](y)|\,\vd y \leq	\frac{1}{\varepsilon} 	\int_{|z|>\varepsilon q}|g(z)|\,\vd z,
            \end{equation}
            which converges to $0$ as $q\to\infty$, uniformly in $n$.
			Moreover, if $\rho>0$, then $g_n[T](0)=g(0)=0$. Since $m_{\params{0}{\beta}}$ is equivalent to $\lambda$ with a bounded density, 
			the conditions of Proposition~\ref{lem_result_superrates} are satisfied for $(g_n)_n=(g_n[T])_n $.
			The proof is thus completed.
		\end{proof}

		\subsection{The subcritical regime}
		\label{ssec_un_case3}
		
		The object of this section is to prove Proposition~\ref{prop_result_un2}, which establishes convergence in the subcritical regime. 
		In particular, this includes sequences $u_n = o(\log n)$ that were left out from our previous analysis.
		
		\begin{proposition}\label{prop_result_un2}
			The convergence~\eqref{eq_thm_limit_gnT_inline} holds under the additional assumption that $(\un^2)_n $ satisfies~\eqref{eq_intro_un_condition}, that is: 
			\begin{equation}
				\label{eq_txt_un_quadratic}
				\un^{2}/n \convergence{} 0
				\quad \text{and} \quad 
				\un \convergence{} \infty.
			\end{equation}
		\end{proposition}
		
		The strategy of our proof involves multiple steps that are inspired by the ones of \cite[Theorem~3.1]{Anagnostakis2022} for sticky BM (that was inspired by~\cite[Lemma~4.5]{Jac98} for standard BM), as the semigroups of sticky BM and SkS-BM scale similarly in time, space, and stickiness parameter.
		The main particularity of SkS-BM is the discontinuous character of the local time field in space at $0$, caused by the skew-oscillating behavior of the process. Therefore, whenever we use continuity arguments (as in Lemma~\ref{lem_loct_continuity}), we need to handle the behavior on each half-line separately. Moreover, our detailed proof completes and revises some steps in \cite{Anagnostakis2022}.
		
		Since the proof has a rather complicated structure and relies on auxiliary lemmas, we give a short step-by-step guide of its mechanism that highlights the role of each intermediary result:
		\begin{enumerate}
			\item Consider the following difference between our approximating sum and an integral with respect to the Lebesgue measure:
			\begin{equation}
				\label{eq:integral_continuous_limit}
				\frac{\un}{n} \sum_{i=1}^{[nt]} g_n[T](\un \hfprocess{X}{(\rho,\beta)}{i-1})
				- 
				\int_{\IR} g_n[T] (y) \rloct{X^{\params{\rho}{\beta}}}{y/\un }{[nt]/n} \vd y.
			\end{equation}
			
			\item Using the time-change representation of skew Brownian motion, the continuity of its local time field on each half line, and the relation between its right, left, and symmetric local times at $0$, we show in Lemma~\ref{lem_loct_continuity} that the second summand in~\eqref{eq:integral_continuous_limit}, i.e.~the integral, converges in probability, as $n\to\infty$, to $m_{\params{\rho}{\beta}}(g) \loct{X^{\params{\rho}{\beta}}}{0}{t}$.
			
			\item We re-express the difference~\eqref{eq:integral_continuous_limit} as
			\begin{equation}
				\label{eq:leading_term1}
				\frac{\un}{n}\sum_{i=1}^{[nt]}
				\int_{0}^{1} \braces{g_n[T](\un X^{\params{\rho}{\beta}}_{\frac{i-1}{n}})
					-g_n[T](\un X^{\params{\rho}{\beta}}_{\frac{i-1}{n}+\frac{s}{n}})} \vd s,
			\end{equation}
			plus some vanishing remainder term.
			
			\item We bound the leading term of~\eqref{eq:leading_term1} by
			\begin{equation}
				\frac{\un}{n} \sum_{i=1}^{[nt]} F_n[g](X^{\params{\rho}{\beta}}_{\frac{i-1}{n}})
			\end{equation}
			which will be introduced in Lemma~\ref{lem_dilatation_convergence}. Note that this quantity is of the same form as the statistic~\eqref{eq_thm_limit_gnT_inline}. 
			As $n$ goes to infinity, it vanishes. More precisely, Lemmas~\ref{lem_dilatation_convergence} and~\ref{lem_aggregate_bound} show that, if $g$ is Lipschitz and vanishes on an open interval around $0$, the above sum converges to $0$. By a standard density argument, this also holds for general integrable $g$ that are vanishing at $0$. Ucp convergence follows from Lemma~\ref{lem_ucp_convergence_condition}.
		\end{enumerate}
		
		For this, we use the following preliminary results.
		
		\begin{lemma}\label{lem_loct_continuity}
			Assume $X^{\params{\rho}{\beta}}$ is a SkS-BM on $\mc P_x$, $\un$, $g$, and $T$ as in Theorem~\ref{thm_main_SkSBM}. 
			The following convergence holds
			\begin{equation}
				\label{eq_lem_loct_continuity}
				\int_{\IR} g_n[T] (y) \rloct{X^{\params{\rho}{\beta}}}{y/\un }{[nt]/n} \vd y
				\convergence{\Prob_x}
				m_{\params{\rho}{\beta}}(g) \loct{X^{\params{\rho}{\beta}}}{0}{t}.
			\end{equation}
		\end{lemma}

		\begin{proof}    
			We first prove that the right and left local times $(t,y)\mapsto \rloct{X^{\params{\rho}{\beta}}}{y}{t} $, $(t,y)\mapsto \rloct{X^{\params{\rho}{\beta}}}{y-}{t} $ are both 
			almost surely $(t,y)$-jointly continuous on each half-plane and \begin{equation}
				\label{eq_proof_loctimes_relation}
				\begin{aligned}
					\rloct{X^{\params{\rho}{\beta}}}{0}{t} &= (1+\beta) \loct{X^{\params{\rho}{\beta}}}{0}{t},\qquad 
					&
					\rloct{X^{\params{\rho}{\beta}}}{0-}{t} &= (1-\beta) \loct{X^{\params{\rho}{\beta}}}{0}{t}.
				\end{aligned}
			\end{equation}  
			
			For this we reduce to the case of skew BM. 
			Indeed, as seen in Lemma~\ref{lem_whgn_sqtisfies_condition_38}, there exists a skew BM
			$X^{\params{0}{\beta}} $ defined on an extension of the probability space such that $X^{\params{\rho}{\beta}} =  \process{X_{\gamma(t)}^{\params{0}{\beta}}}$, 
			where $\gamma $ an a.s.~continuous strictly increasing time-change. By \cite[Exercice~VI.1.27]{RevYor}, we have that  $\rloct{X^{\params{\rho}{\beta}}}{0}{t}=\rloct{X^{\params{0}{\beta}}}{0}{\gamma(t)}$ for all $t\geq 0$.
			
			By \cite[Theorem~VI.1.6]{RevYor}, the right local time field $(t,y) \ra \rloct{X^{\params{0}{\beta}}}{y}{t} $ and symmetric local time field 
			$(t,y) \ra \loct{X^{\params{0}{\beta}}}{y}{t} $ 
			are both almost surely $(t,y) $-jointly continuous on each half-plane
			$(\IR_+ \times (-\infty,0)), (\IR_+ \times [0,\infty)) $.
			The relation between right and symmetric local times of $X^{\params{0}{\beta}}$ at the skew threshold $0$ is given by the identities (see e.g.~\cite[Theorem 2.1]{Salminen2019})
			\begin{equation}
				\begin{aligned}
					\rloct{X^{\params{0}{\beta}}}{0}{t} &= (1+\beta) \loct{X^{\params{0}{\beta}}}{0}{t},\qquad 
					&
					\rloct{X^{\params{0}{\beta}}}{0-}{t} &= (1-\beta) \loct{X^{\params{0}{\beta}}}{0}{t}.
				\end{aligned}
			\end{equation}
			
			Since time-change's inverse $\gamma^{-1} $ is a.s.~continuous, the first step of the proof is demonstrated. 
			
			From the a.s.~piece-wise continuity of $(t,y)\rightarrow \rloct{X^{\params{\rho}{\beta}}}{y}{t} $, for all $t\ge 0$ and $y>0 $,
			\begin{equation}\label{eq_local_time_pointwise_convergence}
				|\rloct{X^{\params{\rho}{\beta}}}{y/\un }{[nt]/n} -  \rloct{X^{\params{\rho}{\beta}}}{0}{t}| \convergence{\Prob_x\text{-}a.s.} 0
				\quad \text{and} \quad
				|\rloct{X^{\params{\rho}{\beta}}}{-y/\un }{[nt]/n} -  \rloct{X^{\params{\rho}{\beta}}}{0-}{t}| \convergence{\Prob_x\text{-}a.s.} 0.
			\end{equation} 
			
			Let us fix $t>0$ and an event $\Omega_0$ of full probability such that the  convergences~\eqref{eq_local_time_pointwise_convergence} hold for all $\omega \in \Omega_0$. 
			We now prove that, on $\Omega_0$,
			\[
			\int_{0}^{\infty} g_n[T](y)  \rloct{X^{\params{\rho}{\beta}}}{y/\un }{[nt]/n} \vd y \convergence{} \rloct{X^{\params{\rho}{\beta}}}{0}{t}  \int_{0}^{\infty} g(y)  \vd y  
			\]
			and 
			\[
			\int_{-\infty}^0 g_n[T](y)  \rloct{X^{\params{\rho}{\beta}}}{y/\un }{[nt]/n} \vd y \convergence{} \rloct{X^{\params{\rho}{\beta}}}{0-}{t}  \int_{-\infty}^0 g(y)  \vd y. 
			\]
			We focus on the first convergence. The second convergence can be proven with similar arguments. 
			We observe that
			\begin{align}
				\int_{0}^{\infty} g_n[T](y)  \rloct{X^{\params{\rho}{\beta}}}{y/\un }{[nt]/n} \vd y &- \rloct{X^{\params{\rho}{\beta}}}{0}{t}  \int_{0}^{\infty} g(y)  \vd y
				\\
				\qquad={} &  
				\int_{0}^{\infty} g_n[T](y)   \braces{\rloct{X^{\params{\rho}{\beta}}}{y/\un }{[nt]/n}
					- \rloct{X^{\params{\rho}{\beta}}}{0}{t}} \vd y 
				\\ & + \rloct{X^{\params{\rho}{\beta}}}{0}{t}   
				\int_{0}^{\infty} \braces{g_n[T](y) - g(y)}  \vd y.
				\label{eq_integral_loctime_continuity_additive}
			\end{align}
			
			Regarding the second additive term of the right-hand-side of~\eqref{eq_integral_loctime_continuity_additive},
			note that
			\begin{equation}
				\label{eq_int_gTn_scaling}
				\int_0^{+\infty} g_n[T](y) \vd y 
				= \int_0^{+\infty} g(y) \frac1{T'(T^{-1}(y/\un))} \vd y.
			\end{equation}
            We have operated a change of variable $z=\un T(y/\un)$ and used that $\lim_{y\to+\infty}T(y)=+\infty$. 
			And, since $T' $ is lower and upper bounded by a positive constant, $T'(0)=1 $, and $g$ is bounded, by dominated convergence we have
			\begin{align}
				\abs{\int_{0}^{\infty} \braces{g_n[T](y) - g(y)}  \vd y}
				&= 
				\abs{\int_{0}^{\infty} g(y) \sqbraces{\frac1{T'(T^{-1}(y/\un))}-1}  \vd y}
				\\ & \leq \int_0^{+\infty} |g(y)| \abs{\frac1{T'(T^{-1}(y/\un))}-1} \vd y \convergence{} 0.
			\end{align}
			
			We now consider the first additive term of the right-hand-side of~\eqref{eq_integral_loctime_continuity_additive}.
			By the non-explosive nature of the process, it holds a.s.~that for every $s>0$ 
			\begin{equation}
				\label{eq_proof_process_finiteness}
				-\infty < \inf_{u\le s} X_u < \sup_{u\le s} X_u < \infty.
			\end{equation} 
			Observe that $\rloct{X^{\params{\rho}{\beta}}}{y/\un}{[nt]/n} \leq \rloct{X^{\params{\rho}{\beta}}}{y/\un}{t} =0$ for all $y > \un \sup_{u\leq t} X^{(\rho,\beta)}_u$, and by continuity in space $K:=\sup_{y\in [0,+\infty)} \rloct{X^{\params{\rho}{\beta}}}{y/\un}{t}= \sup_{y\in [0, \sup_{u\leq t} X^{(\rho,\beta)}_u]} \rloct{X^{\params{\rho}{\beta}}}{y}{t}<\infty$. 
			Note that $K$ does not depend on $n$.
			For all $y\geq 0$, 
			\begin{align}
				& \abs{g_n[T](y)   \braces{\rloct{X^{\params{\rho}{\beta}}}{y/\un }{[nt]/n}
						- \rloct{X^{\params{\rho}{\beta}}}{0}{t}}} \\
				& \leq \abs{g_n[T](y)}  \braces{\rloct{X^{\params{\rho}{\beta}}}{y/\un }{[nt]/n}
					+  \rloct{X^{\params{\rho}{\beta}}}{0}{t}} 
				\leq \braces{K + \rloct{X^{\params{\rho}{\beta}}}{0}{t}} \abs{g_n[T](y)}.
				\label{eq: integrant_bound}
			\end{align}
			By~\eqref{eq_int_gTn_scaling} and lower boundedness of $|T'|$, the above expression is integrable, and 
			$\int_0^\infty \abs{g_n[T](y)} \vd y \leq \xnorm{1/T'}{\infty} \int_0^\infty \abs{g(y)} \vd y<+\infty$.
			Moreover, for every $y\in (0,+\infty)$ we have 
			\begin{equation} 
				\abs{g_n[T](y)   \braces{\rloct{X^{\params{\rho}{\beta}}}{y/\un }{[nt]/n} - \rloct{X^{\params{\rho}{\beta}}}{0}{t}}} 
				\leq \xnorm{g}{\infty} \abs{\rloct{X^{\params{\rho}{\beta}}}{y/\un }{[nt]/n} - \rloct{X^{\params{\rho}{\beta}}}{0}{t}},
			\end{equation}
			which converges to 0 by joint continuity of the local time. 
			We now decompose the term as
			\begin{equation}
				\int_{0}^{R} g_n[T](y)   \braces{\rloct{X^{\params{\rho}{\beta}}}{y/\un }{[nt]/n}- \rloct{X^{\params{\rho}{\beta}}}{0}{t}} \vd y
				+ \int_{R}^{\infty} g_n[T](y)   \braces{\rloct{X^{\params{\rho}{\beta}}}{y/\un }{[nt]/n}- \rloct{X^{\params{\rho}{\beta}}}{0}{t}} \vd y.
			\end{equation}
			By dominated convergence the first term vanishes as $n\to \infty $.

            Regarding the second term, using the bound~\eqref{eq: integrant_bound} and the same change of variables as in~\eqref{eq_int_gTn_scaling}, we have 
            \begin{equation}
        	   \int_R^\infty |g_n[T](y)|\,\vd y 
               = 	\int_{u_nT(R/u_n)}^\infty |g(y)| \frac{1}{T'(T^{-1}(y/u_n))}\,\vd y.
            \end{equation}
            Since $u_nT(R/u_n)\to R$ and $T'(T^{-1}(y/u_n))\to T'(0)=1$, dominated convergence yields
            \begin{equation}
        	   \int_R^\infty |g_n[T](y)|\,\vd y \convergence{} \int_R^\infty |g(y)|\,\vd y.
            \end{equation}
        Therefore,
        \begin{equation}
	           \limsup_{n\to\infty} \abs{\int_R^\infty g_n[T](y) \braces{\rloct{X^{\params{\rho}{\beta}}}{y/u_n}{[nt]/n} -\rloct{X^{\params{\rho}{\beta}}}{0}{t} } \,\vd y }
	           \leq
	           \braces{K+\rloct{X^{\params{\rho}{\beta}}}{0}{t}} \int_R^\infty |g(y)|\,\vd y.
        \end{equation}
        Since $g$ is integrable, letting $R\to\infty$, and combining this estimate with the convergence of the integral over $[0,R]$, gives 
        \begin{equation}
            \int_0^\infty g_n[T](y)
	           \braces{\rloct{X^{\params{\rho}{\beta}}}{y/u_n}{[nt]/n} -\rloct{X^{\params{\rho}{\beta}}}{0}{t}} \,\vd y \xrightarrow[n\to\infty]{\Prob_x\text{-a.s.}} 0.
        \end{equation}

			An identical argument applied to the interval $(-\infty, 0)$ yields
			\begin{equation}
				\int_{-\infty}^{0} g_n[T](y) \rloct{X^{\params{\rho}{\beta}}}{y/\un}{[nt]/n} \vd y
				\convergence{\Prob_{x}\text{-}a.s.}
				\braces{\int_{-\infty}^{0} g(y) \vd y} \rloct{X^{\params{\rho}{\beta}}}{0-}{t}.
			\end{equation}
			Adding these two results together, we conclude that the $\Prob_x $-almost sure
			limit of the sequence of integrals $n \mapsto \int_{\IR} g_n[T] (y) \rloct{X^{\params{\rho}{\beta}}}{y/\un }{[nt]/n} \vd y $,
			as $n\to \infty $ is
			\begin{equation}
				\braces{\int_{-\infty}^{0} g(y)\vd y} \rloct{X^{\params{\rho}{\beta}}}{0-}{t}
				+
				\braces{\int_{0}^{\infty} g(y)\vd y} \rloct{X^{\params{\rho}{\beta}}}{0}{t}.
			\end{equation}
			By the relation between right and symmetric local times~\eqref{eq_proof_loctimes_relation}, the proof is complete. 
		\end{proof}

		\begin{lemma}
			\label{lem_dilatation_convergence}
			Let $k$ be a bounded Lipschitz function with compact support that vanishes on an
			open interval around $0$, $(\un)_n $ be a sequence that satisfies~\eqref{eq_txt_un_quadratic}, $T$ a function that satisfies~\eqref{eq_condition_T} and $(k_n[T])_n,(F_{n}[k])_n  $ be the sequences of functions defined for all $x\in \IR $, $n\in \IN $ as
			\begin{align}\label{eq_semigroup_action_intime_def}
				k_n[T](x) := k\braces{\un T(x/\un)},
				\quad
				F_{n}[k](x) := 
				\int_{0}^{1} 
				\braces{P_{s \un^{2}/n}^{\params{\un \rho}{\beta}} k_n[T](x) - k_n[T](x)}
				\vd s.
			\end{align}
			The sequence $\braces{m_{\params{\un \rho}{\beta}}\braces{F_{n}[k]}}_n $ vanishes
			as $n\to \infty $.
		\end{lemma}

		\begin{proof}
			For all $\delta >0 $, let $U_{\delta} = \{x\in \IR: \; |x| \in (\delta,1/\delta)\} $.
			We choose $\delta>0 $ such that for all $ x \not \in U_{\delta}$, 
			$k(x)=0 $. 
			We observe that there exists some $n_0 \in \IN $ such that, for all $n\ge n_0 $
			and $ |x| \not \in U_{\delta/2}$,
			$k_n[T](x)=0 $.
			
			In what follows, we consider separately the absolutely continuous and singular parts of 
			$m_{\params{\un \rho}{\beta}} $.
			By  Fubini's theorem,
			\begin{align}
				\abs{ m_{\params{ \cdot}{\beta}}\braces{F_{n}[k]} } &= 
				\left| \int_{\IR} \braces{ \int_{0}^{1} 
					\braces{P_{s \un^{2}/n}^{\params{\un \rho}{\beta}} k_n[T](x) - k_n[T](x)}
					\vd s }  m_{\params{ \cdot}{\beta}}(\rd x) \right|
				\\ & \quad \le 
				\int_{0}^{1}
				\braces{ \int_{\IR}  
					\abs{P_{s \un^{2}/n}^{\params{\un \rho}{\beta}} k_n[T](x) - k_n[T](x)}
					m_{\params{ \cdot}{\beta}}(\rd x) }
				{\rd s}.
			\end{align}
			Also, 
			\begin{equation}
				|F_n[k](0)|
				\le 	\int_{0}^{1} 
				\abs{P_{s \un^{2}/n}^{\params{\un \rho}{\beta}} k_n[T](0) - k_n[T](0)}
				\vd s.
			\end{equation}
			Lemma \ref{lem_kernel_bound} yields
			\begin{align}
				& \abs{P_{s \un^{2}/n}^{\params{\un \rho}{\beta}} k_n[T](x) - k_n[T](x)}
				\\ & \qquad \le
				\int_{\IR}\abs{k_n[T](y) - k_n[T](x)} p_{\params{\un \rho}{\beta}}(s \un^{2}/n,x,y)  m_{\params{\un \rho}{\beta}}(\rd y)
				\\ & \qquad \le
				K 
				\int_{\IR}\abs{k_n[T](y) - k_n[T](x)} u_1(s \un^{2}/n,x,y) \vd y
				\\ & \qquad \quad+ K \abs{k_n[T](x)} 
				u_1 (s \un^{2}/n,x,0).
				\label{eq_proof_Lemma611_B1}
			\end{align}
			
			Regarding the second additive term of the right-hand-side of~\eqref{eq_proof_Lemma611_B1}, by the convergence $\lim_{n\to \infty} u^{2}_n/n = 0 $ and boundedness of $k $ (and of $k_n $), it follows that
			\begin{equation}
				\lim_{n\ra \infty} K \abs{k_n[T](x)} 
				u_1 (s \un^{2}/n,x,0) = 0, \quad \text{for all } x\in \IR 
				\text{ and } s\in [0,1].
			\end{equation}
			It also follows that there exists some $n_1 \ge n_0 $
			such that for all $n\ge n_1 $ and $s\in [0,1] $, 
			\begin{align}
				K |k_n[T](x)| u_1(s \un^{2}/n,x,0) 
				\le{} & K \xnorm{k}{\infty}    u_1(s \un^{2}/n,x,0) \indic{U_{\delta/2}}
				\\ \le{} & \xnorm{k}{\infty} \indic{U_{\delta/2}}(x) \sup_{\substack{ n \ge n_1 \\ x\in U_{\delta/2}}} u_1(s \un^{2}/n,x,0)
				\\ \le{} & K \xnorm{k}{\infty} u_1( u_{n_1}^{2}/n_1,0,\delta/2) \indic{U_{\delta/2}}(x),
			\end{align}
			which is $L^{1}(\IR) $.
			Thus, from dominated convergence,
			\begin{equation}
				\lim_{n\ra \infty} \int_{\IR} K |k_n[T](x)| u_1(s \un^{2}/n,x,0) m_{\params{ \rho}{\beta}}(\rd x) =0.
			\end{equation}

			Regarding the second additive term, if $Z$ is a Gaussian $\mathcal{N}(0,1) $ random variable, by assumptions on $k$, for every $n\ge n_0 $ and $s\in [0,1] $, 
			\begin{align}
				\int_{\IR} & \int_{\IR} \abs{k_n[T](y) - k_n[T](x)} u_1(s \un^{2}/n,x,y) \vd y \, m_{\params{ \cdot}{\beta}}(\rd x)
				\\ &=
				\Esp \left( \int_{\IR} \left| k_n[T]\braces{x + \tfrac{\un \sqrt{s}}{\sqrt{n}} Z} - k_n[T](x) \right|\, m_{\params{ \cdot}{\beta}}(\rd x) \right)
				\\ &=
				\Esp \braces{ \int_{\IR} \braces{\indic{U_{\delta/2}}(x) + \indic{U_{\delta/2}}\braces{x + \tfrac{\un \sqrt s}{\sqrt n} Z}} 
					\abs{ k_n[T]\braces{x + \tfrac{\un \sqrt{s}}{\sqrt{n}} Z} - k_n[T](x) }\, m_{\params{ \cdot}{\beta}}(\rd x) }
				\\ &\le 
				\xnorm{k}{\Lip} \xnorm{T}{\Lip}
				\frac{\un}{\sqrt{n}} \Esp \left( |Z| \int_{\IR} \braces{\indic{U_{\delta/2}}(x) + \indic{U_{\delta/2}}\braces{x + \tfrac{\un \sqrt s}{\sqrt n} Z}}\, m_{\params{ \cdot}{\beta}}(\rd x)  \right) 
				\\ &\le 
				\xnorm{k}{\Lip} \xnorm{T}{\Lip}
				\frac{\un}{\sqrt{n}} \frac{8}{\delta} \Esp \braces{ |Z| },  
			\end{align}
			which, since $E(|Z|)<\infty $, converges to $0$ as $n\to \infty $.
			
			Since the bounds do not depend on the variable $s$, from the dominated convergence theorem on $L^{1}([0,1]) $, we obtain
			\begin{equation}
				m_{\params{ \cdot}{\beta}}\braces{F_{n}[k]}  \convergence{} 0.
			\end{equation}
			
			Regarding the singular part of the measures, from~\eqref{eq_def_SkSBM_kernel_factorization},
			\begin{align}
				\rho \un \abs{P_{s \un^{2}/n}^{\params{\un \rho}{\beta}} k_n[T](0) - k_n[T](0)}
				&\le \rho \un K \int_{\IR}\abs{k_n[T](y) - k_n[T](0)} v_{\un \rho}(s \un^{2}/n,0,y) \vd y 
				\\ &= \rho \un K \int_{U_{\delta/2}}\abs{k_n[T](y)} v_{\un \rho}(s \un^{2}/n,0,y) \vd y .
			\end{align}
			Thus, by symmetry of $y\mapsto v_{\rho}(t,0,y) $, and the fact that $\erfc(y) \leq K_{\text{Mills}} e^{-y^2}/y$ (see \cite[p.~98]{GriSti}), 
			\begin{align}
				& \rho\,\un \int_{U_{\delta/2}}\abs{k_n[T](y)}  v_{\un \rho}(s \un^{2}/n,0,y) \vd y
				\\ & \qquad \le 2 \xnorm{k}{\infty} \int_{\delta/2}^{2/\delta} \exp\braces{2\frac{|y|}{\un \rho} + 2\frac{\un^{2} s}{n\un^{2}\rho^2}} \erfc \braces{\frac{|y|}{\sqrt{2s \un^{2}/n}} + \frac{\sqrt{2s\un^{2}/n}}{\un \rho}} \vd y
				\\ & \qquad \le 2 \xnorm{k}{\infty} K_{\Mills}
				\int_{\delta/2}^{2/\delta} 
				\frac{u_n \rho  \sqrt{2s/n}}{|y|\rho + 2 s \un /n}
				e^{- \frac{y^{2}n}{2 s \un^{2}}} \vd y
				\\
				& \qquad \le 4 \xnorm{k}{\infty} K_{\Mills} \frac{\rho \sqrt{2 u^{2}_n/n}}{\delta\rho} \left(\frac{2}{\delta}-\frac{\delta}2\right),
			\end{align}
			which converges to $0$ as $n\to \infty $.
			
			Again, since the bounds do not depend on $s$, from Lebesgue's dominated convergence theorem on $L^{1}([0,1]) $, 
			\begin{equation}
				\lim_{n\ra \infty} \braces{\rho \un} F_n[k] (0) = 0.
			\end{equation} 
			This finishes the proof.
		\end{proof}

		\begin{lemma}\label{lem_aggregate_bound}
			Let $X^{\params{\rho}{\beta}} $ be the $(\rho,\beta)$-SkS-BM on $\mc P_x$. 
			Then, there exists a positive constant $K>0$, such that for all $t>0$ and integrable function $g$, 
			\begin{equation}
				\Esp_{x}\braces{\sup_{s\le t} 
					\abs{\sum_{i=1}^{[n s]}g(\un \hfprocess{X}{\params{\rho}{\beta}}{i-1})}}
				\le 
				|g(\un x)| +  
				K \braces{
					\rho |g(0)| 
					+ 
					\frac{1}{\un} m_{\params{0}{\beta}}(|g|)}
				n \sqrt{2t}.
			\end{equation}
		\end{lemma}
		
		\begin{proof}
			We observe that 
			\begin{equation}
				\int_{\IR} |g(\un y)| m_{\params{\rho}{\beta}}(\rd y)
				= \rho |g(0)| + \frac{1}{\un}\int_{\IR}|g(y)| m_{\params{0}{\beta}}(\rd y) 
				= \rho |g(0)| + \frac{1}{\un} m_{\params{0}{\beta}}(|g|).
			\end{equation}
			From Lemma~\ref{prop_semig_bound} and the fact that $ \sum_{i=1}^{[nt]} i^{-1/2} \le 2\sqrt{nt}$, there exists a positive constant $K>0$ that does not depend on $(t,g)$ such that 
			\begin{align}
				\Esp_{x}\braces{\sup_{s\le t} 
					\abs{\sum_{i=1}^{[n s]}g(\un \hfprocess{X}{\params{\rho}{\beta}}{i-1})}}
				&\le  
				\sum_{i=1}^{[nt]}
				\Esp_{x}\braces{
					\abs{g(\un \hfprocess{X}{\params{\rho}{\beta}}{i-1})}}
				\\
				& = |g(\un x)| + \sum_{i=2}^{[nt]} P^{\params{\rho}{\beta}}_{\frac{i-1}n} |g|(x)
				\\
				&\le |g(\un x)| +  
				K \braces{
					\rho |g(0)| 
					+ 
					\frac{1}{\un} m_{\params{0}{\beta}}(|g|)}
				\sum_{i=2}^{[nt]} \frac{\sqrt n}{\sqrt i}
				\\
				&\le  |g(\un x)| +  
				K \braces{
					\rho |g(0)| 
					+ 
					\frac{1}{\un} m_{\params{0}{\beta}}(|g|)}
				n\sqrt{2t}.
			\end{align}
			This finishes the proof.
		\end{proof}
		
		We are now ready to address the proof of Proposition~\ref{prop_result_un2}.
		
		\begin{proof}
			[Proof of Proposition~\ref{prop_result_un2}]
			For simplicity, in this proof we fix $t>0$. 
			We observe that $X^{\params{\rho}{\beta}} $ is a semimartingale with quadratic variation 
			$\vd \qv{X^{\params{\rho}{\beta}}}_s = \indicB{X^{\params{\rho}{\beta}}_s \not = 0} \vd s $.
			By the occupation times formula 
			\cite[Corollary~VI.1.6]{RevYor},
			for any non-negative measurable function $f$, it holds that
			\begin{equation}
				\int_{0}^{t} f(X^{\params{\rho}{\beta}}_s) \indicB{X^{\params{\rho}{\beta}}_s \ne 0} \vd s  =
				\int_{\IR} f(y) \indicB{y \not = 0} \rloct{X^{\params{\rho}{\beta}}}{y}{t} \vd y.
			\end{equation}	
			In particular, if $\rho=0$ or $f(0)=0 $, then the formula becomes
			\begin{equation}\label{eq_occ_tim_formula_special}
				\int_{0}^{t} f(X^{\params{\rho}{\beta}}_s) \vd s =  \int_{\IR} f(y) \rloct{X^{\params{\rho}{\beta}}}{y}{t}
				\vd y. 
			\end{equation}
			Applying consecutive change of variables and \eqref{eq_occ_tim_formula_special} yields 
			\begin{align}\label{eq_gn_occupation_time_equality}
				&\frac{\un}{n} \int_{0}^{[nt]} g_n[T](\un X^{\params{\rho}{\beta}}_{s/n}) \vd s
				= \un \int_{0}^{[nt]/n} g_n[T](\un X^{\params{\rho}{\beta}}_s) \vd s\\
				&\qquad = \un \int_{\IR} g_n[T]( \un y) \rloct{X^{\params{\rho}{\beta}}}{y}{[nt]/n}\vd y  
				=
				\int_{\IR} g_n[T](y) \rloct{X^{\params{\rho}{\beta}}}{ y/ \un}{[nt]/n} \vd y .
			\end{align}
			
			And, from Lemma~\ref{lem_loct_continuity}, we have
			\begin{equation}
				\label{eq_loct_decomposition_term3}
				\int_{\IR} g_n[T] (x) \rloct{X^{\params{\rho}{\beta}}}{x/\un}{[nt]/n}\vd x 
				\xlra{\Prob_x}{n\ra \infty}
				m_{\params{\rho}{\beta}}(g) \loct{X^{\params{\rho}{\beta}}}{0}{t}.
			\end{equation}
			
			By these assertions, to complete the proof, it suffices to show that 
			\begin{equation}\label{eq_vn_decomposition}
				\frac{\un}{n} \sum_{i=1}^{[nt]} g_n[T] (\un \hfprocess{X}{\params{\rho}{\beta}}{i-1})
				-
				\frac{\un}{n} \int_{0}^{[nt]} g_n[T](\un X^{\params{\rho}{\beta}}_{s/n}) \vd s
				\xlra{\Prob_x}{n\ra \infty} 0.
			\end{equation}
			Indeed, once~\eqref{eq_vn_decomposition} is proven, we obtain
			\begin{equation}
				\frac{\un}{n} \sum_{i=1}^{[nt]} g_n[T] (\un \hfprocess{X}{\params{\rho}{\beta}}{i-1}) 
				\xlra{\Prob_x}{n\to \infty}
				m_{\params{\rho}{\beta}}(g) \loct{X^{\params{\rho}{\beta}}}{0}{t}.
			\end{equation}
			With the same arguments, based on Lemma~\ref{lem_ucp_convergence_condition}, as in the closure of the proof of Proposition~\ref{thm_the_case_a_05}, the $\Prob_x$-ucp convergence is proven.
			
			Let us now prove~\eqref{eq_vn_decomposition}. 
			First, let us note that  
			\begin{align}
				\frac{\un}{n} \sum_{i=1}^{[nt]} g_n[T] & (\un \hfprocess{X}{\params{\rho}{\beta}}{i-1})
				-
				\frac{\un}{n} \int_{0}^{[nt]} g_n[T](\un X^{\params{\rho}{\beta}}_{s/n}) \vd s
				\\
				= & 
				\frac{\un}{n}\sum_{i=1}^{[nt]}
				\int_{0}^{1} \braces{g_n[T](\un X^{\params{\rho}{\beta}}_{\frac{i-1}{n}})
					-g_n[T](\un X^{\params{\rho}{\beta}}_{\frac{i-1}{n}+\frac{s}{n}})} \vd s 
				\\ &+
				\frac{\un}{n} \int^{nt-[nt]}_{0} \braces{g_n[T](\un X^{\params{\rho}{\beta}}_{\frac{[nt]}{n}})
					-g_n[T](\un X^{\params{\rho}{\beta}}_{\frac{[nt]}{n}+\frac{s}{n}})} \vd s.
				\label{eq_vn_decomposition:bis}
			\end{align}
			The second additive term of the right-hand-side of~\eqref{eq_vn_decomposition:bis}
			is bounded by $2 \xnorm{g}{\infty} \un/n $, therefore
			\begin{equation}\label{eq_loct_decomposition_term2}
				\frac{\un}{n} \int^{nt-[nt]}_{0} \braces{g_n[T](X^{\params{\rho}{\beta}}_{\frac{[nt]}{n}})
					-g_n[T](X^{\params{\rho}{\beta}}_{\frac{[nt]}{n}+\frac{s}{n}})} \vd s \xlra{n\to\infty}{a.s.} 0 .
			\end{equation}
			By~\cite[Lemma~9]{GCJ93}, and since $g$ is bounded, to show that the first additive term of the right-hand-side of~\eqref{eq_vn_decomposition:bis} converges to 0 in probability, it suffices to prove that
			\begin{equation}
				D_t^{(n)} := \frac{\un}{n}\sum_{i=1}^{[nt]}
				\int_{0}^{1} \Esp_x \braces{g_n[T](X^{\params{\rho}{\beta}}_{\frac{i-1}{n}})
					-g_n[T](X^{\params{\rho}{\beta}}_{\frac{i-1}{n}+\frac{s}{n}})
					\Big| 
					\bF_{\frac{i-1}{n}}
				} \vd s 
				\convergence{\Prob_x}
				0. 
			\end{equation}
			Using the notation of Lemma~\ref{lem_dilatation_convergence},
			$D_t^{(n)}$ reads 
			\begin{equation}
				D_t^{(n)}
				=  \frac{\un}{n} \sum_{i=1}^{[nt]} F_n[g](X^{\params{\rho}{\beta}}_{\frac{i-1}{n}}).
			\end{equation}
			By Lemma~\ref{lem_aggregate_bound}, we have the bound
			\begin{equation}
				\Esp_x \braces{ | D_t^{(n)} | } 
				\le
				\frac{\un}{n} 
				|F_n[g](\un x)| +  
				\frac{\un}{n} 
				\braces{
					\rho |F_n[g](0)| 
					+ 
					\frac{1}{\un} m_{\params{0}{\beta}}(|F_n[g]|)}
				n \sqrt{2t},
			\end{equation}
			where $|F_n[g](\un x)| \le 2\xnorm{g}{\infty}$.
			
			First assume $g$ is Lipschitz continuous, with compact support and vanishes on an open interval around $0$. Then, by Lemma~\ref{lem_dilatation_convergence}, the upper bound of the above inequality vanishes as $n\ra \infty $. 
			And, in this case, 
			\begin{equation}\label{eq_proof_Bn_convergence}
				D_t^{(n)} \xlra{L^{1}(\Prob_x)}{n\to\infty} 0,
			\end{equation}
			which implies convergence in probability.
			
			Now assume the general form of $g$. 
			We now approximate $g$ by functions which satisfy the assumptions of Lemma~\ref{lem_dilatation_convergence}.    
			Since $g$ is bounded and integrable,
			for every $p>0 $, there exists a Lipschitz-continuous function with compact support $k^{(p)} $ that vanishes in the vicinity of $0$ such that
			$m_{\params{ \rho}{\beta}}(|g-k^{(p)}|)< 1/p $ (see the proof of
			\cite[Lemma 4.5]{Jac98}).
			Moreover, by positive Fubini's theorem, triangle inequality, the fact that $g(0)=k^{(p)}(0)=0$ when $\rho>0$, and the symmetry of the SkS kernel $p_{(\rho,\beta)}(t,x,y)= p_{(\rho,\beta)}(t,y,x)$ (see \cite[p.149]{ItoMcKean96}), we have that
			\begin{multline}
				\label{eq_integral_semigroup_bound_proof}
				m_{\params{ \rho}{\beta}}(|P^{\params{\un \rho}{\beta}}_t g-P^{\params{\un \rho}{\beta}}_t k^{(p)}|)
				\\ \le \int_{\IR} \int_{\IR} \abs{g(y) - k^{(p)}(y)} p_{(\un \rho,\beta)}(t,x,y)\, m_{(\un \rho,\beta)}(\rd y)\,
				m_{(\rho,\beta)}(\rd x)
				\\ 
				\le  \int_{\IR}  \abs{g(y) - k^{(p)}(y)} \braces{\int_{\IR} p_{(\un \rho,\beta)}(t,y,x)
					\, m_{(\un \rho,\beta)}(\rd x)} \, m_{(\rho,\beta)}(\rd y)
				\\
				= \int_{\IR} \abs{g(y) - k^{(p)}(y)} \, m_{(\rho,\beta)}(\rd y)
				= m_{\params{ \rho}{\beta}}(|g-k^{(p)}|) < 1/p.
			\end{multline}
			Similarly, since $T$ was assumed to satisfy $\varepsilon \le T'(y) \le 1/\varepsilon$, for all $y\in \IR $, by a suitable change of variables we obtain for all $t>0$,
			$m_{\params{ \rho}{\beta}}(|g_n[T]-k^{(p)}_n[T]|) < \xnorm{1/T'}{\infty}\braces{1/p} \le 1/(p\varepsilon) $.
			This, combined with similar arguments as in~\eqref{eq_integral_semigroup_bound_proof}, yields that
			\begin{equation}
				\label{eq_semigroup_integral_bound}
				m_{\params{ \rho}{\beta}}(|P^{\params{\un \rho}{\beta}}_t g_n[T]-P^{\params{\un \rho}{\beta}}_t k_n^{(p)}[T]|) < 1/(p\varepsilon).
			\end{equation}
			By the above, Fubini's theorem ensures
			\begin{align}
				m_{\params{ \rho}{\beta}}(|F_n[g]-F_n[k^{(p)}]|)
				\le{}& \int_{0}^{1} m_{\params{ \rho}{\beta}}(
				|P_{s \un^{2}/n}^{\params{\un \rho}{\beta}} g_n[T]
				-
				P_{s \un^{2}/n}^{\params{\un \rho}{\beta}} k^{(p)}_n[T]|)
				\vd s
				\\ &+ \int_{0}^{1} m_{\params{ \rho}{\beta}}(|g_n[T]-k^{(p)}_n[T]|)
				\vd s
				\le 2/(p\varepsilon).
			\end{align}
			Also, from \eqref{eq_def_SkSBM_kernel_factorization}--\eqref{eq_def_SkSBM_kernel_factors}, for $\rho>0$
			\begin{align}
				(\rho \un) p_{\params{\un \rho}{\beta}}(t \un^{2}/n,0,y) = 
				(\rho \un) v_{\un \rho}(t \un^{2}/n,0,y)  
				=  
				e^{\frac{2|y|}{\un\rho} + \frac{2t}{n\rho^2}} \erfc \braces{\frac{|y|\sqrt n}{\sqrt{2t}\un} + \frac{\sqrt{2t}}{\sqrt n \rho}}.
			\end{align}
            For $\rho=0$, the quantity on the left-hand side is identically zero. 
			Therefore
			\begin{equation}
				\lim_{n\to \infty} (\rho \un) p_{\params{\un \rho}{\beta}}(t \un^{2}/n,0,y)
				= 
				\begin{cases}
					0, & y \not = 0 \text{ or } \rho=0,\\
					1, & y = 0 \text{ and } \rho > 0.
				\end{cases}
			\end{equation}
            If $\rho>0$, combining the above, by dominated convergence
			\begin{align}
				&\limsup_{n\to \infty} \un |F_n[g](0)-F_n[k^{(p)}](0)|
				\\
				&\le \limsup_{n\to \infty} \int_{0}^{1} \un
				\abs{ 
					P_{s \un^{2}/n}^{\params{\un \rho}{\beta}}g_n[T](0)
					- P_{s \un^{2}/n}^{\params{\un \rho}{\beta}}k^{(p)}_n[T](0)
				} \vd s
				\\
				& \le
				\limsup_{n\to \infty}  
				\int_{0}^{1} \int_{\IR} 
				\abs{ 
					g_n[T](y) - k^{(p)}_n[T](y)
				} \un p_{\params{\un\rho}{\beta}}(s \un^{2}/n,0,y) \, m_{\params{\un\rho}{\beta}}(\!\vd y) \vd s 
				\\ & \le (1/\rho) m_{\params{\rho}{\beta}}(|g_n[T]-k^{(p)}_{n}[T]|) \leq 1/(\rho p \varepsilon).
				\label{eq_proof_distance_Fn_kp_g}
			\end{align}
            If $\rho=0$, no estimate of $F_n[g](0)-F_n[k^{(p)}](0)$ is needed, since the corresponding term in Lemma~\ref{lem_aggregate_bound} is multiplied by $\rho$. 

			By Lemma~\ref{lem_aggregate_bound}, and the above bounds, 
			\begin{equation}
				\Esp_x \braces{ | D_t^{(n)} | } 
				\le 
				\frac{\un}{n} 
				|F_n[g](\un x)| 
				+ \sqrt{2t} \braces{
					\un \rho |F_n[k^{(p)}](0)|
					+  m_{\params{0}{\beta}}(|F_n[k^{(p)}]|)
					+ \frac{C}{ p \varepsilon}
					+ \frac{2}{p\varepsilon} }.
			\end{equation}
			Thus, by Lemma~\ref{lem_dilatation_convergence}, there exists a constant $C'$ such that for all $p>0$, we have \begin{equation}\label{eq_loct_decomposition_term1}
				\limsup_{n\to \infty} \Esp_x \braces{ | D_t^{(n)} | } \le C'/p
				\quad	
				\Longrightarrow
				\quad
				\lim_{n\to \infty} \Esp_x \braces{ | D_t^{(n)} | } = 0,
			\end{equation}
			the implication is due to arbitrariness of $p>0 $. Thus, we have proven~\eqref{eq_vn_decomposition} and this completes the proof. 
		\end{proof}
		

		\section{Proof of the remaining local time approximation results}
		\label{sec_reduction}

		We first prove Theorem~\ref{thm_main_SOSBM}\ref{thm_item_Refl}, which completes the picture for SkS-BM by covering the limiting case of sticky-reflected BM.
		Next, we use Theorem~\ref{thm_main_SkSBM} and Theorem~\ref{thm_main_SOSBM}\ref{thm_item_Refl} to prove the analogous result for SOS-BM, Theorem~\ref{thm_main_SOSBM}\ref{thm_item_SOS}.
		And finally, we prove Theorem~\ref{thm_main}. 
		The ``reduction'' from SOS-diffusions to SkS-BM with $\beta\in [-1,1]$ relies on the transformations $T_1$ and $T_2$ defined in~\eqref{eq_T1_T2_def}.
		The transformation $T_2 $ (along with local equivalence of laws)  reduces the problem from an SOS-diffusion to an SOS-BM. The transformation $T_1 $ reduces the problem from an SOS-BM to an SkS-BM.
		
		The reduction to SkS-BM, and in particular to non-reflected SkS-BM (i.e.~$\beta\in (-1,1)$), is chosen for analytical convenience. 
		While equivalent results (to the ones derived in Sections~\ref{sec_preliminary} and  \ref{sec_SkSBM}, and in Appendix~\ref{app_proofsSemigroup}) 
		can be obtained directly for SOS-BM (reflected or not), this would require more intricate computations.

		\subsection{From non-reflected SkS-BM to sticky-reflected BM}
		
		Relying on the already proven Theorem~\ref{thm_main_SkSBM}, we now show the result also holds for sticky-reflected BM (i.e.~$|\beta|=1$), and so we prove Theorem~\ref{thm_main_SOSBM}\ref{thm_item_Refl}.  
		Without loss of generality, we consider only the case of positive reflection (i.e.~$\beta=1$). 
		
		Given any non-reflecting SkS-Brownian motion $Y$ of parameters $(\rho,\beta) $, the process $X:=|Y|$ is a sticky-reflected BM starting at $|Y_0|$ with the same stickiness parameter $\rho $. 
		Given functions $(g, T)$ as in the assumptions of Theorem~\ref{thm_main_SkSBM}, let $\wt g(y):=g(|y|)$ and $\wt T(y) := \sgn(y) T(\sgn(y) y)$,
		$y\in \IR $. Observe that $\wt T(y) = T(|y|) $ if $y>0 $, and $\wt T = - T(|y|) $ if $y\le 0 $. This yields the identity 
		\begin{align}
			\frac{\un}{n} \sum_{i=1}^{[n\cdot]} g_n[T](\un \hfprocess{X}{}{i-1}) 
			= \frac{\un}{n} \sum_{i=1}^{[n\cdot]} \wt g_n[\wt T](\un \hfprocess{Y}{}{i-1}).
		\end{align}
		
		Since $g$ and $T$ satisfy the hypothesis of Theorem~\ref{thm_main_SkSBM} 
		($g $ bounded, integrable, vanishing at $0$, and $T$ piecewise twice differentiable so that~\eqref{eq_condition_T} holds), so do the $\wt g$ and $\wt T$. 
		By Theorem~\ref{thm_main_SkSBM} for the process $Y$, the right hand side of the above identity converges in ucp to $\left(\int_{\IR} \wt g \vd m^Y\right) \loct{Y}{0}{}$, where since we consider symmetric local times $\loct{Y}{0}{}=\loct{|Y|}{0}{}=\loct{X}{0}{}$ and 
		\begin{equation}
			\int_{\IR} \wt g \vd m^Y=\int_{-\infty}^{+\infty} g(|x|) (1+\sgn(x)\beta) \vd x= \int_{-\infty}^{+\infty} g(|x|) \vd x = 2 \int_{0}^{+\infty} g(x) \vd x = \int_{\IR} g \vd m^X,
		\end{equation}
		with $m^Y(\!\vd x):=(1+\sgn(x)\beta) \vd x +\rho \delta_0(\!\vd x)$ and $m_0(\!\vd x)=m^X(\!\vd x):= 2 \indic{(0,+\infty)}(x) \vd x + \rho \delta_{0}(\!\vd x)$ (see~\eqref{eq_sm_SOSBM}).
		Therefore, the conclusion of 
		Theorem~\ref{thm_main_SOSBM}\ref{thm_item_Refl} holds for the process $X$ and it implies that, for every $t>0$,
		\begin{equation}
			\label{eq_SRBM_Proba1}
			\frac{\un}{n} \sum_{i=1}^{[n t]} g_n[T](\un \hfprocess{X}{}{i-1}) - \left(\int_{\IR} g \vd m_{0}\right) \loct{X}{0}{t} \xlra{\law}{n\to \infty} 
			0.
		\end{equation}
		This holds when $X$ is the absolute value of the SkS-BM $Y$ as above. 
		
		Now, we deduce that the convergence holds for any sticky-reflected BM. 
		Let $X$ be a sticky-reflected BM on $\mc P_x$ (recall that the space is an element of the (weak) solution). 
		By uniqueness in law, $X$ has the same law as the sticky-reflected BM constructed above from a non-reflected SkS-BM. Hence~\eqref{eq_SRBM_Proba1} holds. 
		
		Since the limit is constant, convergence in law and convergence in probability are equivalent (see e.g.~\cite[Theorem~25.3]{billingsley1995probability}). Therefore, we have that
		\begin{equation}
			\frac{\un}{n} \sum_{i=1}^{[n t]} g_n[T](\un \hfprocess{X}{}{i-1}) 
			- \left(\int_{\IR} g \vd m_{0}\right) \loct{X}{0}{t} \xlra{\Prob_x}{n\to \infty} 0,
		\end{equation}
		for all $t>0$. 
		Arguing as in the final step of the proof of Proposition~\ref{thm_the_case_a_05}, we decompose $g=g^+-g^-$ and apply Lemma~\ref{lem_ucp_convergence_condition} to the two corresponding non-decreasing processes. This yields the desired $\Prob_x$-ucp convergence.
        This completes the proof for sticky-reflected BM. 
		\qed
		
		\subsection{From SkS-BM to SOS-BM}
		
		We now complete the proof of Theorem~\ref{thm_main_SOSBM} by proving its Item~\ref{thm_item_SOS}. 
		Let $X$ be an SOS-BM of parameters $(\rho,\beta,\sigma_-,\sigma_+)$ and let $T_1$ the Lamperti transform given by~\eqref{eq_T1_T2_def}.
		
		\begin{lemma}
			\label{prop_reduction_SOSBM}
			The process $T_1(X) $ is an SkS-BM of
			parameters $(\rho^{(1)},\beta^{(1)}) $, where 
			\begin{equation}
				\begin{aligned}
					\rho^{(1)} &:= \rho \frac{2 \sigma_-\sigma_+}{ \sigma_-(1+\beta) + \sigma_+(1-\beta)},
					& \beta^{(1)} &:= \frac{ \sigma_-(1+\beta) - \sigma_+(1-\beta)}{ \sigma_-(1+\beta) + \sigma_+(1-\beta)}.
				\end{aligned}
			\end{equation}
		\end{lemma}
		
		\begin{proof}
			The function $T_1$ is a piecewise linear function, it is difference of two convex functions, it is strictly increasing, and $T_1(0)=0$.  
			Therefore, applying It\^o--Tanaka formula (Lemma~\ref{lem:ito-tanaka}), the process $Y = T_1(X) $ solves
			\begin{equation}
				\begin{dcases}
					Y_t = Y_0 + \int_0^t \indicB{ Y_s \neq 0 } \vd W_s 
					+ \frac{(\sigma_- + \sigma_+) \beta + (\sigma_- - \sigma_+)}{(\sigma_- - \sigma_+) \beta + (\sigma_- + \sigma_+)} \loct{Y}{0}{t},
					\\
					\int_{0}^{t}\indicB{ Y_s = 0 } \vd s = \rho \frac{2 \sigma_-\sigma_+}{(\sigma_- - \sigma_+) \beta + (\sigma_- + \sigma_+)}   \loct{Y}{0}{t},
				\end{dcases}
			\end{equation}
			which proves the result.
			The fact that $T_1(0)=0$ ensures that the threshold position is preserved by the transformation.  
		\end{proof}
		
		\begin{remark}
			\begin{enumerate}
				\item Similarly, if $\beta \in (-1,1)$, one can prove that with a similar transformation $\wt T_1 := [x \mapsto x/(1+\beta \sgn(x))]$,   
				the process $\wt T_1 (X) $ is an oscillating sticky BM, i.e., an SOS-BM
				with skewness parameter $\beta=0 $, of stickiness and diffusion coefficients
				$(\rho^{(2)},\sigma^{(2)}) $ given by
				\begin{equation}
					\begin{aligned}
						\rho^{(2)} &:= \rho,
						& 
						\sigma^{(2)}(x) &:= 
						\frac{\sigma_-}{1-\beta} \indic{(-\infty,0)}(x) 
						+ \frac{\sigma_+}{1+\beta} \indic{[0,+\infty)}(x),
						&
						&\text{for all } x \in \IR.
					\end{aligned}
				\end{equation}
				
				\item For $\rho=0$, we recover the interplay between skew and oscillating BM (without stickiness), see
				\cite{Harrison1981skew,Maz19}.

                \item Since $\sigma_-,\sigma_+>0$, $\beta^{(1)}=1$ (resp.~$-1$) if and only if $\beta=1$ (resp.~$-1$).
			\end{enumerate}
		\end{remark}
		
		We are now ready to address the proof of Theorem~\ref{thm_main_SOSBM}\ref{thm_item_SOS} for SOS-BM by using Theorem~\ref{thm_main_SkSBM} and Theorem~\ref{thm_main_SOSBM}\ref{thm_item_Refl}. 
		
		\begin{proof}
			[Proof of Theorem~\ref{thm_main_SOSBM}\ref{thm_item_SOS}]
			Let $T$ be a piecewise twice differentiable function which satisfies~\eqref{eq_condition_T}, and let $T_1$ be the function in~\eqref{eq_T1_T2_def}. 
			For all $c>0 $, consider the functions $S_c$ and $\phi_c $, defined as 
			\begin{equation}
				S_c(y) = \frac{T( c  y)}{c}, \quad \phi_c(y)= c y, 
				\quad \text{for all } y \in \IR.
			\end{equation}
			Note that these functions are both sign-preserving and $S_c$ further satisfies~\eqref{eq_condition_T}. 
			
			We consider the process $Y := T_1(X) $. 
			By Lemma~\ref{prop_reduction_SOSBM}, this process is an SkS-BM 
			of parameters $(\rho^{(1)},\beta^{(1)})$ 
			and \[\loct{Y}{0}{} = \frac{(1+\beta)\sigma_-+(1-\beta)\sigma_+}{2 \sigma_+\sigma_-} \loct{X}{0}{}.\]
			Rewriting the statistics in terms of $Y$ instead of $X$, we obtain that
			\begin{align}
				\frac{\un}{n} \sum_{i=1}^{[nt]} g_n[T](\un \hfprocess{X}{}{i-1})
				&=
				\frac{\un}{n} \sum_{i=1}^{[nt]} g_n[T]\braces{\un
					\sigma_0(\hfprocess{Y}{}{i-1}) \hfprocess{Y}{}{i-1}}.
			\end{align}
            Note that $g_n[T]\braces{\un \sigma_0(x) x} = g\braces{\un T(\sigma_0(x) x)}$, but the function $x\mapsto T(\sigma_0(x) x)$ does not satisfy~\eqref{eq_condition_T}. 
			Therefore, we rewrite this terms as follows:
			\begin{align}
				&\frac{\un}{n} \sum_{i=1}^{[nt]} g_n[T](\un \hfprocess{X}{}{i-1})
				\\& \quad = \frac{\un}{n} \sum_{i=1}^{[nt]} (g_{<0} \circ \phi_{\sigma_-} )_n[S_{\sigma_-}]\braces{\un
					\hfprocess{Y}{}{i-1}}
				+
				\frac{\un}{n} \sum_{i=1}^{[nt]} (g_{>0} \circ \phi_{\sigma_+} )_n[S_{\sigma_+}]\braces{\un
					\hfprocess{Y}{}{i-1}},
			\end{align}
			where $g_{>0}= \indic{(0,\infty)} g $ and $g_{<0}= \indic{(-\infty,0)} g$.
			
			Applying Theorem~\ref{thm_main_SkSBM} if $\beta^{(1)}\in(-1,1)$, and Theorem~\ref{thm_main_SOSBM}\ref{thm_item_Refl} for sticky-reflected BM 
			if $\beta^{(1)}\in\{-1,1\}$, with the functions $S_{\sigma_+}, S_{\sigma_-}$, we get
			\begin{align}
				\frac{\un}{n} \sum_{i=1}^{[n\cdot ]} g_n[T](\un \hfprocess{X}{}{i-1})
				\xlra{\Prob_x\text{-ucp}}{n\ra \infty} \ &
				(1-\beta^{(1)}) \int g_{<0} \circ \phi_{\sigma_-} (y)\vd y
				\\ & +
				(1+\beta^{(1)}) \int g_{>0} \circ \phi_{\sigma_+}(y) \vd y\Big) \loct{Y}{0}{}.		
			\end{align}
			A change of variables in the integrals, the relation between the local times, and the expression of $\beta^{(1)}$ yield 
			\begin{equation}
				\label{eq_SOSBM_prelimit}
				\begin{aligned}
					\frac{\un}{n} \sum_{i=1}^{[n\cdot ]} g_n[T](\un \hfprocess{X}{}{i-1})
					\xlra{\Prob_x\text{-ucp}}{n\ra \infty} 
					\Big( \frac{1-\beta}{ \sigma_-^2 } \int_{\IR} g_{<0} (y) \vd y
					+  \frac{1+\beta}{\sigma_+^2} \int_{\IR} g_{>0} (y) \vd y \Big) \loct{X}{0}{}.
				\end{aligned}
			\end{equation}
			Observe that the right-hand side of~\eqref{eq_SOSBM_prelimit} is $\braces{\int_{\IR} g \vd m_0} \loct{X}{0}{}$.
			This completes the proof.
		\end{proof}
		
		\subsection{Proof of the local time approximation for SOS-diffusions, Theorem~\ref{thm_main}}
		\label{subsec_mainproof}
		
		The proof of Theorem~\ref{thm_main} relies on the local equivalence in law between the solutions to~\eqref{eq_SkSID_pathwise_char} and~\eqref{eq_SkSID_equiv}.
		In the entire proof, let $T>0$ be fixed. 
		
		Let us first prove the statement for $Y$ solution to~\eqref{eq_SkSID_equiv} on $\wt{\mc P}_x$. 
		
		First, we suppose that $\sigma $, $\sigma'$, and $1/\sigma $ are all bounded by $1/\varepsilon $ for some $\varepsilon >0 $. 
		In this case, since the function $\sigma/\sigma_0$ is continuously differentiable on both sides of $0$, the transformation $T_2$ defined in~\eqref{eq_T1_T2_def}  
		satisfies condition~\eqref{eq_condition_T}. Also, its inverse $T_2^{-1}$ satisfies  condition~\eqref{eq_condition_T}. In particular $T_2^{-1}(0)=0$, $(T_2^{-1})'(0) = 1 $, 
		\begin{equation}
			\varepsilon^{2} \le (T_2^{-1})'(y) 
			\le \frac{1}{\varepsilon^{2}},
			\quad \text{and} \quad 
			|(T_2^{-1})''(y)| 
			\le \frac{1}{\varepsilon^{4}},
		\end{equation}
		for all $y\in \IR\setminus \{0\}$, with finite one-sided limits $(T_2^{-1})''(0\pm)$. 
		
		Defining the process $Z:=T_2(Y) $ and applying the Itô--Tanaka formula, i.e.~Lemma~\ref{lem:ito-tanaka}, we see that it satisfies $\wt \Prob_x $-a.s.~the system of equations
		\begin{equation}
			\begin{dcases}
				\vd Z_t	= \sigma_0(Z_t) \indicB{Z_t \ne 0} \vd B_t   + \beta \vd\loct{Z}{0}{t}, & t\in[0,T], \\
				\indicB{Z_t = 0}	\vd t = \rho \vd \loct{Y}{0}{t} = \rho \vd \loct{Z}{0}{t}
				, & t\in[0,T].
			\end{dcases}		
		\end{equation}		
		This qualifies $Z $ as an SOS-BM of speed measure $m_0 $ (see~\eqref{eq_sm_SOSBM}). 
		
		Direct application of Theorem~\ref{thm_main_SOSBM} (with function $T_2^{-1}$ and process $Z$) 
		and the a.s.~identity $\loct{Y}{0}{} = \loct{Z}{0}{} $ yield that
		\begin{equation}
			\frac{\un}{n} \sum_{i=1}^{[n\cdot]} g(\un \hfprocess{Y}{}{i-1})
			= \frac{\un}{n} \sum_{i=1}^{[n\cdot]} g_n[T_2^{-1}](\un \hfprocess{Z}{}{i-1}) \xlra{ \wt\Prob_x\text{-ucp}}{n\ra \infty}
			\left( \int_{\IR} g \vd m_0 \right) \loct{Y}{0}{}. 
		\end{equation}
		
		For the general case ($\sigma,\sigma',1/\sigma $ unbounded),  consider an increasing sequence of compacts $(K_m)_m $ of the state space $J$ of $Y$ such that $0 \in K_m $, for all $m\in \IN $, and $\bigcup_{m} K_m  = J$. 
		\begin{itemize}
			\item In case $J=[0,r) $: let $K_m=[0,r_m]$ where $r_m= m$ if $r=\infty$; otherwise $r_m = r-\tfrac{1}{m} $.
			\item In case $J=(l,r) $: let $K_m=[l_m,r_m]$ where $l_m=-m$ if $l=-\infty$; otherwise $l_m = l+ \tfrac{1}{m} $, and $r_m$ as above.
            \item In case $J=(l,0]$: let $K_m=[l_m,0]$, where $l_m=-m$ if $l=-\infty$; otherwise $l_m=l+\tfrac1m$.
		\end{itemize}
		For each $m$, extend $\sigma$ outside $K_m$ to a coefficient $\sigma_m$ such that $\sigma_m$, $1/\sigma_m$, and $\sigma_m'$ are globally bounded, while $\sigma_m=\sigma$ on $K_m$. By uniqueness in law and localization, the law of $Y$ stopped at $\tau_m$ coincides with that of the diffusion with coefficient $\sigma_m$ stopped upon exiting $K_m$. 
        By the boundary behavior of $Y$, we have $\tau_m\to\infty$ almost surely as $m\to\infty$. Hence the previous bounded-coefficient case applies after stopping at $\tau_m$.
        More precisely for all $t\in (0,T]$, $\varepsilon >0$ $m\in \IN$,
		\begin{align}
			&\wt \Prob_x\left(\sup_{s\leq t}\Big|\frac{\un}{n}\sum_{i=1}^{[n s ]} g(\un \hfprocess{Y}{}{i-1})
			-
			\left( \int_{\IR} g \vd m_0 \right) \loct{Y}{0}{s}\Big| \geq \varepsilon\right)
			\\
			&\quad \leq 
			\wt \Prob_x\left(\sup_{s\leq t}\Big|\frac{\un}{n}\sum_{i=1}^{[n s ]} g(\un \hfprocess{Y}{}{i-1})
			-
			\left( \int_{\IR} g \vd m_0 \right) \loct{Y}{0}{s}\Big| \geq \varepsilon, \tau_m \geq t\right)
			+ \wt \Prob_x\left(\tau_m < t\right).
		\end{align}
		For every $m$, the summand of the right hand side of the latter inequality vanishes as $n\to \infty$.
		Therefore, letting $m\to \infty $ yields that 
		\begin{equation}
			\frac{\un}{n}\sum_{i=1}^{[n \cdot ] } g(\un \hfprocess{Y}{}{i-1})
			\xlra{\wt \Prob_x\text{-ucp}}{n\ra \infty}
			\left( \int_{\IR} g \vd m_0 \right) \loct{Y}{0}{ }. 
		\end{equation}
		By local equivalence of the laws of $X$ and $Y$ on $[0,T]$, the convergence in probability holds for all $t\in[0,T]$:
		\begin{equation}
			\frac{u_n}{n}\sum_{i=1}^{[nt]} g(u_n X_{\frac{i-1}{n}})
			\xrightarrow[n\to\infty]{\mathrm{P}_x}
			\left( \int_{\mathbb{R}} g \, dm_0 \right) L_t^0(X),
			\quad \text{for all } t\in[0,T].
		\end{equation}
		Since $T>0$ was arbitrary, the convergence holds for every $t>0$. 
		Arguing as in the final step of the proof of Proposition~\ref{thm_the_case_a_05}, by Lemma~\ref{lem_ucp_convergence_condition}, this implies the locally uniform convergence in time, in probability (ucp). 
		This completes the proof.
		\qed

		\section{Proofs of the estimation results, Propositions~\ref{thm_estimators_SkSBM} and~\ref{prop_SOSBM_estimation}}
		\label{sec_estimation}

		This section develops the proofs for our parameter estimation results, combining occupation time approximations with local time convergence to deduce consistent estimators for the stickiness and skewness parameters.
		
		\begin{proof}
			[Proof of Proposition~\ref{thm_estimators_SkSBM}]
			
			We consider the following statistics 
			\begin{align}
				V_{n}^{(0)}&:= \frac{1}{n} \sum_{i=1}^{[nt]} \indicB{\hfprocess{X}{}{i-1}=0} ,
				&
				V_{n}^{(+)}&:= \frac{\un}{n} \sum_{i=1}^{[nt]} g_{>0}(\un \hfprocess{X}{}{i-1}) ,
				&
				V_{n}^{(-)}&:= \frac{\un}{n} \sum_{i=1}^{[nt]} g_{<0}(\un \hfprocess{X}{}{i-1}),
			\end{align}
			and
			\begin{equation}
				W^{(+)}_n  :=\frac{ S_{n}^{g+}(X) +  S_{n}^{g-}(X)}2, 
				\quad 
				W^{(-)}_n  := \frac{ S_{n}^{g+}(X) -  S_{n}^{g-}(X)}2.
			\end{equation}
			
			By~\eqref{eq_SkSID_pathwise_char} and Lemma~\ref{lem_sosdiff_occtime_approximation}, 
			\begin{equation}
				\label{eq_proof_V0_conv}
				V_{n}^{(0)} \xlra{\Prob_x}{n\ra \infty} \rho\, \loct{X}{0}{t}.
			\end{equation}
			By Theorem \ref{thm_main},
			\begin{align}
				\frac{\un}{n} \sum_{i=1}^{[nt]} g_{>0}(\un \hfprocess{X}{}{i-1}) 
				&\xlra{\Prob_x}{n\ra \infty}
				\frac{1+\beta}{\sigma^{2}(0+)} \braces{\int_{0}^{\infty} g(x) \vd x} \loct{X}{0}{t},
				\label{eq_prop_hf_skew_estimationA}
				\\
				\frac{\un}{n} \sum_{i=1}^{[nt]} g_{<0}(\un \hfprocess{X}{}{i-1} ) 
				&\xlra{\Prob_x}{n\ra \infty}
				\frac{1-\beta}{\sigma^{2}(0-)} \braces{\int_{-\infty}^{0} g(x) \vd x} \loct{X}{0}{t}.
				\label{eq_prop_hf_skew_estimationB}
			\end{align}
			Thus, it holds that
			\begin{equation} \label{eq:conv:sn}
				S_{n}^{g+}(X) \xlra{\Prob_x}{n\ra \infty} (1+\beta)
				\loct{X}{0}{t},
				\quad \text{and} \quad
				S_{n}^{g-}(X) \xlra{\Prob_x}{n\ra \infty} (1-\beta)  \loct{X}{0}{t},  
			\end{equation} 
			and consequently that 	\begin{equation}\label{eq_proof_statistics_prestatistics_convergence}
				\begin{aligned}
					W^{(+)}_n &\xlra{\Prob_x}{n\ra \infty} \loct{X}{0}{t},
					&
					W^{(-)}_n &\xlra{\Prob_x}{n\ra \infty} \beta \loct{X}{0}{t}.
				\end{aligned}
			\end{equation}
			
			We consider the event $\mc L = \{\tau^X_0 < t\} $.
			We observe that on $(\mc L)^{c} $ all the aforementioned statistics converge to $0$. 
			On $\mc L = \{\tau^X_0 < t\} $ they converge to non-zero random quantities. 
			Thus, from \eqref{eq_proof_statistics_prestatistics_convergence},
			\begin{equation}
				\beta_n(X) = W^{(-)}_n \big/ W^{(+)}_n 
				\xlra{\Prob_x}{n\ra \infty} \beta.
			\end{equation}
			Also, from
			\eqref{eq_proof_V0_conv} and \eqref{eq_proof_statistics_prestatistics_convergence},
			\begin{equation}
				\rho_n(X) = V_{n}^{(0)} \big/ W^{(+)}_n 
				\xlra{\Prob_x}{n\ra \infty} \rho.
			\end{equation}
			This finishes the proof. 
		\end{proof}

		\begin{proof}
			[Proof of Proposition~\ref{prop_SOSBM_estimation}]
			
			Let $X$ be the SOS-BM as in the statement of Proposition~\ref{prop_SOSBM_estimation}, and let $B$ be a standard BM on $\mc P_x$ such that $(X,B)$ satisfies~\eqref{eq_SOSBM_pathwise_char}.
			From the It\^o--Tanaka formula, $X^{+} = \max\{X,0\} $ solves
			\begin{equation}
				\begin{cases}
					\vd X^{+}_t =  \sigma_0(X_t) \indicB{X_t >0} \vd B_t + \frac{\beta+1}{2} \vd \loct{X}{0}{t},\\
					\indicB{X_t =0} \vd t = \rho \vd \loct{X}{0}{t}.
				\end{cases}
			\end{equation}
			It is thus a semimartingale with explicit Doob--Meyer decomposition, and quadratic variation 
			given by
			\begin{equation}
				\begin{aligned}
					\vd \qv{X^{+}}_t &= \sigma_{0}^{2}(X_t) \indicB{X_t >0} \vd t,
					&
					\qv{X^{+}}_t &= \int_{0}^{t} \sigma_{0}^{2}(X_s) \indicB{X_s >0} \vd s 
					= \sigma^{2}_{+}  \int_{0}^{t} \indicB{X_s >0} \vd s.
				\end{aligned}
			\end{equation}
			By~\cite[Definition~I.2.3]{RevYor}, 
			\begin{equation}
				\sum_{i=1}^{[nt]} (\hfprocess{X}{+}{i}-\hfprocess{X}{+}{i-1})^{2} 
				\xlra{\Prob_x}{n\ra \infty} \qv{X^{+}}_t = \sigma^{2}_{+}  \int_{0}^{t} \indicB{X_s >0} \vd s.
			\end{equation}
			With similar arguments, one can show that
			\begin{equation}
				\sum_{i=1}^{[nt]} (\hfprocess{X}{-}{i}-\hfprocess{X}{-}{i-1})^{2} 
				\xlra{\Prob_x}{n\ra \infty} \qv{X^{-}}_t = \sigma^{2}_{-}  \int_{0}^{t} \indicB{X_s <0} \vd s.
			\end{equation}
			The above convergences and Lemma~\ref{lem_sosdiff_occtime_approximation} yield the results.
		\end{proof}

		\appendix 
		
		\section{Singular symmetric It\^o calculus}
		\label{app_singIto}

		In this section we derive explicit versions of classical stochastic calculus results
		for the semimartingale $X$ that solves~\eqref{eq_SkSID_pathwise_char}.
		Namely, the It\^o--Tanaka formula and Girsanov theorem. 
		Both results are used in Section~\ref{sec_reduction} for extending the local time approximation from SkS-BM to $X$.
		
		\begin{lemma}[It\^o--Tanaka formula for SOS systems] \label{lem:ito-tanaka}
			Let $X$ be a solution to~\eqref{eq_SkSID_pathwise_char} on $\mc P_x$.
			If $f$ is a difference of two convex functions, such that 
			$f\in C^{2}(\IR\setminus\{0\})$,
			then $ f(X) = \process{f(X_t)}$ satisfies
			\begin{equation}
				\begin{dcases}
					f(X_t) = f(X_0) + \int_0^t \left(f'(X_s) \nu_s + \frac12 f''(X_s) \sigma(X_s)^2 \right) \indicB{X_s \neq 0} \vd s
					\\ 
					\quad + \int_0^t f'(X_s) \sigma(X_s) \indicB{X_s \neq 0} \vd B_s 
					+ \frac{\Sigma f'(0) \beta + \Delta f'(0)}{2} \loct{X}{0}{t},
					\\
					\int_{0}^{t}\indicB{X_s = 0} \vd s = \rho   \loct{X}{0}{t},
				\end{dcases}
			\end{equation}
			where $\Delta f'(0) = f'(0+)-f'(0-) $ and $\Sigma f'(0) = f'(0+)+f'(0-)$. 
			If $f $ is further assumed strictly increasing and $f(0)=0$, it holds that
			\begin{equation}
				\loct{f(X)}{0}{t} = \frac{\Delta f'(0)  \beta 
					+ \Sigma f'(0)}{2}
				\loct{X}{0}{t}.
			\end{equation} 
		\end{lemma}
		
		\begin{proof}
			The process $X $ is a semimartingale as,
			\begin{equation}
				X_t = X_0 + \int_{0}^{t} \nu_s \indicB{X_s \ne 0} \vd s + \int_{0}^{t} \sigma(X_s) \indicB{X_s \ne 0} \vd B_s
				+ \beta \loct{X}{0}{t},
			\end{equation}
			where
			\begin{enumerate}
				\item $\braces{\int_{0}^{t} \nu_s \indicB{X_s \ne 0} \vd s + \beta \loct{X}{0}{t}}_{t\ge 0} $ is a process of bounded variation,
				\item $\braces{\int_{0}^{t} \sigma(X_s) \indicB{X_s \ne 0} \vd B_s}_{t\ge 0} $ is a local martingale.
			\end{enumerate}
			Thus, from \cite[Theorem VI.1.5]{RevYor},
			\begin{align}
				f(X_t) = f(X_0) + \int_0^t \frac{f'(X_s+) + f'(X_s-)}2 \vd X_s + \frac{1}{2} \int_{\IR} \loct{X}{y}{t} f''(\vd y), \label{eq_proof_appendix_SIT_1}
			\end{align}
			where $f' $ is the right derivative of $f$, $f'(x-):= \lim_{h\to 0} f'(x-h) $,
			and $f''(\vd x)$ is the measure defined such that for every $g\in C^1_c$,
			\begin{equation}
				\int_{\IR} f' g \vd x = - \int_{\IR} g f''(\vd x).
			\end{equation}
			Since $f \in C^{2}(\IR\setminus\{0\})$, $f''(\rd x)=f''(x) \indic{\IR\setminus\{0\}}(x) \vd x + \Delta f'(0) \delta_0( \rd x)$
			and 
			from the occupation times formula
			\begin{equation}\label{eq_proof_appendix_SIT_2}
				\frac{1}{2} \int_{\IR} \loct{X}{y}{t} f''(\vd y)=\frac{\Delta f'(0)}{2} \loct{X}{0}{t} + \frac12 \int_0^t f''(X_s) \indicB{X_s \neq 0} \sigma(X_s)^2 \vd s.
			\end{equation}
			From \eqref{eq_proof_appendix_SIT_1} and \eqref{eq_proof_appendix_SIT_2},
			\begin{align} 
				f(X_t) ={} & f(X_0) + \int_0^t \left(f'(X_s) \nu_s + \frac12 f''(X_s) \sigma(X_s)^2 \right) \indicB{X_s \neq 0} \vd s
				\\ &+ \int_0^t f'(X_s) \sigma(X_s) \indicB{X_s \neq 0} \vd B_s 
				+ \frac{\Sigma f'(0)  \beta 
					+ \Delta f'(0)}{2} \loct{X}{0}{t}. \label{eq_proof_IT_intB}
			\end{align}
			This proves the first statement.

			Assume now $f$ is strictly increasing and $f(0)=0$. To complete the proof, it suffices to verify the relationship between the local times.
			Applying \eqref{eq_proof_IT_intB} to the function $|f|$ yields
			\begin{align}
				|f(X_t)| ={} & |f(X_0)| + \int_0^t \sgn(f(X_s)) \left(  f'(X_s) \nu_s + \frac12 f''(X_s) \sigma(X_s)^2 \right) \indicB{X_s \neq 0} \vd s
				\\ & + \int_0^t \sgn(f(X_s)) f'(X_s) \sigma(X_s) \indicB{X_s \neq 0} \vd B_s 
				\\ & + \frac{\Delta f'(0)  \beta 
					+ \Sigma f'(0)}{2} \loct{X}{0}{t}. \label{eq_proof_appendix_SIT_3}
			\end{align}        
			The It\^o--Tanaka formula also yields
			\begin{equation}\label{eq_proof_appendix_SIT_4}
				|f(X_t)|=|f(X_0)| + \int_0^t \sgn(f(X_s)) \vd f(X_s) + \loct{f(X)}{0}{t}.
			\end{equation}
			Thus, since $\sgn(f(0))=\sgn(0)=0$,
			\begin{align}
				|f(X_t)| ={} & |f(X_0)| + \int_0^t \sgn(f(X_s)) \left(  f'(X_s) \nu_s + \frac12 f''(X_s) \sigma(X_s)^2 \right) \indicB{X_s \neq 0} \vd s
				\\ 
				&+ \int_0^t \sgn(f(X_s)) f'(X_s) \sigma(X_s) \indicB{X_s \neq 0} \vd B_s 
				+ \loct{f(X)}{0}{t}. \label{eq_proof_appendix_SIT_5}
			\end{align}     
			From \eqref{eq_proof_appendix_SIT_3} and \eqref{eq_proof_appendix_SIT_5},
			\begin{equation}
				\loct{f(X)}{0}{t} = \frac{\Delta f'(0)  \beta 
					+ \Sigma f'(0)}{2}
				\loct{X}{0}{t}.
			\end{equation}
			The proof is thus completed.
		\end{proof}

		\begin{lemma}[Girsanov theorem for SOS systems]\label{lem_sticky_girsanov}
			Let $X$ be a solution to~\eqref{eq_SkSID_pathwise_char} on
            $\mc P_x$, let $T>0$, and let $\theta$ be a progressively measurable process such that
            $\Prob_x\left(\int_0^T \theta_s^2\,\vd s<\infty\right)=1$  
			and let $\mathcal{E}(\theta) $ be the process 
			\begin{equation}
				\mathcal{E}_t(\theta) = 	\exp \braces{\int_{0}^{t} \theta_s \vd B_s - \frac{1}{2} \int_{0}^{t} \theta^2_s \vd s}, \quad t \in [0,T].
			\end{equation}
			Assume that $(\mc E_t(\theta))_{0\leq t\leq T}$ is a $\Prob_x$-martingale, and define a probability measure $\Qrob_x^T$ on $(\Omega,\bF_T)$ by 
            \[
            \frac{\vd\Qrob_x^T}{\vd\Prob_x^T}
                =
            \mc E_T(\theta),
            \]
            where $\Prob_x^T$ denotes the restriction of $\Prob_x$ to $\bF_T$.
            Then, 
            $\Qrob_x^T\sim\Prob_x^T$ and the process 
			$\widetilde{B} = \braces{ B_t - \int_{0}^{t} \theta_s \vd s  }_{t\in [0,T]}$ is a standard BM on $(\Omega,\bF_T,(\bF_t)_{t\in[0,T]},\Qrob^{T}_x)$.
			Moreover, under $\Qrob^{T}_x$ the process
			$X$ solves, for all $t\in [0,T]$,
			\begin{equation}\label{eq_sde_Q}
				\begin{dcases} 
					\vd X_t	= \braces{\nu_t + \theta_t \sigma(X_t)} \indicB{X_t \ne 0} \vd t + \sigma(X_t) \indicB{X_t \ne 0} \vd \widetilde{B}_t +
					\beta \vd \loct{X}{0}{t}, 
					\\
					\indicB{X_t = 0} \vd t = \rho \vd \loct{X}{0}{t}.
				\end{dcases}
			\end{equation}
		\end{lemma}
		
		\begin{proof}
			By Girsanov's theorem~\cite[Theorem 3.5.1]{KarShr}, the process $(\widetilde{B}_t)_{t\in [0,T]} $ is a standard BM under
			$\Qrob^{T}_x$ and the probability measures $\Prob^{T}_x$ and $\Qrob^{T}_x$ are equivalent.
			By substitution,
			\begin{equation}\label{eq_proof_girs_sde_part}
				\vd X_t	= \braces{\nu_t + \theta_t \sigma(X_t)} \indicB{X_t \ne 0} \vd t + \sigma(X_t) \indicB{X_t \ne 0} \vd \widetilde{B}_t  +
				\beta \vd \loct{X}{0}{t}.
			\end{equation}
			
			Since $\Qrob_x^T\sim\Prob_x^T$, the pathwise identities defining the symmetric local time of $X$ and the sticky relation remain valid $\Qrob_x^T$-a.s.. Hence 
            \[
                \indicB{X_t=0}\,\vd t = \rho\,\vd\loct{X}{0}{t},
            \]
            under $\Qrob_x^T$ as well.
			Thus, $X$ solves \eqref{eq_sde_Q} under $\Qrob^{T}_x$.
		\end{proof}
		
		\begin{remark}
			The transformed drift depends on $\theta_t$ only on the set
			$\{X_t\neq0\}$. Hence, when matching the drift of
			\eqref{eq_SkSID_pathwise_char} with that of
			\eqref{eq_SkSID_equiv}, the process $\theta$ is determined only
			away from the threshold. 
			
			Different progressively measurable extensions on $\{X_t=0\}$ are possible. Provided that the corresponding stochastic exponential is a true martingale, they yield the same transformed equation for $X$. If the target system is unique in law, they therefore induce the same law for $X$, although the resulting probability measures on the underlying space may differ.
		\end{remark}
		
		\section{Occupation times approximation} \label{app_occupation}
        
        \begin{proof}[Proof of Lemma~\ref{lem_sosdiff_occtime_approximation}]

		Without loss of generality, we can always assume that the interval $U$ is contained in $J$. 
		By Lemma~\ref{lem_ucp_convergence_condition}, it suffices to prove for every $t>0$ that
		\begin{equation}
			\label{eq_conv_prob}
			\frac{1}{n} \sum_{i=1}^{[nt]} \indic{U}(\hfprocess{X}{}{i-1})
			\xrightarrow[n\to \infty]{\Prob_x} \occt{X}{U}{t}. 
		\end{equation}
		Let $t>0$ be fixed. 
		
		Let us first assume $X$ is an SkS-BM on $\mc P_x$. 
		
		We need some additional notation. 
		For all $\varepsilon>0 $, let 
		\begin{equation}
			U_{\varepsilon}:=\{y \in \IR: \dist(y,U)<\varepsilon\}
			\quad \text{and} \quad 
			U_{-\varepsilon} := \IR \setminus \{ y \in \IR : \dist(y,U^{c}) \le \varepsilon\}.
		\end{equation}
		Define for every $\varepsilon>0 $ the pair of functions $(\psi_{\varepsilon},\varphi_{\varepsilon}) $ as
		\begin{align}
			\psi_{\varepsilon}(y)&= \indic{U_{-\varepsilon}}(y) + \braces{ \frac{1}{\varepsilon}\dist(y,U^c)}\indic{ U\setminus U_{-\varepsilon}}(y),\quad\forall\, x\in \IR,
			\\
			\varphi_{\varepsilon}(y)&= \indic{U}(y) + \braces{1 - \frac{1}{\varepsilon}\dist(y,U)}\indic{U_{\varepsilon}\setminus U}(y),\quad\forall\, y\in \IR.
		\end{align}
		Note that for every $\varepsilon >0$,
		\begin{equation}
			\indic{U_{-\varepsilon}}(y) \le \psi_{\varepsilon}(y) \le \indic{U}(y) \le \varphi_{\varepsilon}(y) \le \indic{U_{\varepsilon}}.
		\end{equation}

        {\it Step 1: $U$ is a closed interval with non-empty interior. } 
        Assume first that $U=[l,r]$ with $l,r\in\IR$ and $l<r$. 
        For every $\varepsilon\in(0,(r-l)/2)$, we have $U_{\varepsilon} = (l-\varepsilon,r+\varepsilon)$ and $U_{-\varepsilon}=(l+\varepsilon, r-\varepsilon) $, $U\setminus U_{-\varepsilon}=[l,l+\varepsilon) \cup (r-\varepsilon,r]$ and $U_\varepsilon\setminus U =(l-\varepsilon,l) \cup (r,r+\varepsilon)$. 
		Observe that, for $\varepsilon $ small enough, $0$ does not belong to $ \Int(U \setminus U_{-\varepsilon})\cup (U_\varepsilon \setminus U)$.  
		Since we are interested in the limit behavior, we can therefore assume that $ 0 \not\in \Int(U \setminus U_{-\varepsilon})\cup (U_\varepsilon \setminus U)$.
		
		All continuous functions on $[0,t] $ are Riemann integrable, therefore since the functions $[t \mapsto \psi_{\varepsilon}(X_t)] $ and $[t \mapsto \varphi_{\varepsilon}(X_t)] $ are a.s.~continuous, it holds a.s.~that 
		\begin{align}
			\lim_{n\to \infty} \frac{1}{n} \sum_{i=1}^{[nt]} \psi_{\varepsilon}(\hfprocess{X}{}{i-1}) 
			&= \int_{0}^{t}  \psi_{\varepsilon}(X_s) \vd s,
			\quad \text{and} \quad 
			\lim_{n\to \infty} \frac{1}{n} \sum_{i=1}^{[nt]} \varphi_{\varepsilon}(\hfprocess{X}{}{i-1}) 
			= \int_{0}^{t}  \varphi_{\varepsilon}(X_s) \vd s. \label{eq_convergence_phiepsilon}
		\end{align}
		
		From simple computations, we have that
		\begin{align}
			\limsup_{n\to \infty} & \bigg(\frac{1}{n} \sum_{i=1}^{[nt]}  \indic{U}(\hfprocess{X}{}{i-1})
			- \int_{0}^{t}\indic{U}(X_s) \vd s \bigg) 
			\\ &\le 
			\limsup_{n\to \infty} \braces{\frac{1}{n} \sum_{i=1}^{[nt]} \varphi_{\varepsilon}(\hfprocess{X}{}{i-1} )
				- \int_{0}^{t}\indic{U}(X_s) \vd s}
			\le \int_{0}^{t}\varphi_{\epsilon}(X_s)  \vd s
			- \int_{0}^{t}\indic{U}(X_s) \vd s
			\\ & \le \int_{0}^{t} \left( 1 -\frac1\varepsilon \dist(X_s,U) \right) \indic{U_{\varepsilon} \setminus U}(X_s) \vd s \le \int_{0}^{t} \indic{U_{\varepsilon} \setminus U}(X_s) \vd s.
			\label{eq:proof_occt_boundA1}
		\end{align}
		On the other hand, 
		\begin{align}
			\liminf_{n\to \infty} 
			&\bigg(\frac{1}{n} \sum_{i=1}^{[nt]} \indic{U}(\hfprocess{X}{}{i-1})
			- \int_{0}^{t}\indic{U}(X_s) \vd s \bigg)
			\\ &\ge 
			\liminf_{n\to \infty} \braces{\frac{1}{n} \sum_{i=1}^{[nt]} \psi_{\varepsilon}(\hfprocess{X}{}{i-1} )
				- \int_{0}^{t}\indic{U}(X_s) \vd s}
			\\ &\ge \int_{0}^{t}\psi_{\epsilon}(X_s)  \vd s
			- \int_{0}^{t}\indic{U}(X_s) \vd s
			\ge -\frac1\varepsilon \int_{0}^{t} \dist(X_s,U^c) \indic{U \setminus U_{-\varepsilon}}(X_s) \vd s
			\\&
			\ge - \int_{0}^{t} \indic{\Int(U \setminus U_{-\varepsilon})}(X_s) \vd s.
			\label{eq:proof_occt_boundA2}
		\end{align}
		
		Lemma~\ref{prop_semig_bound} ensures that 
		\begin{multline}
			\Esp_x \braces{ \int_0^t \indic{\Int(U \setminus U_{-\varepsilon})\cup (U_\varepsilon \setminus U)}(X_s) \vd s}
			\leq  \int_0^t \Esp_x \braces{ \indic{\Int(U \setminus U_{-\varepsilon})\cup (U_\varepsilon \setminus U)}(X_s)} \vd s 
			\\
			\leq K \lambda(\Int(U \setminus U_{-\varepsilon})\cup (U_\varepsilon \setminus U)) \int_0^t s^{-1/2} \vd s 
			\leq  4 K \varepsilon \sqrt{t},
		\end{multline}
		for some constant $K\in (0,\infty)$ that depends on $\rho,\beta$. 
		Let $\varepsilon_k=2^{-k}$. By the preceding estimate, 
        $
            \sum_{k=1}^{\infty}\Esp_x[Z_k]<\infty,
        $ 
        where $Z_k:= \int_0^t \indic{ \Int(U\setminus U_{-\varepsilon_k})  \cup (U_{\varepsilon_k}\setminus U) }(X_s)\,\vd s$. Hence $Z_k\to0$, $\Prob_x$-a.s..
		By~\eqref{eq:proof_occt_boundA1}-\eqref{eq:proof_occt_boundA2} we have that a.s.~
		\begin{align}
			0 = & \lim_{k\to\infty} \int_{0}^{t} \indic{\Int(U \setminus U_{-\varepsilon_k})}(X_s) \vd s 
			\leq \liminf_{n\to \infty} \bigg(\frac{1}{n} \sum_{i=1}^{[nt]}  \indic{U}(\hfprocess{X}{}{i-1})
			- \int_{0}^{t}\indic{U}(X_s) \vd s \bigg)  
			\\
			& \leq \limsup_{n\to \infty} \bigg(\frac{1}{n} \sum_{i=1}^{[nt]}  \indic{U}(\hfprocess{X}{}{i-1})
			- \int_{0}^{t}\indic{U}(X_s) \vd s \bigg) 
			\leq \lim_{k\to\infty} \int_{0}^{t} \indic{U_{\varepsilon_k} \setminus U}(X_s)\vd s =0.
		\end{align}
		From this, we conclude that $\Prob_{x}$-a.s.
		\begin{equation}
			\label{eq_limit_Uclosed}
			\lim_{n\to \infty} \frac{1}{n} \sum_{i=1}^{[nt]} \indic{U}(\hfprocess{X}{}{i-1})
			= \int_{0}^{t} \indic{U}(X_s) \vd s.
		\end{equation}
        The same argument applies when $U=[l,+\infty)$ or $U=(-\infty,r]$. In these cases there is only one finite endpoint, and the corresponding boundary layers have Lebesgue measure at most $2\varepsilon$. The case $U=\IR$ is immediate. 
		
		{\it Step 2: $U$ is a degenerate closed interval, that is, a singleton. } Let $U=\{a\}$.
		Consider the sets $U^{(0)} := [a-\varepsilon,a+\varepsilon] $, $U^{(1)} := [a-\varepsilon,a] $, and $U^{(2)} := [a,a+\varepsilon] $.
		We observe that $\indic{U}+\indic{U^{(0)}} = \indic{U^{(1)}}+\indic{U^{(2)}} $.
		Since~\eqref{eq_limit_Uclosed} holds for each of $U^{(0)} $, $U^{(1)} $, and $U^{(2)} $, 
		it also holds $\Prob_{x}$-a.s.
		\begin{equation}
			\lim_{n\to \infty} \frac{1}{n} \sum_{i=1}^{[nt]}  \indicB{\hfprocess{X}{}{i-1}=a} 
			= \int_{0}^{t}  \indicB{X_s = a} \vd s.
		\end{equation}
		This completes the proof for the case $U=\{a\} $.
		
		{\it Step 3: $U$ is a half-open or open interval. } 
        Any half-open or open interval differs from its closure by at most two endpoints. Since the result holds for closed intervals by the first step of the proof and for singletons by the second step, the conclusion follows.
		
		{\it Step 4: From SkS-BM to SOS-diffusions. }
		As in the proof of Theorem~\ref{thm_main}, we deduce the general case (SOS-diffusion) from SkS-BM.
		
		First let $Y$ be solution to~\eqref{eq_SkSID_equiv} with $\sigma,\sigma',{1}/{\sigma} $ bounded. The transformations $T_1 $ and $T_2 $ defined in~\eqref{eq_T1_T2_def} are all bijective, continuous, with continuous inverses and we have seen that $Z:=T_1 \circ T_2(Y)$ is an SkS-BM. Since $T_1 \circ T_2 $ is bijective, continuous, with continuous inverse, 
		\begin{itemize}
			\item if $U$ is an open (resp. closed) interval, so is $T_1 \circ T_2(U) = \left(T_1\circ T_2)^{-1}\right)^{-1}(U)$,
			\item $\frac{1}{n} \sum_{i=1}^{[nt]} \indic{U}(\hfprocess{Y}{}{i-1}) = \frac{1}{n} \sum_{i=1}^{[nt]} \indic{ T_1\circ T_2(U)}(\hfprocess{Z}{}{i-1}) $,
			\item and $O^{U}_t (Y)=O^{T_1\circ T_2(U)}_t (Z) $.
		\end{itemize}
		
		By a localization argument, as in the proof of Theorem~\ref{thm_main}, we can remove the boundedness assumption on $\sigma$, $\sigma'$, $1/\sigma$. 
		Local equivalence in law Condition~\ref{cond_drift_volatility_functions}-\ref{item_cond_nu} yields the conclusion for $X$ SOS-diffusion. 
		
		{\it Step 5: Ucp convergence. }  Both the discrete occupation process and its limit are non-decreasing in time, with continuous limit. Hence 	Lemma~\ref{lem_ucp_convergence_condition} yields the ucp convergence. 
		The proof is thus completed. 
		\end{proof}

        \begin{proof}[Proof of Corollary~\ref{cor:occ:time}]
	       It is enough to consider SkS-BM and the convergence for fixed time $t>0$. The extension to SOS-diffusions and to ucp convergence follow exactly as in the proof of Lemma~\ref{lem_sosdiff_occtime_approximation}.
	       Set
	       \begin{equation}
		      V_n(U):=\frac1n\sum_{i=1}^{[nt]} \indic{U}\left(X_{\frac{i-1}{n}}\right).
	       \end{equation}
    	By Lemma~\ref{lem_kernel_bound}, for every compact $K\subset\IR$,
	       \begin{equation}
		      \limsup_{n\to\infty}\Esp_x[V_n(U\cap K)] \le C_t\,m_{\params{\rho}{\beta}}(U\cap K),
	       \end{equation}
    	where $C_t$ does not depend on $U$.
        The same bound holds for $\occt{X}{U\cap K}{t}$. 
        Hence, if $U_j\uparrow U$ or $U_j\downarrow U$ inside a fixed compact set, the corresponding differences vanish in probability, uniformly in $n$. 
        A monotone class argument therefore extends the result from intervals to Borel subsets of compact sets. Finally, letting the compact sets increase to $\IR$ and using the non-explosion of $X$ gives the conclusion for every Borel subset of $\IR$.
    \end{proof}

		\section{Proofs of the preliminary results for SkS-BM} 
		\label{app_proofsSemigroup}

		In this section we prove all postponed proofs of Section~\ref{sec_preliminary}, that is, Proposition~\ref{cor_semi_group_SOS-BM_timescaling_ptk},  the kernel bounds (Lemmas~\ref{lem_kernel_bound}, \ref{prop_semig_bound}, \ref{lem:semig:bound:strong}, and \ref{lem_gamma_bounds_SOS-BM}), and Proposition~\ref{prop_mngn_convergence}. 

        Let us recall that we use Notation~\ref{notation:SkS} and that $a(x)=1+\sgn(x)\beta$. 
		
		\subsection{Proof of Proposition~\ref{cor_semi_group_SOS-BM_timescaling_ptk} (time-space scaling)}
		\label{aapp: proof_TS_scaling}
		
		The occupation time of the positive half-line $[0,\infty) $ by the process $X$ is defined as 
		\begin{equation}
			\begin{aligned}
				\occt{X}{+}{t} &= \int_{0}^{t} \indicB{X_s \ge 0} \vd s,
				& t&\ge 0.
			\end{aligned}
		\end{equation}
		
        \begin{lemma}[Scaling property]
	\label{prop_sticky_skew_scaling}
	Let $\rho\geq 0$, $\beta\in(-1,1)$, $x\in\IR$, and $c>0$.
	Let $X=X^{(\rho,\beta)}$ be the SkS-BM on $\mc P_x$, and let
	$\widetilde X=X^{(\rho/\sqrt c,\beta)}$ be the SkS-BM on
	$\widetilde{\mc P}_{x/\sqrt c}$.
	Then
	\begin{align}
		\label{eq_rescaled_process}
		&\law_{\Prob_x}
		\left(
		X_{ct},
		\loct{X}{0}{ct},
		\occt{X}{+}{ct};
		t\geq0
		\right)
		\\
		&\qquad =
		\law_{\widetilde\Prob_{x/\sqrt c}}
		\left(
		\sqrt c\,\widetilde X_t,
		\sqrt c\,\loct{\widetilde X}{0}{t},
		c\,\occt{\widetilde X}{+}{t};
		t\geq0
		\right).
	\end{align}
\end{lemma}

\begin{proof}
	Define, on $\mc P_x$,
	\[
		Y_t:=\frac{1}{\sqrt c}X_{ct},
		\qquad
		W_t:=\frac{1}{\sqrt c}B_{ct},
		\qquad t\geq0.
	\]
	With respect to the filtration
	$\widetilde{\bF}_t:=\bF_{ct}$, the process $W$ is a Brownian motion.

	We first identify the symmetric local time of $Y$ at $0$.
	Since
	\[
		\qv{Y}_t=\frac1c\,\qv{X}_{ct},
	\]
	the characterization of symmetric local time yields
	\begin{align}
		\loct{Y}{0}{t}
		&=
		\lim_{\varepsilon\downarrow0}
		\frac{1}{2\varepsilon}
		\int_0^t
		\indicB{|Y_s|<\varepsilon}\,
		\vd\qv{Y}_s
		\\
		&=
		\lim_{\varepsilon\downarrow0}
		\frac{1}{2\varepsilon c}
		\int_0^{ct}
		\indicB{|X_u|<\sqrt c\,\varepsilon}\,
		\vd\qv{X}_u
		=
		\frac1{\sqrt c}\loct{X}{0}{ct}.
		\label{eq_local_time_scaling}
	\end{align}

	Using the equation satisfied by $X$ and the definitions of $Y$ and $W$,  
	\begin{align}
		Y_t
		&=
		\frac{x}{\sqrt c}
		+
		\frac1{\sqrt c}
		\int_0^{ct}
		\indicB{X_u\neq0}\,\vd B_u
		+
		\frac{\beta}{\sqrt c}\loct{X}{0}{ct}
		\\
		&=
		\frac{x}{\sqrt c}
		+
		\int_0^t
		\indicB{Y_s\neq0}\,\vd W_s
		+
		\beta\,\loct{Y}{0}{t},
	\end{align}
	where~\eqref{eq_local_time_scaling} was used in the last equality.

	Moreover,
	\begin{align}
		\int_0^t\indicB{Y_s=0}\,\vd s
		&=
		\frac1c
		\int_0^{ct}\indicB{X_u=0}\,\vd u
        =
		\frac{\rho}{c}\loct{X}{0}{ct}
		=
		\frac{\rho}{\sqrt c}\loct{Y}{0}{t}.
	\end{align}
	Hence $Y$ is a weak solution of the SkS-BM system with parameters
	$(\rho/\sqrt c,\beta)$ and initial condition $x/\sqrt c$.

	By uniqueness in law, 
	$
		\law_{\Prob_x}(Y)
		=
		\law_{\widetilde\Prob_{x/\sqrt c}}(\widetilde X).
	$
	Finally,
	\begin{align}
		\occt{Y}{+}{t}
		&=
		\int_0^t\indicB{Y_s\geq0}\,\vd s
		=
		\frac1c
		\int_0^{ct}\indicB{X_u\geq0}\,\vd u
		=
		\frac1c\occt{X}{+}{ct}.
	\end{align}
	Combining the preceding identities proves~\eqref{eq_rescaled_process}.
\end{proof}

    \begin{remark}
	   The scaling property can also be recovered from the explicit joint laws of the position, local time, and occupation time given in 	\cite{Touhami2021,Casteras2023}. 
       The direct proof above avoids the use of these joint-law formulas and the associated parametrization conventions.
      The joint scaling statement above is stronger than what is needed for Proposition~\ref{cor_semi_group_SOS-BM_timescaling_ptk}.
    \end{remark}
		
		\begin{proof}
			[Proof of Proposition~\ref{cor_semi_group_SOS-BM_timescaling_ptk}]
			This follows immediately from the first marginal of Proposition~\ref{prop_sticky_skew_scaling} and the definition of the transition kernel.
		\end{proof}
		
		\subsection{Proof of kernel bounds}
		\label{aapp: proof_kernel_bounds}
		
		\begin{proof}
			[Proof of Lemma~\ref{lem_kernel_bound}]
			Let $u_1,u_2,v_{\rho} $ be the functions defined in \eqref{eq_def_SkSBM_kernel_factorization}--\eqref{eq_def_SkSBM_kernel_factors}.
			We first observe that 
			\begin{equation}\label{eq_proof_kernelbound_step0}
				\frac{1}{a(y)} \braces{u_1(t,x,y) - u_2(t,x,y)} 
				\le \frac{1}{a(y)} u_1(t,x,y) \leq K u_1(t,x,y),
			\end{equation}
			for some positive constant $K$ that depends on $\beta$, since $a(y):=1+\sgn(y)\beta \geq 1-|\beta|$. 
			From the Mills ratio (see \cite[p.~98]{GriSti}), $\erfc(x) \sim e^{-x^2}/x $ and thus there is a constant $K_{\Mills}>0 $ such that 
			$ \erfc(x) \le  K_{\Mills} e^{-x^2}/x$.
			Thus, if $ \rho>0 $,
			\begin{align}
				v_{\rho}(t,x,y)  ={} & 
				\frac{1}{\rho} e^{2(|x|+|y|)/\rho + 2t/\rho^2} \erfc \braces{\frac{|x|+|y|}{\sqrt{2t}} + \frac{\sqrt{2t}}{\rho}} 
				\\ \le{} & K_{\Mills} \frac{1}{\rho}  \frac{\rho \sqrt{2 t}}{\rho (|x|+|y|) + 2t} e^{-(|x|+|y|)^2/2t} 
				\le \sqrt{ \pi} K_{\Mills} u_1(t,x,y).
				\label{eq_proof_kernelbound_step1}
			\end{align}
			If $\rho=0 $, from \eqref{eq_def_SkSBM_kernel_factors},
			\begin{equation}\label{eq_proof_kernelbound_step1b}
				v_{0}(t,x,y)  = u_2(t,x,y) \le u_1(t,x,y).
			\end{equation}
			Combining~\eqref{eq_proof_kernelbound_step0}, \eqref{eq_proof_kernelbound_step1} and \eqref{eq_proof_kernelbound_step1b},
			completes the proof.
		\end{proof}
		
		\begin{proof}
			[Proof of Lemma~\ref{prop_semig_bound}]
			By definition of the probability transition kernel,  
			\begin{equation}
				P^{\params{\rho}{\beta}}_t h(x) 
				= \int_{\IR} h(y) p_{\params{\rho}{\beta}}(t,x,y) \, m_{\params{\rho}{\beta}}(\rd y).
			\end{equation}
			This, Lemma~\ref{lem_kernel_bound}, and  
			the identity $p_{\params{\rho}{\beta}}(t,x,0)= v_{\rho}(t,x,0)$, yield that
			\begin{align}\label{eq_proof_kernel_bound_step0}
				|P^{\params{\rho}{\beta}}_t h(x)| \le{}&
				\int_{\IR} |h(y)| p_{\params{\rho}{\beta}}(t,x,y) \, m_{\params{\rho}{\beta}}(\rd y)
				\\
				={} &
				\rho |h(0)|  p_{\params{\rho}{\beta}}(t,x,0) 
				+ \int_{\IR} |h(y)| p_{\params{\rho}{\beta}}(t,x,y) \, m_{\params{0}{\beta}}(\rd y)
				\\ \le{}&
				\rho K u_1(t,x,0) |h(0)|
				+
				K \xnorm{u_1(t,x,\cdot)}{\infty}  m_{\params{0}{\beta}}(|h|)
				\\ \le{}&
				\frac{K}{\sqrt{2\pi t}} \braces{\rho |h(0)| + m_{\params{0}{\beta}}(|h|)},
			\end{align}
			for some constant $K>0$ that depends on $\beta$ and not on $\rho$. The proof is thus completed.
		\end{proof}

		\begin{proof}
			[Proof of Lemma~\ref{lem:semig:bound:strong}]
			We observe that
			\begin{align}
				& \abs{P^{\params{\rho}{\beta}}h(x) 
					- m_{\params{\rho}{\beta}}(h)p_{\params{\rho}{\beta}}(t,x,0) } 
				\\
				& = \abs{
					\int_{\IR} h(y) \braces{ p_{\params{\rho}{\beta}}(t,x,y) - 
						p_{\params{\rho}{\beta}}(t,x,0) } \, m_{\params{0}{\beta}}(\rd y)
				} \\
				& \le \abs{
					\int_{\IR} h(y) \braces{ u_1(t,x,y) - 
						u_1(t,x,0) }  \vd y
				}
				+ \abs{
					\int_{\IR} h(y) \braces{ u_2(t,x,y) - 
						u_2(t,x,0) }  \vd y
				}
				\\ &
				\qquad + \abs{
					\int_{\IR} h(y) \braces{ v_{\rho}(t,x,y) - 
						v_{\rho}(t,x,0) } \, m_{\params{\cdot}{\beta}}(\rd y)
				}. \label{eq_kernel_integral_decomp_0}
			\end{align}
			For every $\gamma \ge 0 $, there exists a constant $K_{\gamma}''>0 $ such that, for every $x\in \IR $ and $t>0 $,
			\begin{equation}\label{eq_proof_semi_group_bound_2_integral_3_3}
				e^{-x^2/t} \le e^{-x^2/2t} \le K_{\gamma}''\frac{1}{1 + |x/\sqrt{t}|^{\gamma}}.
			\end{equation}
			From \cite[Lemma 3.1]{Jac98}, for every $\gamma\ge 0 $, there exists a $K_{\gamma}'>0 $ such that
			\begin{equation}\label{eq_proof_semi_group_bound_2_integral_1}
				\abs{\int_{\IR} h(y) \sqbraces{u_1(t,x,y) - u_1(t,x,0)} \vd y}  \le \frac{K_{\gamma}'}{t} 
				\braces{\lambda^{(1)}(h) + \frac{\lambda^{(1)}(h)}{1 + |x/\sqrt t|^{\gamma}} + \frac{\lambda^{(\gamma)}(h)}{1+ |x|^{\gamma}}}.
			\end{equation}
			
			For the second additive term on the right-hand-side of \eqref{eq_kernel_integral_decomp_0}, an integration by parts yields that
			\begin{equation}
				\int_{\IR} h(y) \sqbraces{u_2(t,x,y) - u_2(t,x,0)} \vd y
				= \int_{0}^{1}\int_{\IR}h(\zeta)\zeta \frac{\partial}{\partial y}u_2(t,x,y) \Big|_{y = \zeta\theta} \vd \zeta \vd \theta.
			\end{equation}
			Since there exists a constant $c>0 $ such that $|x e^{-x^2}| \le c $, for all $x\in \IR $, we have that
			\begin{equation}
				\abs{\frac{\partial}{\partial y}u_2(t,x, y)} 
				= \abs{\sgn(y)\frac1{t \sqrt{\pi}}\frac{|y|+|x|}{\sqrt{2t}} e^{-\frac{(|x|+|y|)^2}{2t}}} \le \frac{c/ \sqrt{\pi}}{t}, 
				\quad \text{for all }y\in \IR. 
			\end{equation}
			Setting $c' = c / \sqrt{\pi} $ yields 
			\begin{align}
				\Big| \int_{\IR} h(y) &\sqbraces{u_2(t,x,y) - u_2(t,x,0)} \vd y \Big| 
				\\ &\le \int_{0}^{1}\int_{\IR} \abs{h(\zeta)\zeta \frac{\partial}{\partial y}u_2(t,x,y) \Big|_{y = \zeta\theta} } \vd \zeta \vd \theta 
				\le c' \frac{\lambda^{(1)}(h)}{t}.
				\label{eq_proof_semi_group_bound_2_integral_2}
			\end{align}
			
			We now bound the third additive term at the right-hand-side of \eqref{eq_kernel_integral_decomp_0}
			and first consider the case $\rho=0 $, where for all $(t,x,y)\in(0,\infty)\times \IR^{2}$, 
			$v_{0}(t,x,y)=u_2(t,x,y)  $.
			We observe that 
			\begin{align}
				\int_{\IR_+} h(y) \sqbraces{u_2(t,x,y) - u_2(t,x,0)} \vd y&= 
				\int_{0}^{1}\int_{\IR_+}h(\zeta)\zeta \frac{\partial}{\partial y}u_2(t,x,y) \Big|_{y = \zeta\theta} \vd \zeta \vd \theta,\\
				\int_{\IR_-} h(y) \sqbraces{u_2(t,x,y) - u_2(t,x,0)} \vd y&= 
				\int_{0}^{1}\int_{\IR_-}h(\zeta)\zeta \frac{\partial}{\partial y}u_2(t,x,y) \Big|_{y = \zeta\theta} \vd \zeta \vd \theta.
			\end{align}
			Applying the same arguments as for the second additive term on each half-line yields the bound
			\begin{equation}\label{eq_proof_vrho_int_decomp_final_1b}
				\abs{\int_{\IR} h(y) \braces{ v_{0}(t,x,y) - 
						v_{0}(t,x,0) } \, m_{\params{\cdot}{\beta}}(\rd y)}
				\le c' \frac{ m_{\params{\cdot}{\beta}}^{(1)}(h)}{t}.
			\end{equation}
			
			For the case $\rho>0 $, let $M $ be the function defined for all $t\ge 0 $ and $x,y \in \IR $ by
			\[
			M(t,x,y) := \frac{v_\rho (t,x,y)}{u_2(t,x,y)}=\frac{\sqrt {\pi}}{\sqrt{ 2 t}}\frac{2 t}{\rho} f\left(\frac{|x|+|y|+ 2 t/\rho}{\sqrt{2 t}}\right),
			\]
			where $f(z):=e^{z^2} \erfc(z)$.
			This yields that
			\begin{align}
				v_\rho(t,x,y) &= u_2(t,x,y) M(t,x,y) 
				\\ &\quad = u_2(t,x,y) \braces{M(t,x,y)-M(t,x,0) +M(t,x,0)}, 
			\end{align}
			and therefore 
			\begin{align}\label{eq_vrho_M_integrant_decomp_0}
				v_\rho(t,x,y) - v_\rho(t,x,0)  
				=& u_2(t,x,y) \braces{M(t,x,y)-M(t,x,0)}
				\\ &\quad + \braces{u_2(t,x,y)-u_2(t,x,0)} M(t,x,0).
			\end{align}
			We observe that $f$ is decreasing and $|z|f(|z|)\in [0,1/\sqrt{\pi}]$. Hence, for all $t>0$ and $x,y \in \IR $, we have
			\begin{align}
				0 &\leq M(t,x,y) \leq M(t,x,0) \leq M(t,0,0) = \frac{\sqrt \pi}{\sqrt{2 t}}\frac{2 t}{\rho} f\braces{\frac{2 t/\rho}{\sqrt{2 t}}} \leq 1.
			\end{align}
			This and ~\eqref{eq_proof_vrho_int_decomp_final_1b} yield 
			\begin{equation}\label{eq_proof_semi_group_bound_2_integral_2bis}
				\abs{\int_{\IR} h(y) \sqbraces{u_2(t,x,y) - u_2(t,x,0)} |M(t,x,0)|  \, m_{\params{0}{\beta}}(\rd y)}
				\le c' \frac{ m_{\params{0}{\beta}}^{(1)}(h)}{t}.
			\end{equation}
			Moreover,
			\[0 \le M(t,x,0)-M(t,x,y) 
			\le \frac{\sqrt \pi}{ \sqrt{2 t}}\frac{2 t}{\rho} \int_{0}^{|y|/\sqrt{2 t}} \left|f'\left(\frac{|x|+ 2 t/\rho}{\sqrt{2 t}} +\zeta\right)\right| \vd \zeta.
			\]
			Since $|f'(\zeta)|= \frac2{\sqrt{\pi}} (1- \sqrt \pi \zeta f(\zeta))$ is decreasing,
			\begin{align}
				0 \le M(t,x,0)-M(t,x,y) 
				\leq &  |y| \frac{\sqrt{\pi}}{\rho}  \left|f'\left(\frac{2 t}{\rho \sqrt{2 t}}\right)\right| 
				\\ = & |y| \frac{2}{\rho} \left(1- \sqrt \pi \frac{2 t}{\rho \sqrt{2 t}} f\left(\frac{2 t}{\rho \sqrt{2 t}}\right)\right).
			\end{align}
			
			Note that the function $\sqbraces{x \mapsto x \braces{1-\sqrt{\pi} x e^{x^2} \erfc(x)} }:\IR\mapsto \IR$ takes values in $[0,1/4)$, so 
			\begin{equation}
				0 \le M(t,x,0)-M(t,x,y)
				\le 
				\frac{|y|}{2 \sqrt{2t}} .
			\end{equation}
			Therefore, we have that
			\begin{align}\label{eq_proof_semi_group_bound_2_integral_3bis}
				&\abs{\int_{\IR} h(y) \sqbraces{u_2(t,x,y)} \sqbraces{M(t,x,y)- M(t,x,0)} \, m_{\params{0}{\beta}}(\rd y)}
				\\  &\qquad\qquad \le \frac1{2 \sqrt{ 2 t}} \int_{\IR} |y| |h(y)| |u_2(t,x,y)|  \, m_{\params{0}{\beta}}(\rd y)
				\\  &\qquad\qquad \le \frac{1}{4 t \sqrt{\pi}}  m_{\params{0}{\beta}}^{(1)}(h) e^{-\frac{x^2}{2t}}
				\le \frac1t \frac{K_\gamma''}{4 \sqrt{\pi}}  m_{\params{0}{\beta}}^{(1)}(h) \frac{1}{1 + |x/\sqrt{t}|^{\gamma}}, 
			\end{align}
			where the last inequality comes from~\eqref{eq_proof_semi_group_bound_2_integral_3_3}.
			From \eqref{eq_vrho_M_integrant_decomp_0}, \eqref{eq_proof_semi_group_bound_2_integral_2bis}, and \eqref{eq_proof_semi_group_bound_2_integral_3bis},
			\begin{align}
				& \abs{
					\int_{\IR} h(y) \braces{ v_{\rho}(t,x,y)  - 
						v_{\rho}(t,x,0) }  \, m_{\params{0}{\beta}}(\rd y)
				}
				\\ & \qquad\qquad \le \frac{1}{t}\braces{c' + \frac{K_\gamma''}{4 \sqrt{\pi}}} \braces{  m_{\params{0}{\beta}}^{(1)}(h) + \frac{ m_{\params{0}{\beta}}^{(1)}(h)}{1 + |x/\sqrt{t}|^{\gamma}}}.
				\label{eq_proof_vrho_int_decomp_final_1}
			\end{align}
			From \eqref{eq_kernel_integral_decomp_0}, \eqref{eq_proof_semi_group_bound_2_integral_1}, \eqref{eq_proof_semi_group_bound_2_integral_2}, \eqref{eq_proof_vrho_int_decomp_final_1b}, \eqref{eq_proof_vrho_int_decomp_final_1}, and since that for all $\gamma \ge 0 $,  
			\begin{equation}
				(1-|\beta|) \lambda^{(\gamma)}(h)
				\le m_{\params{0}{\beta}}^{(\gamma)}(h) \le 
				(1+|\beta|) \lambda^{(\gamma)}(h),
			\end{equation}
			the desired bound holds for
			\begin{equation}
				K_{\beta,\gamma} = \frac{2}{1-|\beta|}\braces{K_{\gamma}' + c'} + \braces{c' + \frac{K_\gamma''}{4 \sqrt{\pi}}}.
			\end{equation}
			This completes the proof of Lemma~\ref{lem:semig:bound:strong}.
		\end{proof}
		
		\begin{proof}
			[Proof of Lemma~\ref{lem_gamma_bounds_SOS-BM}]
			From Corollary \ref{cor_semi_group_SOS-BM_timescaling}, which is a consequence of the scaling property,
			\begin{equation}
				\gamma^{\params{\rho}{\beta}}_n[h](x,t) = \sum_{i=2}^{[nt]} \Esp_x \sqbraces{ h(\sqrt n \hfprocess{X}{\params{\rho}{\beta}}{i-1})  } = \sum_{i=2}^{[nt]} P^{\params{\sqrt n \rho}{\beta}}_{i-1}h(\sqrt n x).
			\end{equation}
			Thus, from Lemma~\ref{prop_semig_bound} and since $\sum_{i=1}^{[nt]-1} i^{-\frac12} \leq 2 \sqrt{n t}$,
			\begin{align}\label{eq_proof_gamma_bound_alt}
				\abs{\gamma^{\params{\rho}{\beta}}_n[h](x,t)} &\le 
				\sum_{i=2}^{[nt]} \abs{ P^{\params{\sqrt n\rho}{\beta}}_{i-1}h(\sqrt n x) } 
				\\ &\le K \sum_{i=1}^{[nt]-1} \frac{1}{\sqrt{i}}m_{\params{\rho\sqrt n }{\beta}}(|h|)
				\le 	2K m_{\params{\rho \sqrt n}{\beta}}(|h|) \sqrt{n t},
			\end{align}
			for some constant $K>0$. 
            This completes the proof of the first statement.
			
			When $m_{\params{\sqrt n \rho}{\beta}}(h)=0 $,  from Lemma~\ref{lem:semig:bound:strong} (with $\gamma=1$) and the fact that 
			$\sum_{i=1}^{[nt]-1} i^{-1} \leq 1+ \log(nt)$ we obtain, for some constant $K_{\beta,1}>0$ (that depends on $\beta$ and not on $\rho$),  
			\begin{align}
				|\gamma^{\params{\rho}{\beta}}_n[h](x,t)| &\le 
				\sum_{i=2}^{[nt]} \abs{ P^{\params{\rho}{\beta}}_{i-1}h(\sqrt n x) } 
				\\ &\le 3 K_{\beta,1} m^{(1)}_{\params{0}{\beta}}(h)  \sum_{i=2}^{[nt]} \frac1i \leq  3 K_{\beta,1} m^{(1)}_{\params{0}{\beta}}(h)  (1+ \max(0,\log(n t))),
			\end{align}
			which proves the second and last statement. The proof is thus completed.
		\end{proof}
		
		\subsection{Proof of Proposition~\ref{prop_mngn_convergence}}
		\label{aap: proof_lta_limit}
		
		In the proof we use the analytic expression of the probability transition kernel of SkS-BM $(t,x,y)\mapsto p_{(\rho,\beta)}(t,x,y) $ defined in terms of the functions $u_1$, $u_2$, and $v_{\rho} $ (see \eqref{eq_def_SkSBM_kernel_factorization}--\eqref{eq_def_SkSBM_kernel_factors}). 
		
		For simplicity, let $m_n := m_{\params{\rho \sqrt n}{\beta}} $.
		Note also that $u_1 $ is the Gaussian kernel and $u_2(t,x,y) = u_1(t,|x|,-|y|) $, for all $(t,x,y) $.
		Moreover, note that for every $t>0 $, the following identities hold: 
		\begin{equation}\label{eq_exponential_int_equality_1}
			\int_{\IR_{+}}\int_{\IR_{+}} \braces{y-x} u_1(t,x,y) \vd y \vd x = \frac{1}{2}, \quad
			\int_{\IR_{+}}\int_{\IR_{+}} x\, u_2(t,x,y) \vd y \vd x = \frac{1}{4}. 
		\end{equation}	
		
		We first consider the case $\rho=0 $, so that
		\begin{align}
			m_{\params{\rho \sqrt n}{\beta}}  &= m_{\params{0}{\beta}},
			&
			\wh g_n =\wh g &: x \mapsto \Esp_{x} \braces{|X^{\params{0}{\beta}}_{1}| - |x|}.
		\end{align}
		From \eqref{eq_def_SkSBM_kernel_factorization}--\eqref{eq_def_SkSBM_kernel_factors},
		\begin{align}\label{proof_mghat_argsB10}
			m_n(\wh g_n) ={} &  
			\int_{\IR} \wh g (x)\, m_{\params{0}{\beta}}(\rd x) 
			= 
			\int_{\IR} \int_{\IR} \braces{|y| - |x|} p_{\params{0}{\beta}}(1,x,y)\,
			m_{\params{0}{\beta}}(\rd y)\, m_{\params{0}{\beta}}(\rd x)
			\\={} &  
			\int_{\IR} \int_{\IR}  \braces{|y| - |x|}
			u_1(1,x,y)
			a(x)\vd y \vd x
			-
			\int_{\IR} \int_{\IR} \braces{|y| - |x|}
			u_2(1,x,y)
			a(x) \vd y \vd x\\
			&+
			\int_{\IR} \int_{\IR} \braces{|y| - |x|} u_2 (1,x,y)
			a(x)a(y) \vd y \vd x.
		\end{align}
		By \eqref{eq_exponential_int_equality_1}, and since
		for all $t>0$ and $x,y\in \IR $,  
		$u_2(t,x,y)=u_2(t,y,x) $, 
		\begin{equation}\label{proof_mghat_argsB2}
			\begin{aligned}
				&\int_{\IR} \int_{\IR}  \braces{|y| - |x|}
				u_1 (1,x,y) 
				a(x) \vd y \vd x = 1, \\ 
				&\int_{\IR} \int_{\IR}  \braces{|y| - |x|}
				u_2 (1,x,y) a(x) \vd y \vd x 
				= 0,\\
				&\int_{\IR} \int_{\IR} \braces{|y| - |x|} u_2 (1,x,y)
				a(x)a(y) \vd y \vd x
				= 0.
			\end{aligned}
		\end{equation}
		Thus, $m_n(\wh g_n) = 1 $.
		
		We now suppose that $\rho>0 $.
		Then, 
		\begin{align}\label{proof_mghat_argsA1}
			m_n(\wh g_n) &= \int_{\IR} \wh g_n(x) m_n(\vd x) 
			\\ &= \int_{\IR} \int_{\IR} \braces{|y| - |x|} p_{\params{\sqrt n \rho}{\beta}}(1,x,y)
			m_n(\vd y) m_n(\vd x).
		\end{align}
		Since $m_n $ is defined for all $n\in \IN $ and $x\in \IR $ by $m_n(\rd x) =  a(x)\rd x + \rho \sqrt n \delta_0(\rd x) $,
		\begin{align}
			m_n(\wh g_n) ={} & \int_{\IR} \int_{\IR} \braces{|y| - |x|} p_{\params{\sqrt n \rho}{\beta}}(1,x,y)
			a(x)a(y)
			\vd y \vd x\\
			&+ \sqrt n \rho \braces{
				\int_{\IR} |y| p_{\params{\sqrt n \rho}{\beta}} (1,0,y) 
				a(y) \vd y 
				-   \int_{\IR} |x| p_{\params{\sqrt n \rho}{\beta}}(1,x,0)
				a(x) \vd x
			}. \label{proof_mghat_argsA2}
		\end{align}
		We now show that the last additive term of the right-hand-side of \eqref{proof_mghat_argsA2} vanishes. 
		Indeed, we observe that 
		$\int_{\IR} |y| p_{\params{\sqrt n \rho}{\beta}}
		\braces{t,0,y} a(y)
		\vd y < \infty$ and that for all $t>0 $ and $x,y\in \IR $,
		$p_{\params{\sqrt n \rho}{\beta}}(t,x,y)=p_{\params{\sqrt n \rho}{\beta}}(t,y,x) $. 
		Thus,
		\begin{equation}
			\int_{\IR} |y| p_{\params{\sqrt n \rho}{\beta}}
			\braces{t,0,y} a(y)
			\vd y
			=   \int_{\IR} |y| p_{\params{\sqrt n \rho}{\beta}}
			\braces{t,y,0} 
			a(y)
			\vd y
		\end{equation}
		and 
		\begin{equation}\label{proof_mghat_argsA3}
			\sqrt n \rho \braces{
				\int_{\IR} |y| p_{\params{\sqrt n \rho}{\beta}} (1,0,y) 
				a(y) \vd y 
				-   \int_{\IR} |y| p_{\params{\sqrt n \rho}{\beta}}\braces{1,y,0} 
				a(y) \vd y
			} = 0.
		\end{equation}
		
		For the first additive term of the right-hand-side of \eqref{proof_mghat_argsA2}, from \eqref{eq_def_SkSBM_kernel_factorization}, 
		\begin{align}\label{proof_mghat_argsB1}
			\int_{\IR} \int_{\IR} \braces{|y| - |x|} &p_{\params{\sqrt n \rho}{\beta}}\braces{1,x,y} 
			a(x)a(y)
			\vd y \vd x 
			\\ =&
			\int_{\IR} \int_{\IR}  \braces{|y| - |x|}
			\braces{u_1(1,x,y) - u_2(1,x,y)} 
			a(x)\vd y \vd x
			\\ &+
			\int_{\IR} \int_{\IR} \braces{|y| - |x|} v_{\sqrt n \rho }\braces{1,x,y} 
			a(x)a(y) \vd y \vd x.
		\end{align}
		We observe that $v_{\sqrt n \rho }\braces{t,x,y}$ vanishes as $n \to \infty$.
		Moreover, using $\erfc(z)\le e^{-z^2}$ for $z\ge0$, we have, for some constant $K>0$,
        \begin{equation}
	       \big||y|-|x|\big| v_{\sqrt n\rho}(t,x,y)a(x)a(y)
           \le \frac{K}{\rho}(|x|+|y|) e^{-(|x|+|y|)^2/(2t)}.
        \end{equation}
        The right-hand side of the latter inequality is integrable on $\IR^2$ and does not depend on $n$.
        Thus, by dominated convergence, 
		\begin{equation}\label{proof_mghat_argsB3}
			\lim_{n\to \infty}\int_{\IR} \int_{\IR} \braces{|y| - |x|} v_{\sqrt n \rho }\braces{1,x,y} 
			a(x)a(y) \vd y \vd x =0.
		\end{equation}
		Therefore,
		\begin{equation}
			\lim_{n\to \infty} m_n(\wh g_n)
			= 
			\int_{\IR} \int_{\IR}  \braces{|y| - |x|}
			\braces{u_1(1,x,y) - u_2(1,x,y)} 
			a(x)\vd y \vd x,
		\end{equation}
		which by~\eqref{proof_mghat_argsB2} equals $1$.
		This proves the result.
		\qed

		\begin{small}
        \setlength{\itemsep}{0pt}
\setlength{\parsep}{0pt}
\setlength{\parskip}{0pt}
\setlength{\baselineskip}{10pt}  
		\bibliography{bibfile}
		\bibliographystyle{abbrv}
        \end{small}
		
	\end{document}